\documentclass[11pt]{article}

\usepackage{amsmath}
\usepackage{amsfonts}

\newsymbol\Subset 1362

\newcommand{\be}{\begin}
\newcommand{\Ref}[1]{(\ref{#1})}

\newcommand{\cb}{\mathcal B}
\def\cc{\mathcal C}

\def\ce{\mathcal E}
\def\cf{\mathcal F}
\newcommand{\cg}{\mathcal G}
\def\ci{\mathcal I}

\def\cO{\mathcal O}

\def\fl{\mathcal R}
\def\xcs{\mathcal S}

\def\cv{\mathcal V}
\def\cw{\mathcal W}

\newcommand{\al}{\alpha}
\newcommand{\del}{\delta}
\def\udel{\underline\delta}

\newcommand{\eps}{\epsilon}
\def\ga{\gamma}
\def\Ga{\Gamma}
\def\ka{\kappa}
\def\lla{\lambda}
\def\La{\Lambda}
\def\om{\omega}
\def\Om{\Omega}
\newcommand{\si}{\sigma}

\def\ff{\mathbf f}
\def\FF{\mathbf F}
\def\Bg{\mathbf g}
\def\qq{\mathbf q}

\newcommand{\co}{\mathbb C}

\def\bs{\mathbb S}
\def\bt{\mathbb T}

\newcommand{\R}{\mathbb R}

\newcommand{\z}{\mathbb Z}

\def\orp{\oline\R_+}

\newcommand{\oset}{\overset}
\newcommand{\uline}{\underline}
\newcommand{\oline}{\overline}

\newcommand{\la}{\langle}
\newcommand{\ra}{\rangle}
\newcommand{\st}{\,:\,}
\newcommand{\ti}{\tilde}

\newcommand{\prtl}{\partial}
\def\square{\kern20pt{\vbox{\hrule height.4pt
        \hbox{\vrule width.4pt height 6pt\kern6pt
                \vrule width.4pt}
        \hrule height.4pt}}}
\newcommand{\ry}{\R\times Y}
\newcommand{\rpy}{\R_+\times Y}
\def\rmy{\R_-\times Y}

\def\im{\text{im}}
\def\Im{\text{Im}}

\def\coker{\text{coker}}

\def\inv{^{-1}}
\def\hol{\text{Hol}}

\def\const{\text{const}\cdot}
\def\konst{\text{const}}
\def\rf{{\R}^4}
\def\proof{{\em Proof.}\ }
\def\loc{{\text{loc}}}

\def\U#1{\text{U}(#1)}

\def\SO#1{\text{SO}(#1)}

\def\spc{\text{spin}^c}

\def\Spinc#1{\text{Spin}^c(#1)}

\def\ind{\text{ind}}
\def\index{\text{index}}
\def\itr{\text{int}}

\def\hx{{\mathfrak b}}
\def\prt{\mathfrak p}

\def\fq{\mathfrak q}
\def\mm{\mathfrak m}
\def\fo{\mathfrak o}
\def\ft{\mathfrak t}

\def\cd{{\vartheta}}
\def\sa{\mathsf a}
\def\sB{\mathsf B}
\def\sC{\mathsf C}
\def\sd{\mathsf d}

\def\ssf{\mathsf f}
\def\sr{\mathsf r}
\def\sS{\mathsf S}
\def\sv{\mathsf v}
\def\sw{\mathsf w}
\def\sW{\mathsf W}

\def\llw#1#2#3{{L^{#1,#2}_#3}}
\def\lw#1#2{{L^{#1,#2}}}

\def\id{\text{id}}
\def\endt#1{_{:#1}}
\def\Endt#1{^{:#1}}

\def\tu{^{(T)}}
\def\xt{{X\tu}}
\def\mt{M\tu}

\def\Ht{H\tu}
\def\ur{\U1^{r_0}}

\def\bde{\mathbf E}
\def\bdf{\mathbf F}
\def\bdi{\mathbf I}
\def\bds{\mathbf S}
\def\rr{_o}
\def\hxi{\hat\Xi}
\def\hzeta{\hat\zeta}
\def\ual{{\uline\al}}

\def\rp#1{^{(#1)}}

\def\hatf{\hat F}
\def\tmin{\check T}
\def\pso{P_{\text{SO}}}
\def\pglp{P_{\text{GL}^+}}
\def\pglc{P_{\text{GL}^c}}
\def\pspc{P_{\text{Spin}^c}}
\def\hatf{\hat F}

\def\xtn{X^{(T(n))}}
\def\phol{\theta}
\def\enm#1{\|#1\|_{\ce'}}
\def\lnm#1{\|#1\|_{\llw p\ka1}}
\def\gl{^\#}
\def\xx#1{C_{#1}}
\newcommand{\wm}{\hat M}
\newcommand{\tcb}{\breve\cb}
\newcommand{\tcg}{\breve\cg}
\newcommand{\E}{\cb^{**}_\hx(K)}
\def\supp{\text{supp}}
\def\dds{(\del_S)_{\ell,m}}
\def\ddbs{(\del_{\mu,S})_{\ell,b+m}}
\def\tsc{\ti\sC}

\def\obeta{\oline\beta}
\def\ubeta{\uline\beta}
\def\map{\text{Map}}
\def\dth{\dot\Theta}

\begin{document}

\title{Monopoles over $4$--manifolds containing long necks, II}

\author{Kim A.\ Fr\o yshov}

\date{19 April 2006}

\maketitle

\be{abstract}
We establish a gluing theorem for monopoles over 4--manifolds containing
long necks. The theorem is stated in terms of an ungluing map defined
explicitly in terms of data that appear naturally in
applications.
Orientations of moduli spaces are handled using
Benevieri--Furi's concept of orientations of Fredholm operators
of index 0.
\end{abstract}

\bibliographystyle{plain}

\tableofcontents

\section{Introduction}

In this paper we prove a gluing theorem for monopoles suitable for the 
construction of Floer homology groups in the simplest cases and
for establishing certain gluing formulae for Seiberg--Witten
invariants of $4$--manifolds (to be discussed in \cite{Fr12,Fr4}).

There is now a large literature on gluing theory for instantons and monopoles.
The theory was introduced by Taubes~\cite{T4,T5},
who used it to obtain existence
results for self-dual connections over closed $4$--manifolds. 
It was further developed in seminal work of
Donaldson~\cite{D1}, see also Freed--Uhlenbeck~\cite{FU}. 
General gluing theorems for instantons over connected sums were proved by
Donaldson~\cite{D7} and Donaldson--Kronheimer~\cite{DK}.
In the setting of instanton Floer theory there is a highly readable account
in Donaldson~\cite{D5}, see also Floer~\cite{F1} and Fukaya~\cite{Fukaya1}.
Gluing with degenerate asymptotic limits was studied by
Morgan--Mrowka~\cite{Morgan-Mrowka1}; part of their work was adapted to
the context of monopoles by
Safari~\cite{Safari1}. Nicolaescu~\cite{nico1}
established gluing theorems for monopoles in certain situations, including
one involving gluing obstructions. Marcolli--Wang~\cite{Marcolli-Wang1}
discuss gluing theory in connection with monopole Floer homology.
For monopoles over closed $3$--manifolds split along certain tori,
see Chen~\cite{wchen1}.
Finally, gluing theory is a key ingredient in a large programme of 
Feehan--Leness~\cite{Feehan-Leness2} (using ideas of Pidstrigach--Tyurin)
for proving Witten's conjecture
relating Donaldson and Seiberg--Witten invariants.

As should be evident from this brief survey,
there are many different hypotheses under which one can
consider the gluing problem. This paper does not aim at the utmost 
generality, but is an expositary account of gluing in what might be called
the favourable cases. More precisely, we will glue precompact families of
regular monopoles over $4$--manifolds with tubular ends, under similar
general assumptions as in \cite{Fr10}.
Although obstructed gluing is not
discussed explicitly, we will show in \cite{Fr12} how the
parametrized version of our gluing theorem can be used to handle
one kind of gluing obstructions.

One source of difficulty when formulating a gluing theorem is that gluing 
maps are in general not canonical, but rather depend on various choices
hidden in their construction.
We have therefore chosen to express our gluing
theorem as a statement about an {\em ungluing map}, which is explicitly defined
in terms of data that appear naturally in applications.

If $X$ is a $4$--manifold with tubular ends and $\xt$ the glued manifold
as in \cite{Fr10}, then the first component of the ungluing map involves
restricting monopoles over $\xt$ to some fixed compact subset $K\subset X$
(which may also be regarded as a subset of $\xt$ when each $T_j$
is large). In the case of gluing along a reducible critical point, the
ungluing map has an additional component which reads off the $\U1$
gluing parameter by measuring the holonomy along a path running once
through the corresponding neck in $\xt$. 

Ungluing maps of a different kind
were studied already in \cite{D1,FU}, but later authors have mostly formulated
gluing theorems in terms of gluing maps, usually without characterizing these
maps uniquely.

The proof of the gluing theorem is divided into two parts: surjectivity and
injectivity of the ungluing map. In the first part the (quantitative)
inverse function theorem is used to construct a smooth local right inverse
$\hzeta$ of an ``extended monopole map'' $\hxi$. In the second part
the inverse function theorem is applied a second time to show, essentially,
that the image of $\hzeta$ is not too small. There are many similarities
with the proof of the gluing theorem in \cite{DK}, but also some
differences. For instance, we do not use the method of continuity, and
we handle gluing parameters differently.

It may be worth mentioning that the proof does not depend on unique
continuation for monopoles (only for harmonic spinors), as we do not
know whether solutions to our perturbed monopole equations satisfy any
such property. (Unique continuation for genuine monopoles was used in
\cite[Proposition~4.3]{Fr10} in the discussion of perturbations, but this
has little to do with gluing theory.) Therefore, in the injectivity part of
the proof, we argue by contradiction, restricting monopoles to ever larger
subsets $\ti K\subset X$. This is also reflected in the statement of the
theorem, which would have been somewhat simpler if unique continuation
were available.

We also give a detailed account of orientations of moduli spaces, using
Benevieri--Furi's concept of orientations of Fredholm operators
of index~$0$ \cite{Ben-Furi}.
This seems simpler to us than the standard approach using
determinant line bundles. Our main result here,
Theorem~\ref{thm:uglmap-ortn}, tells when the ungluing map preserves resp.\ 
reverses orientation. The length of this part of the paper is much due to the
fact that we allow gluing along reducible critical points and that we work with
(multi)framed moduli spaces (as a means of handling reducibles over the
$4$--manifolds). Although the orientability of monopole moduli spaces
over closed $4$--manifolds is well documented (see \cite{Morgan1,das2}) and
the behaviour of the determinant line bundle under gluing along irreducibles
is described in \cite{D5}, beyond this the existing literature seems short
on details concerning the orientation issue in the type of gluing
problem considered in this paper.

The author gratefully acknowledges the hospitality of the
Institut des Hautes \'Etudes Scientifiques, where
a large part of this work was carried out.


\section{The gluing theorem}
\label{sec:gluing-thm}

\subsection{Statement of theorem}
\label{subsec:gluing-statement1}

Consider the situation of \cite[Subsection~1.4]{Fr10}, but without
assuming any of the conditions (B1),(B2),(C). We now assume that 
every component of $X$ contains an end $\R_+\times Y_j$ or $\R_+\times(-Y_j)$
(ie an end that is being glued).
Fix non-degenerate monopoles $\al_j$ over $Y_j$ and $\al'_j$ over $Y'_j$.
(These should be smooth configurations rather than gauge equivalence classes
of such.)
Suppose $\al_j$ is reducible for $1\le j\le r_0$ and irreducible for
$r_0<j\le r$, where $0\le r_0\le r$. We consider monopoles over $X$ and $\xt$
that are asymptotic to $\al'_j$ over $\rpy'_j$ and (in the case of $X$)
asymptotic to $\al_j$ over $\R_+\times(\pm Y_j)$. These monopoles build
moduli spaces
\[M_\hx=M_\hx(X;\vec\al,\vec\al,\vec\al'),\qquad
\mt_\hx=M_\hx(\xt;\vec\al').\]
Here $\hx\subset X$ is a finite subset to be specified in a moment, and
the subscript indicates that we only divide out by those
gauge transformations that restrict to the identity on $\hx$, see
\cite[Subsection~3.4]{Fr10}. The ungluing map $\ff$ will be a diffeomorphism
between certain open subsets of $M\tu_\hx$ and $M_\hx$ when 
\[\tmin:=\min(T_1,\dots,T_r)\]
is large.

When gluing along the critical point $\al_j$, the stabilizer of $\al_j$
in $\cg_{Y_j}$ appears as a ``gluing parameter''. 
This stabilizer is a copy of $\U1$ if $\al_j$ is reducible
and trivial otherwise. When $\al_j$ is reducible
we will read off the gluing parameter by means of the holonomy
of the connection part of the glued monopole along a path $\ga_j$
in $\xt$ which runs once through the neck $[-T_j,T_j]\times Y_j$.
To make this precise, for $1\le j\le r_0$ fix $y_j\in Y_j$ and smooth paths
\[\ga_j^\pm:[-1,\infty)\to X\]
such that $\ga_j^\pm(t)=\iota^\pm_j(t,y_j)$ for $t\ge0$
and $\ga_j^\pm([-1,0])\subset X\endt0$.
Let $\hx$ denote the collection
of all the start-points $o^\pm_j:=\ga^\pm_j(-1)$. (We do not assume that
these are distinct.)
Note in passing that we then have
\[M^*_\hx=M_\hx.\]
Define the smooth path
\[\ga_j:I_j=[-T_j-1,T_j+1]\to\xt\]
by
\[\ga_j(t)=\be{cases}
\pi_T\ga_j^+(T_j+t), & -T_j-1\le t<T_j,\\
\pi_T\ga_j^-(T_j-t), & -T_j<t\le T_j+1,
\end{cases}\]
where $\pi_T:X^{\{T\}}\to\xt$ is as in \cite[Subsection~1.4]{Fr10}.

Choose a reference configuration $S\rr=(A\rr,\Phi\rr)$ over $X$ with limit
$\al_j$ over $\R_+\times(\pm Y_j)$ and $\al_j'$ over $\R\times Y'_j$.
Let $S'\rr=(A\rr',\Phi\rr')$ denote the reference configuration over $\xt$
obtained from $S\rr$ in the
obvious way when gluing the ends. Precisely speaking,
$S\rr'$ is the unique smooth configuration over $\xt$
which agrees with $S\rr$ over
$\itr(X\endt T)$ (which can also be regarded as a subset of $X$).


If $P\to\xt$ temporarily denotes the principal $\Spinc4$--bundle defining the
$\spc$ structure, then the holonomy of a $\spc$ connection $A$ in $P$
along $\ga_j$ is a $\Spinc4$--equivariant map
\[\text{hol}_{\ga_j}(A):P_{o^+_j}\to P_{o^-_j}.\]
Because $A$ and $A'\rr$ map to the same connection in the tangent bundle
of $\xt$, there is a unique element $\hol_j(A)$ in $\U1$ (identified with
the kernel of $\Spinc4\to\SO4$) such that
\[\text{hol}_{\ga_j}(A)=\hol_j(A)\cdot \text{hol}_{\ga_j}(A'\rr).\]
Explicitly,
\be{equation}\label{hol-defn}
\hol_j(A)=\exp\left(-\int_{I_j}\ga_j^*(A-A'\rr)\right),
\end{equation}
where as usual $A-A'\rr$ is regarded as an imaginary valued $1$--form on
$\xt$. For gauge transformations $u:\xt\to\U1$ we have
\be{equation}\label{eqn:u-hol}
\hol_j(u(A))=u(o^-_j)\cdot\hol_j(A)\cdot u(o^+_j)\inv.
\end{equation}
In particular, there is a natural smooth map
\[\text{Hol}:\mt_\hx\to\ur,\quad[A,\Phi]\mapsto
(\text{Hol}_1(A),\dots,\text{Hol}_{r_0}(A))\] 
which is equivariant with respect to the appropriate action of
\[\bt:=\map(\hx,\U1)\approx\U1^b,\]
where $b=|\hx|$.

Consider for the moment an arbitrary compact codimension~$0$ submanifold
$K\subset X$ containing $\hx$. Let $D\tu$ be the subgroup of $H^1(\xt;\z)$
consisting of those classes whose restriction to each $Y'_j$ is zero.
Let $D_K$ be the cokernel of the restriction map $D\tu\to H^1(K;\z)$.
Here $\tmin$ should be so large that
$K$ may be regarded as a subset of $\xt$,
and $D_K$ is then obviously independent of $T$. In the following
we use the $L^p_1$ configuration spaces etc introduced in
\cite[Subsection~2.5]{Fr10}. Let $\tcg_\hx(K)$ be the kernel of
the (surjective) group homomorphism
\[\cg(K)\to\bt\times D_K,\quad u\mapsto(u|_\hx,[u]),\]
where $[u]$ denotes the image in $D_K$ of the homotopy class of $u$
regarded as an element of $H^1(K;\z)$. Set
\[\tcb_\hx(K)=\cc(K)/\tcg_\hx(K),
\qquad\tcb^*_\hx(K)=\cc^*_\hx(K)/\tcg_\hx(K).\]
On both these spaces there is a natural action of
$\bt\times D_K$. Note that $D_K$ acts freely and
properly discontinuously on the (Hausdorff) Banach manifold $\tcb^*_\hx(K)$
with quotient $\cb^*_\hx(K)$.

It is convenient here to agree once and for all that the Sobolev exponent
$p>4$ is to be an even integer. This ensures that our configuration spaces
admit smooth partitions of unity, which is needed in Subsections~\ref{subsec:model-appl}
and \ref{subsec:q-existence} (but not in the proof of Theorem~\ref{gluing-thm1}).



Fix a $\bt$--invariant open subset $G\subset M_\hx$ whose closure $\oline G$ 
is compact and contains only regular points. (Of course, $G$ is the
pre-image of an open set $G'$ in $M$, but $G'$ may not be a smooth manifold
due to reducibles and we therefore prefer to work with $G$.)

\be{defn}\label{def:kv}
By a {\em kv-pair} we mean a pair $(K,V)$ where
\be{itemize}
\item $K\subset X$ is a compact codimension~$0$ submanifold which contains
$\hx$ and intersects every component of $X$,
\item $V\subset\tcb_\hx(K)$ is a $\bt$--invariant open subset
containing $R_K(\oline G)$, where $R_K$ denotes restriction to $K$.
\end{itemize}
\end{defn}

We define a partial ordering
$\le$ on the set of all kv-pairs, by decreeing that
\[(K',V')\le(K,V)\]
if and only if $K\subset K'$ and
$R_K(V')\subset V$.  

Now fix a kv-pair $(K,V)$ which satisfies the following two additional
assumptions: firstly, that $V\subset\tcb^*_\hx(K)$; secondly, that
if $X_e$ is any component of $X$ which contains a point from $\hx$
then $X_e\cap K$ is connected.
The second condition ensures that the image of
$R_K:M_\hx=M^*_\hx\to\tcb_\hx(K)$ lies in $\tcb^*_\hx(K)$.

Suppose we are given a $\bt$--equivariant smooth map
\be{equation}\label{eqn:q}
q:V\to M_\hx
\end{equation}
such that $q(\om|_K)=\om$ for all $\om\in\oline G$.
(If $\bt$ acts freely on $\oline G$ then such a map always exists when
$K$ is sufficiently large, see Subsection~\ref{subsec:q-existence}.
In concrete applications there is often a natural
choice of $q$, see Subsections~\ref{subsec:model-appl}, \ref{subsec:dd}.)

Let $X\gl$ and the forms $\ti\eta_j,\ti\eta'_j$ be
as in \cite[Subsection~1.4]{Fr10}, and choose $\lla_j,\lla'_j>0$.

\be{thm}\label{gluing-thm1}
Suppose there is class in $H^2(X\gl)$ whose 
restrictions to $Y_j$ and $Y'_j$ are $\lla_j\ti\eta_j$ and $\lla'_j\ti\eta'_j$,
respectively, and suppose the perturbation parameters $\vec\prt,\vec\prt'$
are admissible for $\vec\al'$. Then there exists a kv-pair
$(\ti K,\ti V)\le(K,V)$ such that if 
$(K',V')$ is any kv-pair $\le(\ti K,\ti V)$ then the following holds when
$\tmin$ is sufficiently large. Set
\be{gather*}
\Ht:=\left\{\om\in\mt_\hx\st\om|_{K'}\in V'\right\},\\
\qq:\Ht\to M_\hx,\quad\om\mapsto q(\om|_K).
\end{gather*}
Then $\qq\inv G$ consists only of regular
monopoles (hence is a smooth manifold), and the $\bt$--equivariant map
$\ff:=\qq\times\text{Hol}$ restricts to a diffeomorphism
\[\qq\inv G\to G\times\ur.\]
\end{thm}

{\em Remarks.} 1. When $\tmin$ is large then $K'\subset X$ can also be regarded
as a subset of $\xt$, in which case the expression $\om|_{K'}$ in the
definition of $\Ht$ makes sense.

2. Except for the equivariance of $\ff$, the theorem remains true if 
one leaves out all assumptions on $\bt$--invariance resp.\ --equivariance
on $G$ and $q$, and on $V$ in Definition~\ref{def:kv}, above.
However, it is hard to imagine any application that would
not require equivariance of $\ff$.

3. The theorem remains true
if one replaces $\tcb_\hx(K)$ and $\tcb^*_\hx(K)$ by $\cb_\hx(K)$ and 
$\cb^*_\hx(K)$ above. However, working with $\tcb$
gives more flexibility in the construction of maps $q$, see
Subsection~\ref{subsec:dd}.

4. Concerning admissibility of perturbation parameters, see the remarks
after \cite[Theorem~1.4]{Fr10}. Note that the assumption on 
$\lla_j\ti\eta_j$ and $\lla'_j\ti\eta'_j$ in the theorem above is weaker
than either of the conditions (B1) and (B2) in \cite{Fr10}.
However, in practice the gluing theorem is only useful in conjunction with
a compactness theorem, so one may still have to assume (B1) or (B2).



The proof of Theorem~\ref{gluing-thm2} has two parts.
The first part consists in showing that $\ff$ has a smooth local right inverse
around every point in $\oline G\times\ur$ (Proposition~\ref{prop:surj} below).
In the second part we will prove that $\ff$ is injective on $\qq\inv\oline G$.
(Proposition~\ref{prop:inj} below).

\subsection{Surjectivity}\label{subsec:surj}

The next two subsections are devoted to the proof of Theorem~\ref{gluing-thm1}.
Both parts of the proof make use of the same set-up, which we now
introduce.

We first choose weight functions for
our Sobolev spaces over $X$ and $\xt$. Let $\si_j,\si_j'\ge0$ be small
constants and $w:X\to\R$ a smooth function
which is equal to $\si_jt$ on $\R_+\times(\pm Y_j)$ and equal to
$\si'_jt$ on $\R_+\times Y'_j$. As usual, we require $\si_j>0$ if $\al_j$
is reducible (ie for $j=1,\dots,r_0$), and similarly for $\si_j'$.
For $j=1,\dots,r$ choose a smooth function $w_j:\R\to\R$
such that $w_j(t)=-\si_j|t|$ for $|t|\ge1$. We will always assume $\tmin\ge4$,
in which case we can define a weight function
$\ka:\xt\to\R$ by
\be{align*}
&\ka=w\quad\text{on $\xt\setminus\cup_j\:[-T_j,T_j]\times Y_j$},\\
&\ka(t,y)=\si_jT_j+w_j(t)\quad
\text{for $(t,y)\in[-T_j,T_j]\times Y_j$}.
\end{align*}

Let $\cc$ denote the $\llw pw1$ configuration space over $X$ defined by
the reference configuration $S\rr$, and let $\cc'$ denote
the $\llw p\ka1$ configuration space over $\xt$ defined by $S'\rr$.
Let $\cg_\hx,\cg'_\hx$ be the corresponding groups of gauge transformations
and $\cb_\hx,\cb'_\hx$ the corresponding orbit spaces.


Now fix $(\om_0,z)\in\oline G\times\ur$.
Our immediate goal is to construct
a smooth local right inverse of $\ff$ around this point, but
the following set-up will also be used
in the injectivity part of the proof.

Choose a smooth representative $S_0\in\cc$ for $\om_0$
which is in temporal gauge over the ends of $X$.
(This assumption is made in order to ensure exponential decay of $S_0$.)
Set $d=\dim\,M_\hx$ and let 
$\pi:\cc\to\cb_\hx$ be the projection.
By the local slice theorem we can find a smooth map
\[\bds:\R^d\to\cc\]
such that $\bds(0)=S_0$ and such that
$\varpi=\pi\circ\bds$ is a diffeomorphism onto
an open subset of $M_\hx$.

We will require one more property of $\bds$, involving holonomy.
If $a\in\llw pw1(X;i\R)$ then we define
$\hol^\pm_j(A\rr+a)\in\U1$ by
\[\hol^\pm_j(A\rr+a)=
\exp\left(-\int_{[-1,\infty)}(\ga^\pm_j)^*a\right).\]
The integral exists because, by the Sobolev embedding 
$L^p_1\subset C^0_B$ in $\rf$ for $p>4$, we have
\be{equation}\label{ewa}
\|e^wa\|_\infty\le C\|e^wa\|_{L^p_1}=C\|a\|_{\llw pw1}
\end{equation}
for some constant $C$.
It is clear that $\hol^\pm_j$ is a smooth function on $\cc$.
Because any smooth map $\R^d\to\U1$ factors through $\exp:\R i\to\U1$,
we can arrange, after perhaps modifying $\bds$ by a smooth
family of gauge transformations that are all equal to $1$ outside 
the ends $\R_+\times Y_j$ and constant on $[1,\infty)\times Y_j$, that
\be{equation}\label{z-cond}
\hol^+_j(\bds(v))\cdot(\hol^-_j(\bds(v)))\inv=z_j
\end{equation}
for $j=1,\dots,r_0$ and every $v\in\R^d$.
Here $\hol^\pm_j(\bds(v))$ denotes the holonomy, as defined above, of the
connection part of the configuration $\bds(v)$, and the $z_j$ are the 
coordinates of $z$.

\be{lemma}\label{lemma:efg}
Let $E,F,G$ be Banach spaces, $S:E\to F$ a bounded operator and
$T:E\to G$ a surjective bounded operator such that
\[S+T:E\to F\oplus G,\quad x\mapsto(Sx,Tx)\]
is Fredholm. Then $T$ has a bounded right inverse.
\end{lemma}
\proof Because $S+T$ is Fredholm there is a bounded operator
$A:F\oplus G\to E$ such that $(S+T)A-I$ is compact. Set
$A(x,y)=A_1x+A_2y$ for $(x,y)\in F\oplus G$. Then
\[TA_2-I:G\to G\]
is compact, hence $TA_2$ is Fredholm of index~$0$. Using the surjectivity of
$T$ and the fact that any closed subspace of finite dimension or
codimension in a Banach space is complemented, it is easy to see that there
is a bounded operator $K:G\to E$ (with finite-dimensional image)
such that $T(A_2+K)$ is an isomorphism.\square 

Let
\[\Theta:\cc\to\lw pw\]
we be Seiberg--Witten map over $X$.
By assumption, every point in $\oline G$ is regular, so in particular
$\om_0$ is regular, which means that
$D\Theta(S_0):\llw pw1\to\lw pw$ is surjective.
Let $\Phi$ be the spinor part of $S_0$ and define $\ci_\Phi$ as in
\cite[Subsection~2.3]{Fr10}. Then
\[\ci_\Phi^*+D\Theta(S_0):\llw pw1\to\lw pw\]
is Fredholm, so by Lemma~\ref{lemma:efg}
$D\Theta(S_0)$ has a bounded right inverse $Q$.
(This can also be deduced from \cite[Proposition~2.2~(ii)]{Fr10}.)

Let $\sr:\R\to\R$ be a smooth function such that $\sr(t)=1$ for $t\le0$
and $\sr(t)=0$ for $t\ge1$.
For $\tau\ge1$ set $\sr_\tau(t)=\sr(t-\tau)$ and let 
$S_{v,\tau}$ be the configuration over $X$ which agrees with
$\bds(v)$ away from the ends $\R_+\times(\pm Y_j)$ and satisfies
\[S_{v,\tau}=(1-\sr_\tau)\ual_j+\sr_\tau\bds(v)\]
over $\R_+\times(\pm Y_j)$.
Here $\ual_j$ denotes, as before, the translationary invariant monopole
over $\ry_j$ determined by $\al_j$. For each $v$ we have
\[\|S_{v,\tau}-\bds(v)\|_{\llw pw1}\to0\qquad\text{as $\tau\to\infty$}.\]
Therefore, when $\tau$ is sufficiently large, the operator
\[D\Theta(S_{0,\tau})\circ Q:\lw pw\to\lw pw\]
will be invertible, and we set
\[Q_\tau=Q(D\Theta(S_{0,\tau})\circ Q)\inv:\lw pw\to\llw pw1,\]
which is then a right inverse of $D\Theta(S_{0,\tau})$. It is clear that
the operator norm $\|Q_\tau-Q\|\to0$ as $\tau\to\infty$.

{\em For the remainder of the proof of Theorem~\ref{gluing-thm1},
the term `constant' will always refer
to a quantity that is independent of $\tau,T$, unless otherwise indicated.
The symbols $C_1,C_2,\dots$ and $c_1,c_2,\dots$
will each denote at most one constant, while
other symbols may denote different constants in different contexts.}

Consider the configuration space
\[\cc_j=\ual_j+\llw p{-w_j}1\]
over $\ry_j$ and the Seiberg--Witten map
\[\Theta_j:\cc_j\to\lw p{-w_j}.\]
As explained in \cite[Subsection~3.4]{Fr10} there is an identification
\[\ci^*_{\ual}+D\Theta_j(\ual_j)=\frac d{dt}+P_\al.\]
By the results of \cite{D5} the operator on the right hand side defines a
Fredholm operator $\llw p{w_j}1\to \lw p{w_j}$, and this must be surjective
because of the choice of weight function $w_j$. In particular,
\be{equation}\label{eqn:DThj}
D\Theta_j(\ual_j):\llw p{w_j}1\to \lw p{w_j}
\end{equation}
is surjective, hence has a bounded right inverse $P_j$ by
Lemma~\ref{lemma:efg}.
(Here one cannot appeal to \cite[Proposition~2.2~(ii)]{Fr10}.)

Let $\Theta':\cc'\to\lw p\ka$ be the Seiberg--Witten map over $\xt$.
When $\tmin>\tau+1$ then by splicing 
$S_{v,\tau}$ in the natural way one obtains a smooth configuration
$S_{v,\tau,T}$ over $\xt$. There is a constant
$\xx0\gg0$ such that if
\be{equation}\label{T-bound}
\tmin>\tau+\xx0
\end{equation}
then we can splice the right inverses $Q_\tau$ and $P_1,\dots,P_r$ to
obtain a right
inverse $Q_{\tau,T}$ of\[D\Theta'(S_{0,\tau,T}):\llw p\ka1\to\lw p\ka\]
which satisfies
\[\|Q_{\tau,T}\|\le C\left(\|Q_\tau\|+\sum_j\|P_j\|\right)\]
for some constant $C$, see the appendix.
Since $\|Q_\tau\|$ is bounded in $\tau$ (ie as a function of $\tau$),
we see that $\|Q_{\tau,T}\|$ is bounded in $\tau,T$.

{\em The inequality \Ref{T-bound} will be assumed from now on.}

We now introduce certain $1$--forms that will be added to the configurations
$S_{v,\tau,T}$ in order to make small changes to
the holonomies $\hol_j$. For any 
$c=(c_1,\dots,c_{r_0})\in\R^{r_0}$ define the $1$--form
$\phol_{c,\tau}$ over $\xt$ by
\[\phol_{c,\tau}=\be{cases}
0 & \text{outside $\bigcup_{j=1}^{r_0}[-T_j,T_j]\times Y_j$},\\
ic_j\sr'_{\tau+1-T_j}dt & \text{on $[-T_j,T_j]\times Y_j,
       \quad j=1,\dots,r_0$.}
\end{cases}\]
where $\sr'_s(t)=\frac d{dt}\sr_s(t)$. Set
\[E=\R^d\times\R^{r_0}\times\lw p\ka(\xt;i\La^+\oplus\bs^-).\]

For $0<\eps<1$ let $B_\eps\subset E$
be the open $\eps$--ball about $0$. Define a smooth map $\zeta:E\to\cc'$ by
\be{equation}\label{zeta-def}
\zeta(v,c,\xi)=S_{v,\tau,T}+\phol_{c,\tau}+Q_{\tau,T}\xi,
\end{equation}
where $\phol_{c,\tau}$ is added to the connection part of
$S_{v,\tau,T}$.

When deciding where to add the perturbation $1$--form $\phol_{c,\tau}$ one has to balance
two concerns. One the one hand, because the weight function $\ka$ increases
exponentially as one approaches the middle of the necks $[-T_j,T_j]\times Y_j$,
$j=1,\dots,r_0$,
it is desirable to add $\phol_{c,\tau}$ as close to the boundaries of these necks as
possible. On the other hand, in order for Lemma~\ref{Xi-bounds1} below to work,
the spinor field of $S_{v,\tau,T}$ needs to be ``small'' in the perturbation
region. We have chosen to add $\phol_{c,\tau}$ at the negative end of the
cut-off region, where the spinor field is zero.

Although we will sometimes use the notation $\zeta(x)$, we shall think
of $\zeta$ as a function of three variables $v,c,\xi$, and $D_j\zeta$
will denote the derivative of $\zeta$ with respect to the $j$'th variable.
Similarly for other functions on (subsets of) $E$ that we will define later.
Set
\[\si=\max(\si_1,\dots,\si_r).\]
Notice that if $r_0=0$, ie if we are not gluing along any reducible 
critical point, then we may take $\si=0$.

\be{lemma}\label{zeta-bounds}
There exists a constant $\xx1>0$ such that for $x\in E$ the following hold:
\be{description}
\item[(i)]$\|D_1\zeta(x)\|,\|D^2\zeta(x)\|<\xx1$ if $\|x\|<1$,
\item[(ii)]$\|D_2\zeta(x)\|<\xx1e^{\si\tau}$,
\item[(iii)]$\|D_3\zeta(x)\|<\xx1$.
\end{description}
\end{lemma}

\proof To prove (ii), note that if $r_0>0$ and $c=(c_1,\dots,c_{r_0})$ then
\[\left\|\frac{\prtl\zeta(v,c,\xi)}{\prtl c_j}\right\|_{\llw p\ka1}=
\const e^{\si_j\tau}.\]
The other two statements are left to the reader.\square

Let $\cc'_1$ be the set of all $S\in\cc'$ such that 
$[S|_K]\in V$, $q(S|_K)\in\varpi(\R^d)$, and $\hol_j(S)\neq-z_j$
for $j=1,\dots,r_0$. Then $\cc'_1$ is an open subset of $\cc'$, and
there are unique smooth functions
\[\eta_j:\cc'_1\to(-\pi,\pi)\]
such that $\hol_j(S)=z_j\exp(i\eta_j(S))$.
Set $\eta=(\eta_1,\dots,\eta_{r_0})$ and define
\[\hat\Xi=(\varpi\inv\circ q\circ R_K,\eta,\Theta'):\cc'_1\to E.\]

A crucial point in the proof of Theorem~\ref{gluing-thm1} will be
the construction
of a smooth local right inverse of $\hat\Xi$, defined in a neighbourhood of
$0$. The map $\zeta$ is a first approximation to such a local right inverse.
The construction of a genuine local right inverse will involve an
application of the quantitative inverse function theorem (see
Lemma~\ref{Xi-inv} below).

From now on we will take $\tau$ so large that $K\subset X\endt \tau$
and
\[\hol^+_j(S_{0,\tau})\cdot(\hol^-_j(S_{0,\tau}))\inv\neq-z_j\]
for $j=1,\dots,r_0$. Note that the left hand side of this equation is
equal to $\hol_j(S_{0,\tau,T})$ whenever $T_j>\tau+1$.
There is then a constant $\eps>0$ such that
$\zeta(B_\eps)\subset\cc'_1$, in which case we have a composite map
\[\Xi=\hat\Xi\circ\zeta:B_\eps\to E.\]

Choose $\lla>0$ so that none of the operators $\ti H_{\al_j}$ ($j=1,\dots,r$)
and $\ti H_{\al'_j}$ ($j=1,\dots,r'$) has any eigenvalue of absolute value
$\le\lla$. (The notation $\ti H_\al$ was introduced in
\cite[Subsection~6.1]{Fr10}.)
Recall that we assume the $\si_j$ are small and non-negative, so
in particular we may assume $6\si<\lla$.

\be{lemma}\label{Xi-conv}
There is a constant $\xx2<\infty$ such that
\[\|\Xi(0)\|\le \xx2e^{(\si-\lla)\tau}.\]
\end{lemma}

\proof The first two components of $\Xi(0)=\hat\Xi(S_{0,\tau,T})$ are in fact
zero: the first one because $S_{0,\tau,T}=S_0$ over $K$, the second one
because the $dt$--component of 
\[S_{0,\tau}-S_0=(1-\sr_\tau)(\ual_j-S_0)\]
vanishes on $[1,\infty)\times(\pm Y_j)$ since $S_0$ and $\ual_j$ are both
in temporal gauge there.

The third component of $\Xi(0)$ is $\Theta'(S_{0,\tau,T})$. It suffices to
consider $\tau$ so large that the $\prt$--perturbations do not contribute to
$\Theta'(S_{0,\tau,T})$, which then vanishes outside the two bands of
length~$1$ in $[-T_j,T_j]\times Y_j$ centred at $t=\pm(T_j-\tau-\frac12)$,
$j=1,\dots,r$.
Our exponential decay result from
\cite{Fr10} says that for every $k\ge0$ there is a
constant $C'_k$ such that for every $(t,y)\in\R_+\times(\pm Y_j)$ we have
\[|\nabla^k(S_0-\ual_j)|_{(t,y)}\le C'_ke^{-\lla t}.\]
Consequently,
\[\|\Theta'(S_{0,\tau,T})\|_\infty\le\const(e^{-\lla\tau}+e^{-2\lla\tau})
\le\const e^{-\lla\tau}.\]
This yields
\[\|\Xi(0)\|=\|\Theta'(S_{0,\tau,T})\|_{\lw p\ka}
\le\const e^{(\si-\lla)\tau}.\square\]

\be{lemma}\label{Xi-bounds1}
There is a constant $\xx3<\infty$ such that for sufficiently large $\tau$
the following hold:
\be{description}
\item[(i)]$\|D\Xi(0)\|\le \xx3$,
\item[(ii)]$D\Xi(0)$ is invertible and
$\|D\Xi(0)\inv\|\le \xx3$.
\end{description}
\end{lemma}


\proof
By construction, the derivative of $\Xi$ at $0$ has the
form
\[D\Xi(0)=
\left(\be{array}{ccc}
I & 0 & \beta_1 \\
\del_2 & I & \beta_2 \\
\del_3 & 0 & I
\end{array}\right),\]
where the $k$'th column is the $k$'th partial derivative and
$I$ the identity map.

The middle top entry in the above matrix is zero because $\phol_{c,\tau}$
vanishes on $K$. The middle bottom entry is zero because $S_{0,\tau,T}=\ual_j$
on the support of $\phol_{c,\tau}$ and the spinor field of $\ual_j$
is zero (for $j=1,\dots,r_0$). Adding $\phol_{c,\tau}$ to
$S_{0,\tau,T}$ therefore has the effect of altering the latter by
a gauge transformation over $[-T_j+\tau+1,-T_j+\tau+2]\times Y_j$,
$j=1,\dots,r_0$.

We claim that $\beta_k$ is bounded in $\tau,T$ for $k=1,2$. For $k=1$
this is obvious from the boundedness of $Q_{\tau,T}$. For $k=2$
note that the derivative of $\eta_j:\cc'_1\to(-\pi,\pi)$ at any $S\in\cc'_1$
is
\be{equation}\label{D-eta}
D\eta_j(S)(a,\phi)=i\int_{I_j}\ga^*_ja
\end{equation}
where $a$ is an imaginary valued $1$--form and $\phi$ a positive spinor.
Because of the weights
used in the Sobolev norms it follows that $D\eta(S)$ is 
(independent of $S$ and) bounded in $\tau,T$
(see \Ref{ewa}). This together with the bound on $Q_{\tau,T}$
gives the desired bound on $\beta_2$.

Note that, for $k=2,3$, $\|\del_k\|$ is independent of $T$ when $\tau\gg0$,
and routine calculations show that $\|\del_k\|\to0$ as $\tau\to\infty$.
(In the case of $\del_2$ this depends on the normalization~\Ref{z-cond}
of the holonomy of $\bds(v)$.)

Write $D\Xi(0)=x-y$, where
\[x=\left(\be{array}{ccc}
I & 0 & \beta_1 \\
0 & I & \beta_2 \\
0 & 0 & I
\end{array}\right),
\qquad
x\inv=\left(\be{array}{ccc}
I & 0 & -\beta_1 \\
0 & I & -\beta_2 \\
0 & 0 & I
\end{array}\right).\]
When $\tau$ is so large that $\|y\|\,\|x\inv\|<1$ then
of course $\|yx\inv\|\le\|y\|\,\|x\inv\|<1$, hence
$x-y=(I-yx\inv)x$ is invertible. Moreover,
\[(x-y)\inv-x\inv=x\inv[(I-yx\inv)\inv-I]
=x\inv\sum_{k=1}^\infty(yx\inv)^k,\]
which gives
\[\|(x-y)\inv-x\inv\|\le
\frac{\|x\inv\|^2\,\|y\|} {1-\|x\inv\|\,\|y\|}
\to0\quad\text{as $\tau\to\infty$}.\square\]

We now record some basic facts that will be used in 
the proof of Lemma~\ref{Xi-bounds2} below.

\be{lemma}\label{lemma:chain-rule}
If $E_1,E_2,E_3$ are Banach spaces, $U_j\subset E_j$ an open set
for $j=1,2$, and $f:U_1\to U_2$, $g:U_2\to E_3$ smooth maps then the
second derivate of the composite map $g\circ f:U_1\to E_3$ is given by
\be{equation*}
\be{aligned}
D^2(g\circ f)(x)(y,z)&=D^2g(f(x))(Df(x)y,Df(x)z)\\
&+Dg(f(x))(D^2f(x)(y,z))
\end{aligned}
\end{equation*}
for $x\in U_1$ and $y,z\in E_1$.
\end{lemma}

\proof Elementary.\square

It is also worth noting that embedding and multiplication theorems for
$L^q_k$ Sobolev spaces on $\rf$ ($k\ge0$, $1\le q<\infty$) carry over to
$\xt$, and that the embedding and multiplication constants are bounded
functions of $T$.

Furthermore, a differential operator of degree $d$ over $\xt$ which is 
translationary invariant over necks and ends induces a bounded operator
$L^q_{k+d}\to L^q_k$ whose operator norm is a bounded function of $T$.

\be{lemma}\label{Xi-bounds2}
There is a constant $\xx4>0$ such that
$\|D^2\Xi(x)\|\le \xx4$ whenever $\|x\|\le \xx4\inv$ and $\tau\ge \xx4$.
\end{lemma}

\proof We will say a quantity depending on $x,\tau$
is {\em s-bounded} if the lemma holds with this quantity in place of $D^2\Xi$.

Let $\Xi_1,\Xi_2,\Xi_3$ be the components of $\Xi$.

The assumption $K\subset X\endt \tau$ ensures that
$\Xi_1(v,c,\xi)$ is independent of $c$. It
then follows from Lemma~\ref{lemma:chain-rule} and the bound on
$Q_{\tau,T}$ that
$D^2\Xi_1$ is s-bounded.

When $v,c,\xi$ are small we have
\[\Xi_2(v,c,\xi)=\eta(S_{v,\tau,T}+\phol_{c,\tau}+Q_{\tau,T}\xi)
=c+\eta(S_{v,\tau,T}+Q_{\tau,T}\xi).\]
Since $D\eta$ is constant, as noted above, we have $D^2\eta=0$. 
From the bounds on $D\eta$ and $Q_{\tau,T}$ we then deduce that
$D^2\Xi_2$ is s-bounded.

To estimate $\Xi_3$, we fix $h\gg0$ and consider only $\tau\ge h$.
It is easy to see that
\[\Xi_3(x)|_{X\endt h}\]
is s-bounded. By restricting to small $x$ and choosing $h$ large we may
arrange that the $\prt$--perturbations do not contribute to
\[\Xi_{3,j}:=\Xi_3|_{[-T_j+h,T_j-h]\times Y_j}\]
for $j=1,\dots,r$. We need to show that each $D^2\Xi_{3,j}$ is s-bounded, but
to simplify notation we will instead prove the same for $D^2\Xi_3$ under
the assumption that the $\prt$ perturbations are zero.

First observe that for any configuration $(A,\Phi)$ over $\xt$
and any closed, imaginary valued $1$--form $a$ we have
\[\Theta'(A+a,\Phi)=\Theta'(A,\Phi)+(0,a\cdot\Phi).\]
Moreover,
\[\|a\cdot\Phi\|_{\lw p\ka}=C\|a\cdot e^\ka\Phi\|_p\le
C\|a\|_{2p}\,\|e^\ka\Phi\|_{2p}\le C'\|a\|_{2p}\,\|\Phi\|_{\llw p\ka1}\]
for some constants $C,C'<\infty$.
Taking $(A,\Phi)=S_{v,\tau,T}+Q_{\tau,T}\xi$ and $a=\phol_{c,\tau}$
we see that $D_kD_2\Xi_3(x)$ is s-bounded for $k=1,2,3$.

Next note that the derivative of the Seiberg--Witten
map $\Theta':\cc'\to\lw p\ka$ at a point $S'\rr+s_1$ has the form
\[D\Theta'(S'\rr+s_1)s_2=Ls_2+B(s_1,s_2)\]
where $B$ is a pointwise bilinear operator, and $L$
a first order operator which is independent of $s_1$ and translationary
invariant over necks and ends. This yields
\be{equation}\label{DTheta-bound}
\|D\Theta'(S\rr'+s)\|\le \const(1+\|s\|_{2p}).
\end{equation}
Moreover, $D^2\Theta'(S)=B$
for all $S$, hence there is a constant $C''<\infty$ such that
\be{equation*}
\|D^2\Theta'(S)\|\le C''
\end{equation*}
for all $T$.

Combining the above results on $\Theta'$ with Lemma~\ref{lemma:chain-rule}
we see that
$D_jD_k\Xi_3$ is s-bounded also when $j,k\neq2$.\square

\be{lemma}\label{Xi-inv}
There exist $c_5>0$ and $\xx6<\infty$
such that if $0<\eps'<c_5\eps<c_5^2$ then for 
sufficiently large $\tau$ the following hold:
\be{description}
\item[(i)]$\Xi:B_\eps\to E$ is injective,
\item[(ii)]There is a (unique)
smooth map $\Xi\inv:B_{\eps'}\to B_\eps$ such that
$\Xi\circ\Xi\inv=I$,
\item[(iii)]$\|D(\Xi\inv)(x)\|\le \xx6$ for all $x\in B_{\eps'}$,
\item[(iv)]$\|D^2(\Xi\inv)(x)\|\le \xx6$ for all $x\in B_{\eps'}$,
\item[(v)]$\|\Xi\inv(0)\|\le \xx6e^{(\si-\lla)\tau}$.
\end{description}
\end{lemma}

\proof For sufficiently large $\tau$ we have
\[\eps'+\|\Xi(0)\|<c_5\eps\]
by Lemma~\ref{Xi-conv}. Statements (i)--(iv) now follow from the inverse
function theorem, \cite[Proposition~B.1]{Fr10}, applied to the function
$x\mapsto\Xi(x)-\Xi(0)$, together with Lemmas~\ref{Xi-bounds1} and
\ref{Xi-bounds2}. To prove (v), set $h=\Xi\inv$,
$x=\Xi(0)$ and take $\tau$ so large that $x\in B_{\eps'}$.
Since $\Xi$ is injective on $B_\eps$ we must have $h(x)=0$, so
\[\|h(0)\|=\|h(x)-h(0)\|\le\|x\|\sup_{\|y\|\le\eps'}\|Dh(y)\|.\]
Now (v) follows from (iii) and Lemma~\ref{Xi-conv}.
\square


From now on we assume that $\eps,\eps',\tau$ are chosen so that the conclusions
of the lemma are satisfied.
Define 
\[\hzeta=\zeta\circ\Xi\inv:B_{\eps'}\to\cc'_1.\]
Then clearly
\[\hxi\circ\hzeta=I.\]
Thus, $(v,c)\mapsto\hzeta(v,c,0)$
is a ``gluing map'', ie for small $v,c$ it solves the problem
of gluing the monopole $\bds(v)$ over $X$
to get a monopole over $\xt$ with prescribed
holonomy $z_je^{ic_j}$ along the path $\ga_j$ for $j=1,\dots,r_0$.

\be{lemma}\label{hzeta-bounds1}
There is a constant $\xx7<\infty$ such that for $x\in B_{\eps'}$ one has
\[\|D\hzeta(x)\|,\|D^2\hzeta(x)\|\le \xx7e^{\si\tau}.\]
\end{lemma}

\proof This follows from Lemmas~\ref{zeta-bounds} and \ref{Xi-inv}
and the chain rule.\square

The following proposition refers to the situation of
Subsection~\ref{subsec:gluing-statement1} and uses the notation
of Theorem~\ref{gluing-thm1}. 

\be{prop}\label{prop:surj}
If $(K',V')$ is any kv-pair $\le(K,V)$
then $\oline G\times\ur$ can be covered by finitely many connected
open sets
$W$ in $M_\hx\times\ur$ such that if $\tmin$ is sufficiently
large then for each $W$
there exists a smooth map $\mathbf h:W\to\Ht$ whose image consists only of 
regular points and which satisfies $\ff\circ\mathbf h=I$.
\end{prop}

Here we do not need any
assumptions on $\ti\eta_j,\ti\eta'_j$ or on $\prt_j,\prt'_j$.

\proof Let $(\om_0,z)\in\oline G\times\ur$ and consider the set-up above,
with $\tau$ so large that $K'\subset X\endt \tau$ and $\eps$ so small that
\be{equation}\label{zkv}
\zeta(x)|_{K'}\in V'\qquad\text{for every $x\in B_\eps$}.
\end{equation}
Note that taking $\eps$ small may require taking $\tau$ (and hence $\tmin$)
large, see Lemma~\ref{Xi-inv}.
For any sufficiently small open neighbourhood $W\subset M_\hx\times\ur$
of $(\om_0,z)$ we can define a smooth map $\nu:W\to\cc'_1$ by the formula
\[\nu(\om,a)=
\hat\zeta(\varpi\inv(\om),-i\log\frac az,0).\]
Here $\log e^u=u$ for any complex number $u$ with $|\Im\,u|<\pi$, and
$i\log\frac az\in\R^{r_0}$ denotes the vector whose $j$'th component is
$i\log\frac{a_j}{z_j}$.
Because
$\hat\Xi\circ\hat\zeta=I$ and the Seiberg--Witten map is the third 
component of $\hat\Xi$, the image of $\nu$ consists of regular monopoles.
Let $\mathbf h:W\to\cb'_\hx$ be the composition of $\nu$ with the projection
$\cc'\to\cb'_\hx$. Unravelling the definitions involved
and using \Ref{zkv} one finds that $\mathbf h$ has the required properties.

How large $\tmin$ must be for this to work might depend on $(\om_0,z)$.
But $\oline G\times\ur$ is compact, hence it can be covered by
finitely many such open sets $W$. If $\tmin$ is sufficiently large
then the above construction will work for each of these $W$.
\square

\subsection{Injectivity}\label{subsec:inj}

We now continue the discussion that was interrupted by
Proposition~\ref{prop:surj}. 
Set
\[\ti S=S_{0,\tau,T},\quad\hat S=\hzeta(0).\]

\be{lemma}\label{hzeta-bounds2}
There is a constant $\xx8<\infty$ such that for sufficiently large $\tau$
one has
\[\|\hat S-\ti S\|_{\llw p\ka1}\le \xx8e^{(2\si-\lla)\tau},\quad
\|\hat S-S'\rr\|_{\llw p\ka1}\le \xx8.\]
\end{lemma}

\proof Set $\Xi\inv(0)=(v,c,\xi)\in B_\eps$.
For sufficiently large $\tau$ we have
\be{align*}
\|\hat S-\ti S\|_{\llw p\ka1}
&\le\|S_{v,\tau,T}-S_{0,\tau,T}\|_{\llw p\ka1}+
     \|\phol_{c,\tau}+Q_{\tau,T}\xi\|_{\llw p\ka1}\\
&\le\const(\|v||+e^{\si\tau}\|c\|+\|\xi\|)\\
&\le\const e^{(2\si-\lla)\tau},
\end{align*}
where we used 
Lemma~\ref{Xi-inv}\,(v) to obtain the last inequality.
Because $\|\ti S-S'\rr\|_{\llw p\ka1}$
is bounded in $\tau,T$, and we assume $6\si<\lla$,
the second inequality of the lemma follows as well.\square

For positive spinors $\Phi$ on $\xt$ it is convenient to extend the definition
of $\ci_\Phi$ to complex valued functions on $\xt$:
\[\ci_\Phi f=(-df,f\Phi).\]
(However, $\ci^*_\Phi$ will always refer to the formal adjoint
of $\ci_\Phi$ acting on imaginary valued functions.)
When $\Phi$ is the spinor part of $S'\rr,\ti S,\hat S$ 
then the corresponding operators $\ci_\Phi$ will be denoted
$\ci'\rr,\ti\ci,\ci$, respectively.
(We omit the \^{} on $\ci$ to simplify notation.)
As in \cite[Subsection~2.2]{Fr10} we define
\[\ce'=\{f\in L^p_{2,\loc}(\xt;\co)\st \ci'\rr f\in\llw p\ka1\}.\]
We can take the norm to be
\[\enm f=\|\ci'\rr f\|_{\llw p\ka1}+\sum_{x\in\hx}|f(x)|.\]

\be{lemma}\label{sup-ci}
There is a constant $\xx9<\infty$ such that if $\bdi$ is any of the operators
$\ci'\rr,\ti\ci,\ci$ then for all $f\in\ce'$ one has
\[\|f\|_\infty
\le \xx9\left(\|\bdi f\|_{\lw p\ka}+\sum_{x\in\hx}|f(x)|\right).\]
\end{lemma}

\proof We first prove the inequality for $\bdi=\ti\ci$ (the case of $\ci'\rr$
is similar, or easier). If $X_e$ is any component of $X$ and
$0\le\bar\tau\le\tau$ then for some constant $C_{\bar\tau}<\infty$ one has
\[\|f\|_\infty\le\const\|f\|_{L^p_1}\le
C_{\bar\tau}\left(\|\ti\ci f\|_p+\sum_{x\in\hx\cap X_e}|f(x)|\right)\]
for all $L^p_1$ functions $f:(X_e)\endt{\bar\tau}\to\co$.
Here the Sobolev inequality holds because
$p>4$, whereas the second inequality follows from \cite[Lemma~2.1]{Fr10}.
We use part~(i) of that lemma if the spinor field of $S_0$ is not
identically zero on $X_e$, and part~(ii) otherwise. (In the latter case
$\hx\cap X_e$ is non-empty.)

When $\bar\tau,\tau$ are sufficiently large we can apply part~(i) of the same
lemma in a similar fashion to the band $[t,t+1]\times Y'_j$ provided
$t\ge\bar\tau$ and $\al'_j$ is irreducible, and to the band
$[t-1,t+1]\times Y_j$ provided $|t|\le T_j-\bar\tau-1$ and $\al_j$ is
irreducible. To estimate $|f|$ over these bands
when $\al'_j$ resp.\ $\al_j$ is reducible one can use 
\cite[Lemma~2.2\,(ii)]{Fr10}. This proves the lemma for $\bdi=\ti\ci$ (and
for $\bdi=\ci'\rr$).

We now turn to the case $\bdi=\ci$.
Let $\phi$ denote the spinor part of $\hat S-\ti S$. Then
\be{align*}
\|\ti\ci f\|_{\lw p\ka}
&\le\|\ci f\|_{\lw p\ka}+\|f\phi\|_{\lw p\ka}\\
&\le \|\ci f\|_{\lw p\ka}
  +\const\left(\|\ti\ci f\|_{\lw p\ka}+\sum_{x\in\hx}|f(x)|\right)\,
  \cdot\|\phi\|_{\lw p\ka}.
\end{align*}
By Lemma~\ref{hzeta-bounds2}
we have $\|\phi\|_{\llw p\ka1}\to0$ as $\tau\to0$,
so for sufficiently large $\tau$ we get
\[\|\ti\ci f\|_{\lw p\ka}
\le\const\left(\|\ci f\|_{\lw p\ka}+\sum_{x\in\hx}|f(x)|\right).\]
Therefore, the lemma holds with $\bdi=\ci$ as well.
\square

\be{lemma}\label{ci-prop}
There is a constant $\xx{10}<\infty$ such that for all $f,g\in\ce'$ and
$\phi\in\llw p\ka1(\xt;\bs^+)$ one has
\be{description}
\item[(i)]$\|fg\|\le \xx{10}\|f\|\,\|g\|$,
\item[(ii)]$\|f\phi\|\le \xx{10}\|f\|\,\|\phi\|$,
\end{description}
where we use the $\llw p\ka1$ norm on spinors and the $\ce'$ norm on elements
of $\ce'$.
\end{lemma}

\proof By routine calculation using
Lemma~\ref{sup-ci} with $\bdi=\ci'\rr$ one easily proves (ii) and
the inequality
\[\lnm{d(fg)}\le\const\enm f\enm g.\]
Now observe that by definition
$g\Phi'\rr\in\llw p\ka1$, where as before $\Phi'\rr$ denotes the spinor
field of the reference configuration $S'\rr$.
Applying (ii) we then obtain
\[\lnm{fg\Phi'\rr}\le\const\enm f\,\lnm{g\Phi'\rr}\le\const\enm f\,\enm g,\]
completing the proof of (i).\square

Recall from \cite[Subsection~2.4]{Fr10} that the Lie algebra $L\cg'_\hx$ 
is the space of imaginary valued functions in $\ce'$ that vanish on $\hx$.

\be{lemma}\label{ci-compare}
There is a constant $\xx{11}>0$ such that for $\tau>\xx{11}$ and
all $f\in L\cg'_\hx$ one has
\[\xx{11}\inv\|\ci'\rr f\|_{\llw p\ka1}\le\|\ci f\|_{\llw p\ka1}\le
\xx{11}\|\ci'\rr f\|_{\llw p\ka1}.\]
\end{lemma}

\proof 
Let $\psi$ denote the spinor part of $\hat S-S'\rr$. Then
\be{align*}
\|f\psi\|_{\llw p\ka1}
&\le\const(\|f\|_\infty\|\psi\|_{\llw p\ka1}
  +\|df\|_{\lw{2p}\ka}\|\psi\|_{2p})\\
&\le\const\|\ci f\|_{\llw p\ka1}\|\psi\|_{\llw p\ka1},
\end{align*}
and similarly with $\ci'\rr$ instead of $\ci$.
The lemma now follows from Lemma~\ref{hzeta-bounds2}.\square

We are going to use the inverse function theorem a second time,
to show that the image of the smooth map 
\be{align*}
\Pi:L\cg'_\hx\times B_{\eps'}&\to\cc'_1,\\
(f,x)&\mapsto\exp(f)(\hat\zeta(x)).
\end{align*}
contains a ``not too small'' neighbourhood of $\hat S$.
The derivative of $\Pi$ at $(0,0)$ is
\be{align*}
D\Pi(0,0):L\cg'_\hx\oplus E&\to\llw p\ka1,\\
(f,x)&\mapsto\ci f+D\hat\zeta(0)x.
\end{align*}
To be concrete, let $L\cg'_\hx\oplus E$ have the norm
$\|(f,x)\|=\enm f+\|x\|_E$.

\be{lemma}\label{Pi-Fredholm}
$D\Pi(0,0)$ is a linear homeomorphism.
\end{lemma}

\proof By \cite[Proposition~2.2]{Fr10}, $\ci^*\ci:L\cg'_\hx\to\lw p\ka$
is a Fredholm operator
with the same kernel as $\ci$. Now, $\ci$ is injective on $L\cg'_\hx$,
because $\hzeta$ maps into $\cc'_1$ and therefore
$[\hat S|_K]\in V\subset\cb^*_\hx$. Since
\[W=\ci^*\ci(L\cg'_\hx)\]
is a closed subspace of $\lw p\ka$ of finite codimension, we can choose
a bounded operator
\[\pi:\lw p\ka\to W\]
such that $\pi|_W=I$. Set
\[\ci^\#=\pi\ci^*:\llw p\ka1\to W.\]
Then
\[\ci^\#\ci:L\cg'_\hx\to W\]
is an isomorphism. Furthermore, 
\[\index(\ci^\#+D\Theta'(\hat S))=\dim\,\mt_\hx=\dim M_\hx+r_0,\]
where `dim' refers to expected dimension (which in the case of $M_\hx$
is equal to the actual dimension of $G$), and the second equality follows from
the addition formula for the index (see Corollary~\ref{cor:add-formula}).
Consequently,
\[\index(\ci^\#+D\hat\Xi(\hat S))=0.\]
We now compute
\be{equation}\label{Pi-matrix}
(\ci^\#+D\hat\Xi(\hat S))\circ D\Pi(0,0)=
\left(\be{array}{cc}
\ci^\#\ci & B\\
0 & I
\end{array}\right)
:L\cg'_\hx\oplus E\to W\oplus E,
\end{equation}
where $B:E\to W$. The zero in the matrix above is due to the fact that
\[D\hat\Xi(\hat S)\ci f=\left.\frac d{dt}\right|_0\hat\Xi(e^{tf}(\hat S))=0,\]
which holds because $\hxi_1,\hxi_2$ are 
$\cg'_\hx$--invariant, $\hxi_3$ is $\cg'_\hx$--equivariant,
and $\hxi(\hat S)=0$.

Since the right hand side of
\Ref{Pi-matrix} is invertible, it follows that $\ci^\#+D\hat\Xi(\hat S)$
is a surjective Fredholm operator of index~$0$, hence invertible. Of course,
this implies that $D\Pi(0,0)$ is also invertible.\square

\be{lemma}\label{Pi-bounds1}
There is a constant $\xx{12}<\infty$ such that for sufficiently large $\tau$,
\[\|D\Pi(0,0)\inv\|\le \xx{12}e^{\si\tau}.\]
\end{lemma}

\proof In this proof all unqualified norms are $\llw p\ka1$ norms.
It follows from \Ref{D-eta}, \Ref{DTheta-bound} and
Lemma~\ref{hzeta-bounds2} that $D\hxi(\hat S)$ is bounded in $\tau,T$.
Therefore there exists a constant $C<\infty$ such that
\[\|x\|_E=\|D\hxi(\hat S)(\ci f+D\hzeta(0)x)\|\le C\|D\Pi(0,0)(f,x)\|\]
for all $f\in L\cg'_\hx$ and $x\in E$. From Lemma~\ref{ci-compare}
and Lemma~\ref{hzeta-bounds1} we get
\be{align*}
\xx{11}\inv\|\ci'\rr f\|&\le\|\ci f\|\\
&\le\|D\Pi(0,0)(f,x)\|+\|D\hzeta(0)x\|\\
&\le\|D\Pi(0,0)(f,x)\|+\xx7e^{\si\tau}\|x\|_E\\
&\le(1+C\xx7e^{\si\tau})\|D\Pi(0,0)(f,x)\|.
\end{align*}
This yields
\[\enm f+\|x\|_E\le\const e^{\si\tau}\|D\Pi(0,0)(f,x)\|.\square\]

\be{lemma}\label{Pi-bounds2}
There is a constant $\xx{13}<\infty$
such that for sufficiently large $\tau$ one has
\[\|D^2\Pi(f,x)\|\le \xx{13}e^{\si\tau}\]
for all $f\in L\cg'_\hx$, $x\in E$ such that $\|f\|<1$ and $\|x\|<\eps'$.
\end{lemma}

\proof For the purposes of this proof it is convenient to rescale the norm
on $\ce'$ so that we can take $\xx{10}=1$ in Lemma~\ref{ci-prop}.

If $f,g\in\ce'$ then $e^fg\in\ce'$, and from Lemma~\ref{ci-prop} we obtain
\[\|e^fg\|\le\sum_{n=0}^\infty\frac1{n!}\|f^ng\|\le e^{\|f\|}\|g\|,\]
and similarly with a spinor $\phi\in\llw p\ka1$ instead of $g$.

The first two derivatives of $\exp:\ce'\to1+\ce'$ are
\[D\exp(f)g=g\exp(f),\qquad D^2\exp(f)(g,h)=gh\exp(f),\]
so
\be{equation*}
\|D\exp(f)\|,\|D^2\exp(f)\|\le \exp(\|f\|).
\end{equation*}

Let $\hzeta_1,\hzeta_2$ be the connection and spinor
parts of $\zeta$, respectively, and define $\Pi_1,\Pi_2$ similarly.
Then
\[\Pi(f,x)=(\hzeta_1(x)-df,e^f\cdot\hzeta_2(x)).\]
We regard $\Pi(f,x)$ as a function of the the two variables $f,x$.
Let $D_j\Pi$ denote the derivative of $\Pi$ with respect to the $j$'th
variable. Similarly for the second derivatives $D_jD_k\Pi$.

Applying Lemmas~\ref{hzeta-bounds1} and \ref{ci-prop} we now find that
\be{gather*}
\|D_2^2\Pi_1(f,x)\|=\|D^2\hzeta_1(x)\|\le\const e^{\si\tau},\\
\|D_jD_k\Pi_2(f,x)\|\le\const e^{\si\tau},\quad j,k=1,2
\end{gather*}
for $\|f\|<1$ and $\|x\|<\eps'$.
Since $D_jD_1\Pi_1=0$ for $j=1,2$, the lemma is proved.\square

In the following, $B(x;r)$ will denote the open $r$--ball about $x$ (both in
various Banach spaces and in $\cc'$).

\be{lemma}\label{Pi-image}
There exist constants $r_1,r_2>0$ such that for sufficiently large $\tau$
the image of $\Pi$ contains
the ball of radius $r_2e^{-3\si\tau}$ about $\hat S$ in $\cc'$;
more precisely one has
\be{equation*}
B(\hat S;r_2e^{-3\si\tau})\subset\Pi(B(0;r_1e^{-2\si\tau})).
\end{equation*}
\end{lemma}

\proof We wish to apply the inverse function theorem 
\cite[Proposition~B.1]{Fr10} to the map $\Pi$ 
restricted to a ball $B(0;R_1)$, where $R_1\in(0,\eps']$ is to be chosen.
For the time being let $M,L,\ka$ have the same meaning as in that proposition.
By Lemma~\ref{Pi-bounds2} we can take $M=\xx{13}e^{\si\tau}$, and by
Lemma~\ref{Pi-bounds1} we have $\|L\inv\|\le \xx{12}e^{\si\tau}$. We need
\[0\le\ka=\|L\inv\|\inv-R_1M.\]
This will hold if
\[R_1\le(\xx{12}\xx{13})\inv e^{-2\si\tau}.\]
When $\tau$ is large we can take $R_1$ to be the right hand side of this
inequality. By \cite[Proposition~B.1]{Fr10}, $\Pi(B(0;R_1))$
contains the ball $B(\hat S;R_2)$ where
\[R_2=\frac12R_1\xx{12}\inv e^{-\si\tau}=
\frac12\xx{12}^{-2}\xx{13}\inv e^{-3\si\tau}.
\square\]


Theorem~\ref{gluing-thm1} is a consequence of Proposition~\ref{prop:surj}
and the following proposition:

\be{prop}\label{prop:inj}
Under the assumptions of Theorem~\ref{gluing-thm1}, and using the
same notation, there is a kv-pair $(K',V')\le (K,V)$ such that
$\ff$ is injective on $\qq\inv(\oline G)$ for all
sufficiently large $\tmin$.
\end{prop}

The proof of Proposition~\ref{prop:inj} occupies the remainder of
this subsection.


For any natural number $m$ which is so large that $K\subset X\endt m$,
let $V'_m$ be the set of all $\om\in\tcb_\hx(X\endt m)$ such that
there exist a representative $S$ of $\om$, and a
configuration $\bar S=(\bar A,\bar\Phi)$ over $X$ representing
an element of $\oline G$, such that
\be{equation}\label{eqn:dm-defn}
d_m(S,\bar S):=\int_{X\endt m}|\bar S-S|^p+|\nabla_{\bar A}(\bar S-S)|^p
<\frac1m.
\end{equation}
Note that
\[d_m(u(S),u(\bar S))=d_m(S,\bar S)\]
for any gauge transformation $u$ over $X\endt m$. 
In particular, $V'_m$ is $\bt$--invariant. 


\be{lemma}\label{vp-conv1}
Let $\om_n\in V'_{m_n}$ for $n=1,2,\dots$, where $m_n\to\infty$.
Then there exists for each $n$ a representative $S_n$ of $\om_n$
such that a subsequence of $S_n$ converges locally in $L^p_1$ over $X$ to a
smooth configuration representing an element of $\oline G$.
\end{lemma}


\proof By assumption there exist for each $n$ a representative $\sS_n$
of $\om_n$ and a configuration $\bar{\sS}_n$ over $X$ representing an
element of $\oline G$ such that
\be{equation}\label{dn-assump}
d_{m_n}(\sS_n,\bar{\sS}_n)<\frac1{m_n}.
\end{equation}
After passing to a subsequence we may assume (since $\oline G$ is compact)
that $[\bar{\sS}_n]$ converges in $\oline G$ to
some element $[\bar S]$, and we can choose $\bar S$ smooth.
Since $M_\hx=M^*_\hx$, the local slice theorem guarantees that
for large $n$ we can find $u_n\in\cg_\hx$ such that
$\bar S_n=u_n(\bar{\sS}_n)$ satisfies
\[\|\bar S_n-\bar S\|_{\llw pw1}\to0.\]
Set $S_n=u_n(\sS_n)$, which is again a representative of $\om_n$.
Let $\bar A,\bar A_n$ be the
connection parts of $\bar S,\bar S_n$,
respectively. Then \Ref{dn-assump} implies that
$\bar S_n-S_n\to0$ and $\nabla_{\bar A_n}(\bar S_n-S_n)\to0$
locally in $L^p$ over $X$,
hence also $S_n\to \bar S$ locally in $L^p$ over $X$. Now
\[\nabla_{\bar A}(S_n-\bar S)=
\nabla_{\bar A_n}(S_n-\bar S_n)+\nabla_{\bar A}(\bar S_n-\bar S)
+(\bar A-\bar A_n)(S_n-\bar S_n),\]
and each of the three terms on the right hand side converges to $0$ locally
in $L^p$ over $X$ (the third term because of the continuous multiplication
$L^p_1\times L^p\to L^p$ in $\rf$ for $p>4$).
Hence $S_n\to \bar S$ locally in $L^p_1$ over $X$.
\square

\be{cor}For sufficiently large $n$ one has that $R_K(V'_n)\subset V$.\square
\end{cor}

\be{lemma}\label{vp-conv2}
Let $\om_n\in V'_{m_n}$ for $n=1,2,\dots$, where $m_n\to\infty$.
Suppose
$q(\om_n|_K)$ converges in $M_\hx$ to an element $g$ as $n\to\infty$.
Then $g\in\oline G$, and there exists
for each $n$ a representative $S_n$ of $\om_n$
such that the sequence $S_n$ converges locally in $L^p_1$ over $X$ to a smooth
configuration representing $g$.
\end{lemma}

\proof Let $\sS_n,\bar{\sS}_n$ be as in the proof of Lemma~\ref{vp-conv1}.
First suppose that $[\bar{\sS}_n]$ converges in $\oline G$ to
some element $[\bar S]$, where $\bar S$ is smooth. Choosing $S_n,\bar S_n$
as in that proof we find again that $S_n\to \bar S$ locally
in $L^p_1$ over $X$, hence
\[g=\lim_n q(S_n|_K)=q(\bar S|_K)=[\bar S].\]
We now turn to the general case when $[\bar{\sS}_n]$
is not assumed to converge.
Because $\oline G$ is compact, every subsequence of
$[\bar{\sS}_n]$ has a convergent
subsequence whose limit must be $g$ by the above argument. Hence 
$[\bar{\sS}_n]\to g$.\square

Suppose we are given a sequence $\{m_n\}_{n=1,2,\dots}$ of natural numbers
tending to infinity, and for each $n$
an $r$--tuple $T(n)$ of real numbers such that
\[\tmin(n):=\min_jT_j(n)>m_n.\]
Define $\qq_n$ and $\ff_n$ as in Theorem~\ref{gluing-thm1}, with
$K'=X\endt{m_n}$ and $V'=V'_{m_n}$.

\be{lemma}\label{s-approx}
For $n=1,2,\dots$ suppose $S_n$ is a smooth configuration over
$X^{(T(n))}$
representing an element $\om_n\in \qq_n\inv(\oline G)$, and such that
\[\ff_n(\om_n)\to(\om_0,z)\in\oline G\times\ur\]
as $n\to\infty$. There exists a constant $\xx{14}<\infty$ such that for
sufficiently large $\tau$ the following holds for sufficiently large
$n$. Let the map $\hzeta=\hzeta_n$ be
defined as above and set $\hat S_n =\hzeta_n (0)$. Then there exists
a smooth gauge transformation $u_n\in\cg'_\hx$ such that
\[\|u_n(S_n)-\hat S_n \|_{\llw p\ka1}
\le \xx{14}e^{(3\si-\lla)\tau}.\]
\end{lemma}

Note: This constant $\xx{14}$ depends on $(\om_0,z)$ but not on the sequence
$S_n$.

Before proving the lemma, we will use it to show that $\ff_n$ is injective
on $\qq_n\inv(\oline G)$
for some $n$. This will prove Proposition~\ref{prop:inj}. Suppose
$\om_n,\om'_n\in \qq_n\inv(\oline G)$ and
$\ff_n(\om_n)=\ff_n(\om'_n)$, $n=1,2,\dots$.
After passing to a subsequence we may assume that $\ff_n(\om_n)$ converges to
some point $(\om_0,z)\in\oline G\times\ur$.
Combining Lemmas~\ref{Pi-image},\,\ref{s-approx} and the assumption
$6\si<\lla$ we conclude that if $\tau$ is sufficiently large
then for sufficiently large $n$ we can represent $\om_n$ and $\om'_n$ by
configurations $\hzeta(x_n)$ and $\hzeta(x'_n)$,
respectively, where $x_n,x'_n\in B_{\eps'}$.
Now recall that $\hxi\circ\hzeta=I$, and that
the components $\hxi_1,\hxi_2$ are $\cg'_\hx$--invariant whereas $\hxi_3$
is the Seiberg--Witten map. Comparing the definitions of $\ff_n$ and
$\hxi$ we conclude that
\[x_n=\hxi(\hzeta(x_n)=\hxi(\hzeta(x'_n)=x'_n,\]
hence $\om_n=\om'_n$ for large $n$. To complete the
proof of Proposition~\ref{prop:inj} it therefore only remains to prove
Lemma~\ref{s-approx}.

{\em Proof of Lemma~\ref{s-approx}:}
In this proof, constants will be independent of the sequence $S_n$
(as well as of $\tau$ as before).

By Lemma~\ref{vp-conv2} we can find for each $n$ an $L^p_{2,\loc}$
gauge transformation $v_n$ over $X$ with $v_n|_\hx=1$
such that $S'_n=v_n(S_n)$ converges locally in $L^p_1$ over $X$ to a smooth
configuration $S'$ representing $\om_0$.
A moment's thought shows that we
can choose the $v_n$ smooth, and we can clearly arrange that $S'=S_0$.
Then for any $t\ge0$ we have
\be{equation}\label{S'-bound}
\limsup_n\|S'_n-\hat S_n \|_{\llw p\ka1(X\endt t)}=
\limsup_n\|S_0-\hat S_n \|_{\llw p\ka1(X\endt t)}
\le\const e^{(2\si-\lla)\tau}
\end{equation}
when $\tau$ is so large that Lemma~\ref{hzeta-bounds2} applies.

For $t\ge0$ and any smooth configurations $S$ over $X\endt t$
consider the functional
\be{gather*}
E(S,t)=\sum_{j=1}^r\lla_j
\left(\cd(S|_{\{t\}\times Y_j})
+\cd(S|_{\{t\}\times(-Y_j)})\right)\\
+\sum_{j=1}^{r'}\lla'_j
\left(\cd(S|_{\{t\}\times Y'_j})-\cd(\al'_j)\right),
\end{gather*}
where in this formula $\{t\}\times(\pm Y_j)$ has the boundary orientation
inherited from $X\endt t$. (Recall that the Chern-Simons-Dirac functional
$\cd$ changes sign when the orientation of the $3$--manifold in question
is reversed.)
The assumption on $\lla_j,\lla'_j$ and
$\ti\eta_j,\ti\eta'_j$ in Theorem~\ref{gluing-thm1} implies that
$E(S,t)$ depends only on the gauge equivalence class of $S$. Since
$\cd$ is a smooth function on the $L^2_{1/2}$
configuration space by \cite[Lemma~3.1]{Fr10}, we obtain
\[E(S_n,t)=E(S'_n,t)\to E(S_0,t)\]
as $n\to\infty$. By our exponential decay results (see the proof of
\cite[Theorem~6.1]{Fr10}),
\[E(S_0,t)<\const e^{-2\lla t}\quad\text{for $t\ge0$.}\]
It follows that
\[E(S_n,t)<\const e^{-2\lla t}\quad\text{for $n>N(t)$}\]
for some positive function $N$.
By assumption the perturbation parameters $\vec\prt,\vec\prt'$ are admissible,
hence there is a constant $C<\infty$ such that when $\tmin(m_n)>C$, each of the
$(r+r')$ summands appearing in the definition of $E(S_n,t)$ is non-negative.
Explicitly, this yields
\be{gather*}
0\le\cd(S_n|_{\{-T_j(n)+t\}\times Y_j})-\cd(S_n|_{\{T_j(n)-t\}\times Y_j})
<\const e^{-2\lla t},\\
0\le\cd(S_n|_{\{t\}\times Y'_j})-\cd(\al'_j)<\const e^{-2\lla t},
\end{gather*}
where the first line holds for $0\le t\le T_j(n)$ and $j=1,\dots,r$,
the second line for $t\ge0$ and $j=1,\dots,r'$, and in both cases we
assume $\tmin(n)>C$ and $n>N(t)$.

In the following we will ignore the ends $\rpy'_{j'}$ of $X$,
ie we will pretend
that $X\gl$ is compact. If $\al'_{j'}$ is irreducible then the
argument for dealing with the end $\rpy'_{j'}$ is completely analogous to
the one given below for a neck $[-T_j,T_j]\times Y_j$, while if $\al'_{j'}$
is reducible it is simpler. (Compare the proof of
\cite[Proposition~6.3\,(ii)]{Fr10}.)

For the remainder of the proof of this lemma we will focus on one particular
neck $[-T_j(n),T_j(n)]\times Y_j$ where $1\le j\le r$. To simplify notation
we will therefore mostly omit $j$ from notation and write
$T(n),Y,\al$ etc instead of $T_j(n),Y_j,\al_j$. 

For $0\le t\le T(n)$ set
\[B_t=[-T(n)+t,T(n)-t]\times Y,\]
regarded as a subset of $\xtn$.
By the above discussion there is a constant $t_1>0$ such that when
$n$ is sufficiently large, 
$S_n$ will restrict to a genuine monopole over the band $B_{t_1+3}$ by
\cite[Lemmas~4.1,\,4.2,\,4.3]{Fr10} and will 
have small enough energy over this band for
\cite[Theorem~6.2]{Fr10} to apply. That theorem then provides a smooth
\[\ti v_{n}:B_{t_1}\to\U1\]
such that $S''_n=\ti v_{n}(S_n|_{B_{t_1}})$ is in temporal gauge and
\[\|S''_n-\uline\al\|_{\llw p\ka1(B_{t})}
\le\const e^{(\si-\lla)t},\quad t\ge t_1.\]
Writing
\[S''_n-\hat S_n =(S''_n-\ual)+(\ual-\ti S)+(\ti S-\hat S_n )\]
we get
\be{equation}\label{S''-bound}
\limsup_{n\to\infty}\|S''_n-\hat S_n \|_{\llw p\ka1(B_t)}
\le\const\left(e^{(2\si-\lla)\tau}+e^{(\si-\lla)t}\right)
\end{equation}
when $t\ge t_1$ and $\tau$ is so large that Lemma~\ref{hzeta-bounds2} applies.

To complete the proof of the lemma we interpolate between $v_n$ and $\ti v_n$
in the overlap region $\cO_\tau =X\endt \tau\cap B_{\tau-1}$. 
(This requires $\tau\ge t_1+1$.) The choice of this overlap region is somewhat
arbitrary but simplifies the exposition. Define
\[w_n=\ti v_nv_n\inv:\cO_\tau\to\U1.\]
Then
\[w_n(S'_n)=S''_n\quad\text{on $\cO_\tau$.}\]
Set $x^\pm=\ga(\pm(T-\tau))$, where $\ga=\ga_j$ is the path introduced in
the beginning of this section. If $\al$ is reducible then by multiplying
each $\ti v_n$ by a constant and redefining $w_n,S''_n$ accordingly
we can arrange that $w_n(x^+)=1$ for all $n$. These changes have no
effect on the estimates above.

Lemma~\ref{s-approx} is a consequence of the estimates~\Ref{S'-bound},
\Ref{S''-bound} together with the following sublemma (see the proof
of \cite[Proposition~6.3\,(ii)]{Fr10}.)

\be{sublemma}There is a constant $\xx{15}<\infty$ such that if
$\tau\ge \xx{15}$ then
\[\limsup_{n\to\infty}\|w_n-1\|_{L^p_2(\cO_\tau)}
\le \xx{15} e^{(2\si-\lla)\tau}.\]
\end{sublemma}

{\em Proof of sublemma:} If $\al$ is irreducible then the sublemma
follows from inequalities~\Ref{S'-bound}, \Ref{S''-bound} and
\cite[Lemmas~6.9,\,6.11]{Fr10}. (In this case the sublemma holds with
$\xx{15}e^{(\si-\lla)\tau}$ as upper bound.)

Now suppose $\al$ is reducible. We will show that
\be{equation}\label{wnx}
\limsup_{n\to\infty}|w_n(x^-)-1|\le\const e^{(2\si-\lla)\tau}
\end{equation}
for large $\tau$.
Granted this, we can prove the sublemma by applying \cite[Lemma~6.9]{Fr10} and
\cite[Lemma~6.10\,(ii)]{Fr10} to each component of $\cO_\tau$.

In the remainder of the proof of the sublemma we will omit $n$ from
subscripts.
To prove \Ref{wnx}, define intervals
\[J_0=[-T-1,-T+\tau],\quad J_1=[-T+\tau,T-\tau],\quad J_2=[T-\tau,T+1]\]
and for $k=0,1,2$ set $\ga\rp k=\ga|_{J_k}$. Let $\hol\rp k$ denote holonomy
along $\ga\rp k$ in the same sense as \Ref{hol-defn}, ie $\hol\rp k$ is
the result of replacing the domain of integration $I_j$ in that formula
with $J_k$. Define $\del\rp k\in\co$ by
\be{align*}
\hol\rp k(\hat S)&=\hol\rp k(S')(1+\del\rp k),\quad k=0,2,\\
\hol\rp 1(\hat S)&=\hol\rp 1(S'')(1+\del\rp 1)
\end{align*}
where as usual we mean holonomy with respect to the connection parts of the
configurations. For large $\tau$ the estimates
\Ref{S'-bound} and \Ref{S''-bound} give
\[|\del\rp k|\le\const e^{(2\si-\lla)\tau}\]
when $n$ is sufficiently large.

Writing $h=\prod_{k=0}^2(1+\del\rp k)$ we have
\[z=\hol(\hat S)=\prod_{k=0}^2\hol\rp k(\hat S)
=h\,\hol\rp0(S')\,\hol\rp1(S'')\,\hol\rp2(S').\]
Now, by the definition of holonomy,
\[\hol\rp1(S'')=\frac{\ti v(x^+)}{\ti v(x^-)}\,\hol\rp1(S),\]
and there are similar formulas for $\hol\rp k(S')$.
Because $w(x^+)=1$ we obtain
\[z=h\,\hol(S)\,w(x^-)\inv.\]
Setting $a=\hol(S)z\inv$ we get
\[w(x^-)-1=ah-1=(a-1)h+h-1.\]
Since by assumption $a\to1$ as $n\to\infty$, we have
\[|w(x^-)-1|\le\const\left(|a-1|+\sum_k|\del\rp k|\right)
\le\const e^{(2\si-\lla)\tau}\]
for large $n$, proving the sublemma and hence also Lemma~\ref{s-approx}.\square

This completes the proof of Proposition~\ref{prop:inj} and thereby also
the proof of Theorem~\ref{gluing-thm1}.

\subsection{Existence of maps $q$}
\label{subsec:q-existence}

Let $G\subset M_\hx$ be as in Subsection~\ref{subsec:gluing-statement1}. 
In this subsection we will show that there is always a map
$q$ as in \Ref{eqn:q} provided $\bt$ acts freely on $\oline G$ and
$K$ is sufficiently large.
It clearly suffices to prove the same with $\cb^*_\hx(K)$ in place of
$\tcb^*_\hx(K)$. 

Let $\cb,M$ denote the configuration and moduli spaces over $X$ with
the same asymptotic limits as $\cb_\hx,M_\hx$, but using the full group of
gauge transformations $\cg$ rather than $\cg_\hx$. 

Because $\cg_\hx$ acts freely on $\cc$, an element in $\cb_\hx$
has trivial stabilizer in $\bt$ if
and only if its image in $\cb$ is irreducible, ie when its spinor field 
does not vanish identically on any component of $X$.

Throughout this subsection, $K$ will be a compact codimension~$0$ submanifold
of $X$ which contains $\hx$ and intersects every component
of $X$.

\be{prop}\label{prop:q-existence}
If $\bt$ acts freely on $\oline G$ then for sufficiently large $K$ there exist
a $\bt$--invariant open neighbourhood
$V\subset\cb^*_\hx(K)$ of $R_K(\oline G)$ and a $\bt$--equivariant
smooth map $q:V\to M_\hx$
such that $q(\om|_K)=\om$ for all $\om\in G$.
\end{prop}


We first prove three lemmas. Let $H\subset M^*$ be the image of $G$. 
Because $\bt$ is compact, the projection $\cb_\hx\to\cb$ is a closed map
and therefore maps $\oline G$ to $\oline H$. Let $H_0\subset M^*$ be any
precompact open subset which contains $\oline H$ and whose closure
consists only of regular points.

\be{lemma}\label{lemma:immersion}
If $K$ is sufficiently large then $R_K:M^*\to\cb^*(K)$ restricts to
an immersion on an open neighbourhood of $\oline H_0$.
\end{lemma}

By `immersion' we mean the same as in \cite{Lang1}. Since a finite-dimensional
subspace of a Banach space is always complemented, the condition in our case
is simply that the derivative of the map be injective at every point.

\proof Fix $\om=[S]\in\oline H_0$. We will show that 
$R_K$ is an immersion at $\om$ (hence in a
neighbourhood of $\om$) when $K$ is large enough.
Since $\oline H_0$ is compact, this will prove the lemma.

Let $W\subset\llw pw1$ be a linear subspace such that the derivative at $S$
of the projection $S+W\to\cb^*$ is a linear isomorphism onto
the tangent space of $M$ at $\om$. Let $\del$ denote that derivative.
For $t\ge0$ so large that $\hx\subset X\endt t$ let
$\del_t$ be the derivative at $S$ of the natural map
$S+W\to\cb^*(X\endt t)$.
We claim that $\del_t$ is injective for $t\gg0$. 
For suppose $\{w_n\}$ is a sequence in $W$ such that $\|w_n\|_{\llw pw1}=1$ 
and $\del_{t_n}(w_n)=0$ for each $n$, where $t_n\to\infty$.
Set $K_n=X\endt{t_n}$. Then
\[w_n|_{K_n}=\ci_\Phi f_n\]
for some $f_n\in L\cg(K_n)$, where $\Phi$ is the spinor field of $S$.
After passing to a subsequence we may assume that $w_n$ converges in $\llw pw1$
to some $w\in W$ (since $W$ is finite-dimensional).
By \cite[Lemma~2.1]{Fr10} there exists for each $n$
a constant $C_n<\infty$ such that
for all $h\in L\cg(K_n)$ one has
\[\|h\|_{L^p_2}\le C_n\|\ci_\Phi h\|_{L^p_1}.\]
It follows that 
$f_n$ converges in $L^p_2$ over compact subsets of $X$ to some function $f$.
We obviously have $\ci_\Phi f=w$, hence $f\in L\cg$ and
$\del(w)=0$. But this is 
impossible, since $w$ has norm $1$. This proves the lemma.\square

\be{lemma}\label{lemma:K-emb}
If $K$ is sufficiently large then the restriction map
$H_0\to\cb^*(K)$ is a smooth embedding.\square
\end{lemma}

\proof Because of Lemma~\ref{lemma:immersion} it suffices to show
that $R_K$ is injective on $\oline H_0$ when $K$ is large.

Suppose $\om_n,\om'_n\in\oline H_0$ restrict to the same element
in $\cb(X\endt{t_n})$ for $n=1,2,\dots$, where $t_n\to\infty$. Since
$\oline H_0$ is compact we may assume, after passing to a subsequence,
that $\om_n,\om'_n$ converge in $\oline H_0$ to $\om,\om'$ respectively.
As in the proof of Lemma~\ref{choose-tiK} one finds that $\om=\om'$.
When $n$ is large then $\oline H_0\to \cb(X\endt{t_n})$ will be injective
in a neighbourhood of $\om$ by Lemma~\ref{lemma:immersion}, hence
$\om_n=\om'_n$ for $n$ sufficiently large.\square

For the time being, we will call a Banach space $E$ {\em admissible}, if
$x\mapsto\|x\|^r$ is a smooth function on $E$ for some $r>0$. (The
examples we have in mind are $L^p_k$ Sobolev spaces where $p$ is an
even integer.)

\be{lemma}\label{lemma:tub-nbhd}
Let $B$ be any second countable (smooth) Banach manifold modelled on an
admissible Banach space. Then any submanifold $Z$ of $B$ possesses
a tubular neighbourhood (in the sense of \cite{Lang1}).
\end{lemma}


\proof According to \cite[p\,96]{Lang1}, if a Banach manifold
admits partitions of unity then any {\em closed} submanifold
possesses a tubular neigbourhood. Now observe that
$Z$ is by definition locally closed, hence $C=\oline Z\setminus Z$
is closed in $B$. But then $Z$ is a closed submanifold of $B\setminus C$.
In general, any second countable, regular $T_1$--space is metrizable, hence
paracompact (see \cite{Kelley}). Because $B\setminus C$ is modelled on 
an admissible Banach space, the argument in \cite{Lang1} carries over to show
that $B\setminus C$ admits partitions of unity. Therefore, $Z$ possesses
a tubular neighbourhood in $B\setminus C$, which also serves as a tubular
neighbourhood of $Z$ in $B$.\square

{\em Proof of Proposition~\ref{prop:q-existence}:}
Choose $K$ so large that $H_0\to\cb^*(K)$ is an embedding, with image
$Z$, say.  Let $G_0$ denote the preimage of $H_0$ in $M_\hx$.

Let $\E$ be the open subset of $\cb_\hx(K)$ consisting of those elements
whose spinor does not vanish identically
on any component of $K$. Then the projection
$\pi:\E\to\cb^*(K)$ is a principal $\bt$--bundle, and restriction to $K$
defines a diffeomorphism
\[\iota:G_0\to\pi\inv Z.\]
By Lemma~\ref{lemma:tub-nbhd} there is an open
neighbourhood $U$ of $H_0$ in $\cb^*(K)$ and a smooth map
\[\rho:U\times[0,1]\to\cb^*(K)\]
such that $\rho(x,1)\in Z$ for all $x$, and $\rho(x,t)=x$ if
$x\in Z$ or $t=0$. (In other words, $\rho$ is a strong deformation retraction
of $U$ to $Z$.)
After choosing a connection in the $\bt$--bundle $\E$ we can 
then construct a $\bt$--invariant smooth retraction
\[\ti\rho:\pi\inv(U)\to\pi\inv(Z)\]
by means of holonomy along the paths
$t\mapsto\rho(t,x)$. Now set
\[q=\iota\inv\circ\ti\rho:\pi\inv U\to G_0.\square\]


\section{Applications}
\label{sec:applications}

\subsection{A model application}
\label{subsec:model-appl}

In this subsection we will show in a model case how the gluing theorem 
may be applied in combination with the compactness results of \cite{Fr10}.
Here we only consider gluing along irreducible critical points. Examples
of gluing along reducible critical points will appear in \cite{Fr12,Fr4}.
The main result of this subsection, Theorem~\ref{gluing-thm3},
encompasses both the simplest
gluing formulae for Seiberg--Witten invariants (in situations where reducibles
are not encountered) and, as we will see in the next subsection, the formula
$d\circ d=0$ for the standard Floer differential.

Recall that the Seiberg--Witten invariant of a closed $\spc$ $4$--manifold
(with $b^+_2>1$) can be defined as the number of points (counted with sign)
in the zero-set of a generic section of a certain
vector bundle over the moduli space. To
obtain a gluing formula, this vector bundle and its section
should be expressed as
the pull-back of a vector bundle $E\to\tcb^*(K)$ with section $s$, where
$K\subset X$. In the proof of Theorem~\ref{gluing-thm3} below we will see
how the section $s$ gives rise in a natural way to a map $q$ as in
Theorem~\ref{gluing-thm1}. Thus, the section $s$ is being incorporated into the
equations that the gluing map is required to solve.
(We owe this idea to \cite[p~99]{D5}.)

We will now describe the set-up for our model application.
Let $X$ be as in \cite[Subsection~1.4]{Fr10} with $r=1$ and $r'\ge0$, and
set $Y=Y_1$. In other words, we will be gluing one single pair of ends
$\R_+\times(\pm Y)$ of $X$, but $X$ may have other ends $\R_+\times Y'_j$
not involved in the gluing. We assume $X^\#$ is connected, which means
that $X$ has one or two connected components. 
For $j=1,\dots,r'$ fix a critical point $\al'_j\in\ti\fl_{Y'_j}$.
Let $\mu$ be a $2$--form and
$\prt$ a perturbation parameter for $Y$, and let $\mu'_j,\prt'_j$ be 
similar data for $Y'_j$. Let each $\prt,\prt'_j$ have small $C^1$ norm.
To simplify notation we write, for $\al,\beta\in\ti\fl_Y$,
\[M_{\al,\beta}=M(X;\al,\beta,\vec\al'),\qquad M\tu=M(\xt;\vec\al').\]
We make the following assumptions:

\be{itemize}
\item (Compactness)
At least one of the conditions (B1), (B2) of \cite{Fr10} holds for some
$\lla_j,\lla'_j>0$,
\item (Regularity)
All moduli spaces over $\ry$, $\ry'_j$, and $X$ contain only regular 
points, and
\item (No reducibles)
Given $\al_1,\al_2\in\ti\fl_Y$ and $\al'_j\in\ti\fl_{Y'_j}$,
if there exist a broken gradient line over $\ry$ from $\al_1$
to $\al_2$ and for each $j$ a broken gradient line over $\ry'_j$ from
$\al'_j$ to $\beta'_j$ then $M(X;\al_1,\al_2,\vec\al')$ contains no reducible.
(It then follows by compactness that $\mt$ contains no reducible when $T$
is large.)
\end{itemize}
The regularity condition is stronger than necessary, because there are energy
constraints on the moduli spaces that one may encounter in the situation 
to be considered, but we will not elaborate on this here.

Note that we have so far only developed a full transversality theory
in the case when $Y$ and each $Y'_j$ are rational homology spheres; in the
remaining cases the discussion here is therefore somewhat theoretical
at this time.

Let $K\subset X$ be a compact codimension~$0$ submanifold which intersects
every component of $X$. When $T\gg0$ then $K$ may also be regarded as a 
submanifold of $\xt$, and we have restriction maps
\[R_{\al,\beta}:M^*_{\al,\beta}\to\tcb^*(K),\qquad R':M\tu\to\tcb^*(K).\]
These take values in $\tcb^*(K)$ rather than just in 
$\tcb(K)$ because of the unique continuation property of harmonic spinors.

Suppose $E\to\tcb^*(K)$ is an oriented smooth real vector bundle whose rank
$d$ is equal to the (expected)
dimension of $\mt$. Choose a smooth section $s$ of $E$
such that the pull-back section $s_{\al,\beta}=R^*_{\al,\beta}s$ is
transverse to the zero-section of the pull-back bundle
$E_{\al,\beta}=R^*_{\al,\beta}E$ over $M^*_{\al,\beta}$ for
each pair $\al,\beta$. (Here the Sobolev exponent $p>4$ should be an even
integer to ensure the existence of smooth partitions of unity.)
Set $s'=(R')^*s$, which is a section of $E'=(R')^*E$.
We write $M_\al=M_{\al,\al}=M^*_{\al,\al}$ and $s_\al=s_{\al,\al}$ etc.
Let $\wm_\al$, $\wm\tu$ denote the zero-sets of $s_\al$, $s'$ respectively.
By index theory we have
\[0=\dim\,\hat M\tu=\dim\,\hat M_\al+n_\al,\]
where $n_\al=0$ if $\al$ is irreducible and $n_\al=1$ otherwise.
Thus, $\wm_\al$ is empty if $\al$ is reducible.

\be{lemma}\label{dim-count}
If $\om_n\in\wm^{(T(n))}$ for $n=1,2,\dots$, where $T(n)\to\infty$, then a
subsequence of $\om_n$ chain-converges to an element of $\wm_\al$
for some $\al\in\ti\fl^*_Y$. Moreover, if $\om_n=[S_n]$ chain-converges to
$[S]\in\wm_\al$ then there exists for each $n$ a smooth
$u_n:\xtn\to\U1$ whose restriction
to each end $\rpy'_j$ is null-homotopic and such that the sequence
$u_n(S_n)$ c-converges over $X$ to $S$.
\end{lemma} 

\proof The statement of the first sentence follows
from \cite[Theorem~1.4]{Fr10} by dimension counting. Such maps $u_n$ exist
in general for chain-convergent sequences when the $\om_n$ all have the
same asymptotic limits over the ends $\rpy'_j$, see \cite{Fr10}.\square


Let $J\subset H^1(Y;\z)$ be the subgroup consisting of elements of the form
$z|_Y$ where $z$ is an element of $H^1(X^\#;\z)$ satisfying $z|_{Y'_j}=0$
for $j=1,\dots,r'$.
This group $J$ acts on the disjoint union
\[\wm^u=\bigcup_{\al\in\ti\fl^*_Y}\wm_\al,\]
permuting the summands. 

\be{lemma}
The quotient $\wm=\wm^u/J$ is a finite set.
\end{lemma}

\proof By \cite[Theorem~1.3]{Fr10} any sequence $\om_n\in\hat M_{\al_n}$,
$n=1,2,\dots$ has a chain-convergent subsequence, and for dimensional reasons
the limit (well-defined up to gauge equivalence) must lie in some moduli space
$\hat M_\beta$. Furthermore, if $\om_n$ chain-converges to an element in
$\hat M_\beta$ then $\wm_{\al_n}$ is contained in the orbit
$J\cdot\wm_\beta$ for $n\gg0$. Therefore, each $\wm_\al$ is a finite set, and
only finitely many orbits $J\cdot\wm_\al$ are non-empty. This is 
equivalent to the statement of the lemma.\square

Note that $J$ is the largest subgroup of $H^1(Y;\z)$ which acts on
$\wm^u$ in a natural way. On the other hand, if $\wm^u$ is non-emtpy then,
since $H^1(Y;\z)$ acts freely on $\ti\fl_Y$, only subgroups $J'\subset J$
of finite index have the property that $\wm^u/J'$ is finite.

\be{lemma}\label{choose-tiK}
There is a compact codimension~$0$ submanifold $K_0\subset X$ such that
the restriction map $\wm\to\cb(K_0)$ is injective.
\end{lemma}

\proof Let $[S_j]\in \wm_{\beta_j}$, $j=1,2$, where each $S_j$ is in temporal 
gauge over the ends of $X$ (and therefore decays exponentially). 
Suppose there exists a sequence of smooth gauge transformations
$u_n:X\endt{t_n}\to\U1$ where $t_n\to\infty$, such that $u_n(S_1)=S_2$
over $X\endt{t_n}$. After passing to a subsequence we can arrange
that $u_n$ c-converges over $X$ to some gauge transformation $u$
with $u(S_1)=S_2$. If $t\gg0$ then $u|_{\{t\}\times(\pm Y)}$ will both be
homotopic to a smooth $v:Y\to\U1$ with $v(\al_1)=\al_2$. Hence
$\wm_{\al_1},\wm_{\al_2}$ lie in the same $J$--orbit, and $S_1,S_2$ represent
the same element of $\wm$ by \cite[Proposition~2.5~(iii)]{Fr10}.

Thus we can take $K_0=X\endt t$ for $t\gg0$.\square

Now fix $K_0$ as in Lemma~\ref{choose-tiK} and with $K\subset K_0$. Let
$\{b_1,\dots,b_m\}$ be the image of the restriction map
$R_{K_0}:\wm\to\tcb(K_0)$.
Choose disjoint open neighbourhoods $W_j\subset\tcb(K_0)$
of the points $b_j$. If $T\gg0$ then
\[R_{K_0}(\wm\tu)\subset\bigcup_jW_j\]
by Lemma~\ref{dim-count}. For such $T$ we get a natural map
\[g:\wm\tu\to\wm.\]
It is clear that if $g'$ is the map corresponding to a different choice
of $K_0$ and neighbourhoods $W_j$ then $g=g'$ for $T$ sufficiently large.

\be{thm}\label{gluing-thm3}
For sufficiently large $T$ the following hold:
\be{description}
\item[(i)]Every element of $\wm\tu$ is a regular point in $\mt$ and a regular
zero of $s'$.
\item[(ii)]$g$ is a bijection.
\end{description}
\end{thm}

\proof If $\wm$ is empty then, by Lemma~\ref{dim-count},
$\wm\tu$ is empty as well for $T\gg0$, and there is nothing left to prove.

We now fix $b_j$ and for the remainder of the proof omit $j$ from notation.
(Thus $b=b_j$, $W=W_j$ etc.)
We will show that for $T\gg0$ the set
\[\hat B\tu=\{\om\in\wm\tu\st\om|_{K_0}\in W\}\]
consists of precisely one element, and that this element is regular in the
sense of (i). This will prove the theorem.

By definition, $b$ is the restriction of some $\om_0\in\wm_\al$.
Choose an open neighbourhood $V\subset\tcb^*(K)$ of $b|_K$ and a smooth map
\[\pi:E|_V\to\R^d\]
which restricts to a linear isomorphism on every fibre. Choose an open
neighbourhood $V_0\subset W$ of $b$ such that $R_K(V_0)\subset V$.
Let $G_+\subset M_\al$ be a precompact open neighbourhood of $\om_0$ such
that $R_{K_0}(\oline G_+)\subset V_0$.
The assumption that $\om_0$ be a regular
zero of $s_\al$ means that the composite map
\[G_+\oset{R_K}\to V\oset{\pi\circ s}\to\R^d\]
is a local diffeomorphism at $\om_0$. We can then find an injective smooth map
\[p:\R^d\to M_\al\]
such that $p\circ\pi\circ s\circ R_K=\id$ in some open neighbourhood
$G\subset G_+$ of $\om_0$. In particular, $p\inv(\om_0)=\{0\}$ and $p$ is
a local diffeomorphism at $0$. Set 
\[q=p\circ\pi\circ s:V\to M_\al.\]
By Theorem~\ref{gluing-thm1} there is a kv-pair $(K',V')\le(K_0,V_0)$
such that if $T\gg0$ then $\qq\inv G$ consists only of regular monopoles and
\[\ff=q\circ R_K:\qq\inv G\to G\]
is a diffeomorphism. By Lemma~\ref{dim-count} one has
\[\hat B\tu=\qq\inv G\cap(s')\inv(0)=\ff\inv(\om_0)\]
for $T\gg0$. For such $T$ the set $\hat B\tu$ consists of precisely one point,
and this point is regular in the sense of (i).\square


\subsection{The Floer differential}\label{subsec:dd}

Consider the situation of \cite[Subsection~1.2]{Fr10}. Suppose a perturbation
parameter $\prt$ of small $C^1$ norm
has been chosen for which all moduli spaces $M(\al_1,\al_2)$
over $\ry$ are regular. (This is possible at least when $Y$ is a rational
homology sphere, see \cite{Fr10}.) Fix $\al_1,\al_2\in\ti\fl^*_Y$ with
\[\dim\,M(\al_1,\al_2)=2.\]
We will show that the disjoint union
\[\check M:=
\bigcup_{\beta\in\ti\fl^*_Y\setminus\{\al_1,\al_2\}}\check M(\al_1,\beta)\times
\check M(\beta,\al_2)\]
is the boundary of a compact $1$--manifold. (In other words, the standard
Floer differential $d$ satisfies $d\circ d=0$ at least with $\z/2$
coefficients.) To this end we will apply Theorem~\ref{gluing-thm3}
to the case when $X$ consists of two copies of $\ry$, say
\[X=\ry\times\{1,2\},\]
and we glue $\rpy\times\{1\}$ with $\rmy\times\{2\}$. Thus $r=1,\:r'=2$.
We take $K=K_1\cup K_2$, where $K_j=[0,1]\times Y\times\{j\}$.
In this case, $\tcb^*(K)$ is the quotient of $\cc^*(K)$ by the
null-homotopic gauge transformations. The bundle $E$ over $\tcb^*(K)$
will be the product bundle with fibre $\R^2$. To define the section $s$ of $E$,
choose $\del_1,\del_2>0$ such that $\cd$ has no critical value in the set
\[(\cd(\al_2),\cd(\al_2)+\del_2]\cup[\cd(\al_1)-\del_1,\cd(\al_1)).\]
This is possible because we assume Condition~(O1) of \cite{Fr10}.
For any configuration $S$ over $[0,1]\times Y$ set
\[s_j(S)=\int_0^1\cd(S_t)\,dt-\cd(\al_j)-(-1)^j\del_j.\]
Note that $s_j(S)$ does not change if we apply a null-homotopic gauge
transformation to $S$.
 
A configuration over $K$ consists of a pair $(S_1,S_2)$ of
configurations over $[0,1]\times Y$.
Define a smooth function $s:\tcb^*(K)\to\R^2$ (ie a section of $E$) by
\[s([S_1],[S_2])=(s_1(S_1),s_2(S_2)).\]
If $[S]$ belongs to some moduli space $M(\beta_1,\beta_2)$ over $\ry$ 
with $\beta_1\neq\beta_2$ in $\ti\fl_Y$ then $\frac d{dt}\cd(S_t)<0$
for all $t$ by choice of $\prt$. Since $J=0$, 
the natural map $\wm\to\check M$ is therefore a bijection.

Let $s'_j$ be the pull-back of $s_j$ to $\mt$. Here $\mt$ is defined using
\cite[Equation~17]{Fr10} with $\fq=0$, and so can be identified
with $M(\al_1,\al_2)$. By \cite[Theorem~1.3]{Fr10} the set
\[Z\tu=\{\om\in\mt\st s'_1(\om)=0,\quad s'_2(\om)\le0\}\]
is compact for all $T>0$. If $T\gg0$ then, by 
Theorem~\ref{gluing-thm3}, $Z\tu$ is a smooth submanifold of $\mt$, and
the composition of the two bijections
\[\prtl Z\tu=\wm\tu\oset g\to\wm\to\check M\]
yields the desired identification of $\check M$ with the boundary of a compact
$1$--manifold.

\section{Orientations}\label{sec:orientations}

In this section we discuss orientations of moduli spaces and determine
when the ungluing map of Theorem~\ref{gluing-thm1} preserves resp.\ 
reverses orientation.

We will adopt the approach to orientations of Fredholm operators
(and families of such) introduced by Benevieri--Furi \cite{Ben-Furi},
which was brought to our attention by Shuguang Wang \cite{ShuguangWang}. This
approach is more economical than the
traditional one using determinant line bundles 
(see \cite{DK,das2,Fr7}) in the sense that it produces the
orientation double cover directly. It also 
fits well in with gluing theory.

After reviewing Benevieri--Furi orientations in
Subsection~\ref{subsec:ben-furi}
we study orientations of unframed and (multi)framed moduli spaces, and the
relationship between these, in Subsection~\ref{subsec:orient-mod}. The
framings require some extra care because of reducibles. The orientation
cover $\lla\to\cb=\cb(X;\vec\al)$ is defined by the family of Fredholm
operators $\ci^*_S+D\Theta_S$ parametrized by $S\in\cc$
(cf.\ \cite[p\ 130]{D5}). Any section of
$\lla$ (which is always trivial, see Proposition~\ref{prop:llatriv2} below)
defines an orientation of the regular part of the moduli space $M^*_\hx$
for any finite, oriented subset $\hx\subset X$.
If all limits $\al_j$ are reducible then any homology orientation of $X$
determines a section of $\lla$, see Proposition~\ref{prop:llatriv} below.
To relate ungluing maps
to orientations we show that, in the notation of Subsection~\ref{subsec:surj},
any section of $\lla\to\cb$ determines a section of the
orientation cover $\lla'\to\cb'$. (Here $\cb,\cb'$ are configuration spaces
over $X,\xt$, respectively.) This is explained in
Subsection~\ref{subsec:orient-gluing} after some preparation in
Subsection~\ref{subsec:gl-Lapl} concerning framings. 
With this background material in place, the result on ungluing maps,
Theorem~\ref{thm:uglmap-ortn}, is an easy consequence of earlier estimates.
Subsection~\ref{subsec:homortngl} addresses the question of whether gluing
of orientations in the above sense is compatible with gluing
of homology orientations in the case when all limits $\al_j,\al'_j$ 
are reducible.

\subsection{Benevieri--Furi orientations of Fredholm operators}
\label{subsec:ben-furi} 

We first review Benevieri--Furi's concept of orientability of a Fredholm
operator $L:E\to F$ of index~$0$ between real Banach spaces. A {\em corrector}
of $L$ is a bounded operator $A:E\to F$ with finite dimensional image such that
$L+A$ is an isomorphism. We introduce the following equivalence relation
in the set $\cc(L)$ of correctors of $L$. Given $A,B\in\cc(L)$ set 
\[P=L+A,\quad Q=L+B.\]
Let $F_0$ be any finite dimensional subspace of $F$
containing the image of $A-B$.
Then $QP\inv$ is an automorphism of $F$
which maps $F_0$ into itself. We call $A$ and $B$ {\em equivalent}
if the map $F_0\to F_0$ induced by $QP\inv$ is orientation
preserving (which holds by convention if $F_0=0$).
This condition is independent of $F_0$. The set $\cc(L)$ is now
partitioned into
two equivalence classes (unless $E=F=0$), and we define an {\em orientation}
of $L$ to be a choice of an equivalence class, the elements of which
are then called {\em positive} correctors. A corrector which is not positive
is called {\em negative}.

Benevieri--Furi
consider $Q\inv P$ instead of $QP\inv$, but it is easy to see that 
this yields the same equivalence relation.

Note that the equivalence classes are open and closed subsets
of $\cc(L)$ with respect to the operator norm. To see this, observe that
$\cc(L)$ is open among the bounded operators $E\to F$ of finite rank.
Therefore,
if $B$ is a corrector sufficiently close to a given corrector $A$, then
$A_t=(1-t)A+tB$ is a corrector for $0\le t\le1$. Since
$\im(A_t-A)\subset\im(A)+\im(B)$, it follows by continuity that
the $A_t$ are all equivalent. In particular, $A$ and $B$ are equivalent.

If $L:E\to F$ is a Fredholm operator of arbitrary index then
for any non-negative integers $m,n$ we can form the operator
\be{equation}\label{eqn:lmn}
L_{m,n}:E\oplus\R^m\to F\oplus\R^n,\quad(x,0)\mapsto(Lx,0).
\end{equation}
If $L$ has index~$0$ then for any $m$
there is a canonical correspondence between 
orientations of $L$ and orientations of $L_{m,m}$ such that if
$A$ is a positive corrector of $L$ then $A\oplus I_{\R^m}$ is a positive
corrector of $L_{m,m}$.
If $L$ has
index $k\neq0$ then we define an orientation of $L$ to be an orientation
(in the above sense) of $L_{0,k}$ (if $k>0$) or $L_{-k,0}$ (if $k<0$).

Note that if $A$ is a corrector of $L_{m,n}$ where $n-m=\index(L)$,
and $C$ an automorphism of $\R^m$ then $A$ is equivalent to
$A\circ(I_E\oplus C)$ if and only if $\det(C)>0$, and similarly for
automorphisms of $\R^n$.

It is clear that orientations of two Fredholm operators determine 
an orientation of their direct sum. Also, a complex linear Fredholm
operator carries a canonical orientation (in this case we replace
$\R$ by $\co$ in \Ref{eqn:lmn} and the orientation is then given by
any complex linear corrector).

We now consider families of Fredholm operators. Let $\bde,\bdf$ be Banach
vector bundles over a topological space $T$, with fibres $\bde_t,\bdf_t$
over $t\in T$. (We require that these satisfy the analogues of the
vector bundle axioms VB 1-3 in \cite[pp.\ 41-2]{Lang1}
in the topological category.)
Let $L(\bde,\bdf)$ denote the Banach vector bundle over $T$ whose fibre over 
$t$ is the Banach space of bounded operators $\bde_t\to\bdf_t$.
Suppose $h$ is a (continuous) section
of $L(\bde,\bdf)$ such that $h(t):\bde_t\to\bdf_t$ is a Fredholm operator
of index~$0$ for every $t\in T$. If $\bde_t\neq0$ for every $t$
then there is a natural double cover $\ti h\to T$, the {\em orientation cover}
of $h$, whose fibre over $t$ consists of the two orientations of $h(t)$.
If $U\subset T$ is an open subset and $a$ a section of $L(\bde,\bdf)$ such
that $a(t)$ has finite rank for all $t\in U$ then $a$ defines a 
trivialization of $\ti h$ over the open set of those $t\in U$ for which 
$h(t)+a(t)$ is an isomorphism. An {\em orientation} of $h$ is by definition a
section of $\ti h$. If instead each $h(t)$ has index $k\neq0$
then we define the orientation cover $\ti h$ and orientations of $h$
by first turning $h$ into a family of index~$0$ operators as above and
then applying the definitions just given for such families. 

Wang \cite{ShuguangWang} established a 1--1 correspondence between orientations
of any family of index~$0$ Fredholm operators (between fixed Banach spaces)
and orientations of its
determinant line bundle. While we will make no use of determinant line bundles
in this paper, we need to fix our convention for passing between orientations
of a Fredholm operator $L:E\to F$ of arbitrary index and orientations
of its determinant line,
\[\det(L)=\La^{\max}\ker(L)\otimes\La^{\max}\coker(L)^*.\]
Set $n=\dim\ker(L)$ and $m=\dim\coker(L)$. Choose bounded operators
$A_1:E\to\R^n$ and $A_2:\R^m\to F$ which induce isomorphisms
\[\ti A_1:\ker(L)\to\R^n,\quad
\ti A_2:\R^m\to\coker(L).\]
Then $A=
\left(\be{smallmatrix}
0 & A_2\\
A_1 & 0
\end{smallmatrix}\right)$
is a corrector of $L_{m,n}$ which also defines an isomorphism
$J_A:\det(L)\to\R$. Moreover, two such correctors $A,B$
are equivalent if and only if $J_AJ_B\inv$ preserves orientation.
(To see this, note that after altering $A_j,B_j$ by automorphisms
of $\R^m$ or $\R^n$ as appropriate one can assume that $\ti A_j=\ti B_j$,
in which case $(1-t)A+tB$ is a corrector of $L_{m,n}$ for every $t\in\R$.)
This provides a 1--1 correspondence
between orientations of $L$ and orientations of $\det(L)$.

\subsection{Orientations of moduli spaces}
\label{subsec:orient-mod}

In the situation of \cite[Subsection~3.4]{Fr10} set
\be{equation}\label{eqn:ce-def}
\xcs=\llw pw1(X;i\La^1\oplus\bs^+),\quad\cf_1=\lw pw(X;i\R),\quad
\cf_2=\lw pw(X;i\La^+\oplus\bs^-)
\end{equation}
and consider the family of Fredholm operators
\[\del_S=\ci^*_\Phi+D\Theta_S:\xcs\to\cf:=\cf_1\oplus\cf_2\]
parametrized by $S=(A,\Phi)\in\cc(X;\vec\al)$. This family is gauge equivariant
in the sense that
\[\del_{u(S)}(us)=u\del_S(s)\]
for any $s\in\llw pw1$, $u\in\cg$, where as usual $u$ acts trivially on
differential forms and by complex multiplication on spinors. Therefore,
$\cg$ acts continuously on the orientation cover
$\ti\del$ such that the projection
$\ti\del\to\cc(X;\vec\al)$ is $\cg$--equivariant. The local slice theorem
and Lemma~\ref{lemma:cov-descend} below
then show that $\ti\del$ descends to
a double cover $\lla\to\cb(X;\vec\al)$.

(Note that in the situation of \cite[Subsection~2.4]{Fr10}, the local slice
theorem for the group $\cg$ at a reducible point $(A,0)\in\cc$ 
is easily deduced from the version of
\cite[Proposition~2.6]{Fr10} with $\hx$ consisting of one point from each
component of $X$ where $\Phi$ vanishes a.e.)

\be{lemma}\label{lemma:cov-descend}  
Let the topological group $G$ act continuously on the spaces
$Z,\ti Z$, and let $\pi:\ti Z\to Z$ be a $G$--equivariant covering map.
Suppose any point in $Z$ has arbitrarily small open neighbourhoods $U$
such that for any $z\in U$ the set
\[\{g\in G\st gz\in U\}\]
is connected. Then the natural map $\ti\pi:\ti Z/G\to Z/G$ is a covering whose
pull-back to $Z$ is canonically isomorphic to $\pi$. Pull-back defines
a 1-1 correspondence between (continuous) sections of $\ti\pi$ and
$G$--equivariant sections of $\pi$.
If in addition $G$ is connected then any section of $\pi$ is $G$--equivariant.
\end{lemma}

\proof 
Let $p:Z\to Z/G$ and $q:\ti Z\to\ti Z/G$.
If $U$ is as in the lemma and $s$ is a section of $\pi$ over $U$
then for all $z\in U$, $g\in G$ with $gz\in U$ one has
\[s(gz)=gs(z).\]
Hence $s$ descends to a section of $\ti\pi$ over $p(U)$.
If in addition $\pi\inv U$ is the disjoint union of open sets
$V_j$ each of which is mapped homeomorphically onto $U$ by $\pi$ then
$q(V_j)\cap q(V_k)=\emptyset$
when $j\neq k$. Moreover, $\cup_j q(V_j)=\ti\pi\inv p(U)$.
\square


The following proposition extends a well known result in the case when
$X$ is closed (see \cite{Morgan1,das2}). 
A related result was proved in \cite[Prop.\ 4.4.18]{nico1}.

\be{prop}\label{prop:llatriv}
If each $Y_j$ is a rational homology sphere and each $\al_j$ is reducible
then any homology orientation of $X$ canonically determines a section
of $\lla\to\cb(X;\vec\al)$.
\end{prop} 

\proof We may assume $\vec\prt=0$, since rescaling $\vec\prt$ yields
a homotopy of families $\del$. For any $(A,0)\in\cc(X;\vec\al)$
the operator $\del_{(A,0)}$ is the connected sum of the operators
$-d^*+d^+$ and $D_A$. While the homology orientation of $X$ determines
an orientation of $-d^*+d^+$ (whose cokernel we identify with
$H^0\oplus H^+$ rather than $H^+\oplus H^0$, where $H^+$ now denotes the space
of self-dual closed $L^2$ $2$--forms on $X$),
the family of complex linear operators
$D_A$ carries a natural orientation which is preserved by the action
of $\cg$. This yields a section of $\lla$ over the reducible part
$\cb^{\text{red}}\subset\cb=\cb(X;\vec\al)$. Since
the map $([A,\Phi],t)\mapsto[A,t\Phi]$, $0\le t\le1$ is a deformation
retraction of $\cb$ to $\cb^{\text{red}}$, we also obtain
a section of $\lla$ over $\cb$.\square

Returning to the situation discussed before Lemma~\ref{lemma:cov-descend},
a section of $\lla$ determines an orientation of the regular part of
the moduli space $M^*(X;\vec\al)$. As we will now explain, it also 
determines an orientation of the regular part of $M^*_\hx(X;\vec\al)$
for any finite oriented subset $\hx\subset X$. (By an orientation of $\hx$
we mean an equivalence class of orderings, two orderings being equivalent
if they differ by an even permutation.)

Let $\cw$ be the space of spinors that may occur in elements of
$\cb^*_\hx(X;\vec\al)$; more precisely, $\cw$ is the open subset of
$\Phi\rr+\llw pw1$ consisting of those elements $\Phi$ such that
$\hx\cup\text{supp}(\Phi)$ intersects every component of $X$.
For any such $\Phi$ the operator
\be{equation}\label{eqn:Delta-Phi}
\Delta_\Phi:=\ci^*_\Phi\ci_\Phi=\Delta+|\Phi|^2:L\cg\to\cf_1
\end{equation}
is injective on $L\cg_\hx$, hence 
\[\cv_\Phi:=\cf_1/\Delta_\Phi(L\cg_\hx)\]
has dimension $b:=|\hx|$ by \cite[Proposition~2.2~(i)]{Fr10}. Since
$\Phi\mapsto\Delta_\Phi$ is a smooth map from $\Phi\rr+\llw pw1$
into the space of
bounded operators $L\cg\to\cf_1$, 
the spaces $\cv_\Phi$ form a smooth vector bundle $\cv$ over $\cw$. Because
$\cw$ is simply-connected, $\cv$ is orientable.
To specify an orientation it suffices to consider those
$\Phi$ that do not vanish identically on any component of $X$. Given such a 
$\Phi$, the operator \Ref{eqn:Delta-Phi} is an isomorphism, and
we decree that a $b$--tuple $g_1,\dots,g_b\in\cf_1$ spanning
a linear complement of $\Delta_\Phi(L\cg_\hx)$ is {\em positive}
if the determinant of the matrix 
\[\left(-i\Delta_\Phi\inv(g_j)(x_k)\right)_{j,k=1,\dots,b}\]
is positive, where $(x_1,\dots,x_b)$ is any positive ordering of $\hx$.

It is natural to ask what it means for $\{g_j\}$ to be a positive
basis for $\cv_\Phi$ when $\Phi=0$. Subsection~\ref{subsec:orient-cv0}
answers this question in the case $b=1$.

For the purpose of understanding ungluing maps it is convenient to
introduce local slices for the action of $\cg_\hx$ that are defined by
compactly supported functions on $X$. Given
$S=(A,\Phi)\in\cc^*_\hx(X;\vec\al)$, choose
compactly supported smooth functions $g_j,h_j:X\to i\R$, $j=1,\dots,b$,
such that
\be{equation}\label{eqn:gjhj}
\int_Xg_jh_k=\del_{jk},
\end{equation}
where $\del_{jk}$ is the Kronecker symbol, and such that
$(g_1,\dots,g_b)$ represents a positive basis for $\cv_\Phi$. (Note that
there is a preferred choice of $h_k$, which lies in the linear span of the
$g_j$'s.)
We define the operator $\mu:\cf_1\to\cf_1$ by
\be{equation}\label{eqn:muf}
\mu f=f-\sum_{j=1}^b g_j\int_Xfh_j.
\end{equation}
Clearly, this is a projection operator whose kernel is spanned by
$g_1,\dots,g_b$.
Furthermore, $\mu$ restricts to an isomorphism
\be{equation}\label{eqn:Delta-mu}
\Delta_\Phi(L\cg_\hx)\to\im(\mu).
\end{equation}
Set $\ci^\#_\Phi=\mu\circ\ci^*_\Phi$ and
\be{equation}\label{eqn:delmusdefn}
\del_{\mu,S}:=\ci^\#_\Phi+D\Theta_S:
\xcs\to\im(\mu)\oplus\cf_2.
\end{equation}
After composing with the inverse
of \Ref{eqn:Delta-mu}, $\ci^\#_\Phi$ becomes an operator of the same kind
as considered in \cite[Subsection~3.4]{Fr10}.
Therefore, the local slice theorem
\cite[Proposition~2.6]{Fr10} applies, and if $S$ represents 
a regular point of $M^*_\hx(X;\vec\al)$ then an orientation of 
$\del_{\mu,S}$ defines an orientation of the tangent space
$T_{[S]}M^*_\hx(X;\vec\al)$.

We will now relate orientations of $\del_S$ to orientations of
$\del_{\mu,S}$. For any imaginary valued
function $f$ on $X$ let $\mu'f\in\R^b$ have coordinates
$\int_Xfh_j$, $j=1,\dots,b$. Choose non-negative integers $\ell,m$ with
$m-\ell=\index(\del_S)$ and set
\be{align*}
\nu:\cf_1\oplus\cf_2\oplus\R^m
&\oset\approx\to\im(\mu)\oplus\cf_2\oplus\R^b\oplus\R^m,\\
(x_1,x_2,y)&\mapsto(\mu x_1,x_2,\mu'x_1,y).
\end{align*}
To any corrector $\sC$ of $\dds$ we associate a corrector $\sC_\hx$ of
$\ddbs$ given by
\[\ddbs+\sC_\hx=\nu\circ(\dds+\sC).\]
The map $\sC\mapsto\sC_\hx$ clearly respects the equivalence relation
for correctors and therefore defines a 1--1 correspondence between orientations
of $\del_S$ and orientations of $\del_{\mu,S}$. Moreover, for gauge 
transformations $u$ one has
\[(u\sC u\inv)_\hx=u\sC_\hx u\inv,\]
where $u$ acts by multiplication on spinors and trivially on the other
components.

If $[S]$ is a regular point of $M^*_\hx(X;\vec\al)$ then the above defined
correspondence between orientations of $\del_S$ and orientations of
$M^*_\hx(X;\vec\al)$ at $[S]$ 
does not depend on the choice of the
$2b$--tuple $g_1,\dots,g_b,h_1,\dots,h_b$,
because the space of such $2b$--tuples supported in a given compact subset
of $X$ is
path-connected in the $C^\infty$--topology.

The relationship between the orientations of $M^*=M^*(X;\vec\al)$ and
$M^*_\hx=M^*_\hx(X;\vec\al)$ can be described explicitly as follows 
(assuming $M^*$ is regular). Let $M^{**}_\hx$ be the
part of $M^*_\hx$ that lies above $M^*$. Then $\pi:M^{**}_\hx\to M^*$
is a principal $\U1^b$--bundle whose fibres inherit orientations from
$\U1^b$. If $(v_1,\dots,v_b)$ is a positive basis for the vertical tangent
space of $M^{**}_\hx$ at $\om$ and $(v_{b+1},\dots,v_d)$ a $(d-b)$--tuple of
elements of $T_\om M^{**}_\hx$ which
maps to a positive basis for $T_{\pi(\om)}M^*$ then
$(v_1,\dots,v_d)$ is a positive basis for $T_\om M^{**}_\hx$.

\subsection{Gluing and the Laplacian}
\label{subsec:gl-Lapl}


We continue the discussion of the previous subsection, but we now consider
the situation of Subection~\ref{subsec:gluing-statement1}, so that
the ends of $X$ are labelled as in \cite[Subsection~1.4]{Fr10}
and $\hx\subset X$ is the set of start-points of the paths $\ga^\pm_j$,
$j=1,\dots,r_0$.

We define the function spaces $\xcs'$ and
$\cf'=\cf'_1\oplus\cf'_2$ over $\xt$ just as
the corresponding spaces $\xcs,\cf$ etc over $X$, replacing the weight function
$w$ by $\ka$. We also define the space $\cw'$ of spinors over $\xt$ and
the oriented vector bundle $\cv'\to\cw'$ in the same way as $\cv\to\cw$,
using the same set $\hx$.


Let $S=(A,\Phi)\in\cc$
be a configuration over $X$ such that $S-S\rr$ is compactly supported, where
$S\rr$ is the reference configuration over $X$.
For large $\tmin$ consider the glued configuration $S'=(A',\Phi')$ over
$\xt$; this is the smooth configuration over $\xt$ which agrees with
$S$ over $\itr(X\endt T)$. (This notation will also be used in later subsections.
For the time being we are only interested in the spinors.)
Let $g=(g_1,\dots,g_{b})$
be as in the previous subsection.

\be{lemma}\label{lemma:gposbas}
If $g$ represents a positive basis for $\cv_\Phi$
then $g$ also represents a positive basis for
$\cv'_{\Phi'}$ when $\tmin$ is sufficiently large.
\end{lemma}

(In this lemma it is not essential that $X$ be a $4$--manifold or that
$\bs^+$ be a spinor bundle, one could just as well use the more general set-up
in \cite[Subsection~2.1]{Fr10}, at least if $p>\dim\,X$.)

\proof There is one case where the lemma is obvious, namely when
$\Phi$ does not vanish on any component of $X$ and $g_j=\Delta_\Phi f_j$,
where $f_j$ is compactly supported.
We will prove the general case by deforming
a given set of data $\Phi,g$ to one of this special form. We begin by 
establishing a version of the lemma where `positive basis' is
replaced by `basis' and one considers compact families
of such data $\Phi,g$. To make this precise,
choose $\rho>0$ such that $\supp(g_j)\subset X\endt\rho$ for each $j$, and
let $\tmin>\rho$.
Let $\Ga'\subset L^p_1(X;\bs^+)$
and $\Ga''\subset C^\infty(X;(i\R)^{b})$ be the subspaces
consisting of those elements that vanish outside $X\endt\rho$.
Let $\Ga$ (resp.\ $\Ga_T$) be the set of pairs $(\phi,g)\in \Ga'\times \Ga''$
such that $\Phi\rr+\phi\in\cw$ (whence $\Phi'\rr+\phi\in\cw'$) and such that
$g$ represents a basis for $\cv_{\Phi\rr+\phi}$ 
(resp.\ $\cv'_{\Phi'\rr+\phi}$).

\be{sublemma}\label{sublemma:KRT}
If $K$ is any compact subset of $\Ga$ then $K\subset \Ga_T$ for $\tmin\gg0$.
\end{sublemma}

Assuming the sublemma for the moment, choose $\phi\in \Ga'$ such that 
on each component of $X$ exactly one of $\Phi,\phi$ is zero. Choose
smooth functions $f_j:X\to i\R$, $j=1,\dots,b$, which are supported
in $X\endt\rho$ and satisfy $f_j(x_k)=i\del_{jk}$,
where $(x_1,\dots,x_b)$ is a positive ordering of $\hx$. Choose
a small $\eps>0$ and set $\ti g_j=\Delta_{\Phi+\eps\phi}f_j$.
Choose a path $(\mathbf\Phi(t),\mathbf g(t))$,
$0\le t\le1$ in $\Ga$
from $(\Phi,g)$ to $(\Phi+\eps\phi,\ti g)$ such that $\mathbf g(t)=g$ for
$0\le t\le\eps$ and
\[\mathbf\Phi(t)=\be{cases}
\Phi+t\phi, & 0\le t\le\eps,\\
\Phi+\eps\phi, & \eps\le t\le1.
\end{cases}\]
Let $\mathbf\Phi'(t)$ be the glued spinor over $\xt$ obtained from
$\mathbf\Phi(t)$.
By the sublemma, if $\tmin\gg0$ then for $0\le t\le1$ one has
$(\mathbf\Phi'(t),\mathbf g(t))\in \Ga_T$. Since $\ti g$
represents a positive basis for $\cv'_{\Phi'+\eps\phi}$, it follows
by continuity that $g$ must represent a positive basis for $\cv'_{\Phi'}$.
This proves the lemma assuming the sublemma.

{\em Proof of Sublemma~\ref{sublemma:KRT}:}
Suppose to the contrary that for $n=1,2,\dots$ there are 
$(\phi(n),g(n))\in K\setminus \Ga_{T(n)}$, where $\tmin(n)\to\infty$.
We may assume
$(\phi(n),g(n))\to(\phi,g)$ in $K$. Let $V$ be the linear span of 
$g_1,\dots,g_b$ and $V_n$ the linear span of $g_1(n),\dots,g_b(n)$.
Set
\[\Phi'_n=\Phi'\rr+\phi(n).\]
By assumption there exists a non-zero
$f_n\in L\cg'_\hx$ with $\Delta_{\Phi'_n}f_n\in V_n$. Choose real numbers
$\si,\tau$ with $\rho<\si<\tau$. Since $\Delta_{\Phi'\rr}f_n=0$ outside
$X\endt\rho$, unique continuation implies that $f_n$ cannot vanish identically
on $X\endt\tau$, so we may assume that
\[\|f_n\|_{L^p_2(X\endt\tau)}=1.\]

We digress briefly to consider an injective bounded operator $J:E\to F$
between normed vector spaces and for fixed $m$ a sequence of linear maps 
$P_n:\R^m\to E$ which converges in the operator norm to an injective
linear map $P$. Then there is a constant $C<\infty$ such that
$\|e\|< C\|Je\|$ for all $e$ in a neighbourhood $U$ of $P(S^{m-1})$.
For large $n$ one must have $P_n(S^{m-1})\subset U$, hence
\[\|P_nx\|\le C\|J P_nx\|\]
for all $x\in\R^m$.

We apply this result with $m=b$,
$E=L^p_1(X)$, $F=L^p(X)$, $J$ the inclusion map,
$P_nx=\sum_jx_jg_j(n)$, and $Px=\sum_jx_jg_j$. We conclude that there
is a constant $C<\infty$ such that for sufficiently large $n$ one has
\[\|v\|_{L^p_1}\le C\|v\|_{L^p}\]
for all $v\in V_n$. For such $n$,
\be{align*}
\|f_n\|_{L^p_3(X\endt\si)}
&\le\const\left(\|\Delta f_n\|_{L^p_1}+\|f_n\|_{L^p_2}\right)\\
&\le\const\left(\|\Delta_{\Phi'_n}f_n\|_{L^p_1}
    +\||\Phi'_n|^2f_n\|_{L^p_1}+1\right)\\
&\le\const\left(\|\Delta_{\Phi'_n}f_n\|_{L^p}
    +\|\Phi'_n\|_{L^p_1}^2\|f_n\|_{L^p_1}+1\right)\\
&\le\const\left(\|\Delta f_n\|_{L^p}+1\right)\\
&\le\text{const},
\end{align*}
where except in the first term all norms are taken over $X\endt\tau$.

Let $\psi_j,\psi'_j$ be the spinor
parts of $\al_j,\al'_j$, respectively. Fix $n$ for the moment and write
$\bar\rho=T_j(n)-\rho$. Define $\bar\si$ similarly.
Over $[-\bar\rho,\bar\rho]\times Y_j$ 
we then have
\[(-\prtl_1^2+\Delta_{\psi_j})f_n=0,\]
where $\prtl_1=\frac\prtl{\prtl t}$ and
$\Delta_{\psi_j}=\Delta_{Y_j}+|\psi_j|^2$. If $h$ is any
continuous real function
on $[-\bar\rho,\bar\rho]\times Y_j$
satisfying $(-\prtl_1^2+\Delta_{\psi_j})h=0$ on
$(-\bar\rho,\bar\rho)\times Y_j$ then for any non-negative integer $k$ and
$t\in[-\bar\si+1,\bar\si-1]$ one has
\[\|h\|_{C^k([t-1,t+1]\times Y_j)}\le
\const\left(\|h\|_{L^2(\{-\bar\rho\}\times Y_j)}
+\|h\|_{L^2(\{\bar\rho\}\times Y_j)}\right),\]
where the constant is independent of $n,t$. (To see this, expand $h$ in terms
of eigenvectors of $\Delta_{\psi_j}$ and note that each coefficient function
$c$ satisfies an equation $c''=\lla^2 c$, $\lla\in\R$, which yields
$(c^2)''=2(c')^2+2(\lla c)^2\ge0$. Combine this with the usual elliptic
estimates.)
Similarly, if $h$ is any bounded
continuous function on $[\rho,\infty)\times Y'_j$
satisfying $(-\prtl_1^2+\Delta_{\psi'_j})h=0$ on $(\rho,\infty)\times Y'_j$
then for any non-negative integer $k$ and
$t\ge\si$ one has
\[\|h\|_{C^k([t,t+1]\times Y'_j)}\le\const\|h\|_{L^2(\{\rho\}\times Y'_j)},\]
for some constant independent of $t$. 

After passing to a subsequence, we may therefore assume that 
$f_n$ converges in $L^p_2$ over compact subsets
of $X$ to some function $f$, whose restriction to each end of $X$ must
be the sum of a constant function and an exponentially decaying one, the
constant function being zero if the limiting spinor over that end
($\psi_j$ or $\psi'_j$) is non-zero. In particular, $f\in L\cg_\hx$.
Furthermore,
\[\Delta_{\Phi\rr+\phi}f\in V,\qquad\|f\|_{L^p_2(X\endt\tau)}=1.\]
Since $\Delta_{\Phi\rr+\phi}$ is injective on $L\cg_\hx$ this 
contradicts the assumption that $V$ is a linear complement of
$\Delta_{\Phi\rr+\phi}(L\cg_\hx)$ in $\cf_1$.

This completes the proof of Sublemma~\ref{sublemma:KRT} and thereby
also the proof of Lemma~\ref{lemma:gposbas}.\square

\subsection{Orientations and gluing}
\label{subsec:orient-gluing}


Let $S,S'$ be as in the beginning of Subsection~\ref{subsec:gl-Lapl}.
For the time being
we will consider a map $\mu$ defined by fixed but arbitrary $b$--tuples 
$\{g_j\}$, $\{h_j\}$ of imaginary valued, compactly supported, smooth functions
on $X$ satisfying the duality relation \Ref{eqn:gjhj},
where $b$ is any non-negative integer.
We will show that an orientation of $\del_{\mu,S}$ canonically determines
an orientation of $\del_{\mu,S'}$ for large $\tmin$. Set
\[\cf_\mu=\im(\mu)\oplus\cf_2,\quad\cf'_\mu=\im(\mu)\oplus\cf'_2.\]
Choose $\tau>1$ so large that the functions $g_j,h_j$ are all supported in
$X\endt\tau$, and define
\[\xcs\endt\tau=L^p_1(X\endt\tau;i\La^1\oplus\bs^+).\]
Let $\xcs\Endt\tau$ be the subspace of $\xcs$ consisting of those elements
that are supported in $X\endt\tau$, and define $\cf_\mu\Endt\tau\subset\cf_\mu$
similarly. Set
\[\cc\Endt\tau=S\rr+\xcs\Endt{(\tau-1)}.\]
In other words, $\cc\Endt\tau$ is the set of all $L^p_{1,\loc}$ configurations
$S$ over $X$ such that $S-S\rr$ is supported in $X\endt{(\tau-1)}$.
(The $\tau-1$ is chosen here because of the non-local nature of our 
perturbations.) Suppose we
are given a bounded operator
\be{equation}\label{eqn:scdef}
\sC:\xcs\endt\tau\oplus\R^\ell\to\cf_\mu\Endt\tau\oplus\R^m
\end{equation}
with finite dimensional image, where $m-\ell=\index(\del_{\mu,S})$. Clearly,
$\sC$ induces linear maps
\[\xcs\oplus\R^\ell\to\cf_\mu\oplus\R^m,\quad
\xcs'\oplus\R^\ell\to\cf'_\mu\oplus\R^m\]
(the latter when $\tmin\ge\tau$); these will also be denoted by $\sC$. Fix an
$r_0$--tuple of paths $\ga=(\ga_1,\dots,\ga_{r_0})$ as in
Subsection~\ref{subsec:gluing-statement1}, and for
any imaginary valued $1$--form $a$ on $\xt$
let $H_\ga(a)\in\R^{r_0}$ have coordinates
\[H_{\ga_j}(a):=\int_{\ga_j}ia,\quad j=1,\dots,r_0.\]

\be{lemma}\label{lemma:CDE}
There exists a constant $C<\infty$ with the property that if
$\sC$ is any map as above and $S$ any element of
$\cc\Endt\tau$ such that
\[D:=\del_{\mu,S}+\sC:\xcs\oplus\R^\ell\to\cf_\mu\oplus\R^m\]
is invertible, then
\[E:=\del_{\mu,S'}+H_\ga+\sC:
\xcs'\oplus\R^\ell\to\cf'_\mu\oplus\R^{r_0}\oplus\R^m\]
is invertible when $\tmin>C(\|D\inv\|+1)$.
\end{lemma}

\proof Let $P_j$ be a bounded right inverse of the operator \Ref{eqn:DThj}.
As in Appendix~\ref{app:splicing}, if $\tmin>\const(\|D\inv\|+\sum_j\|P_j\|)$
then we can splice $D\inv,P_1,\dots,P_r$ to obtain a right inverse $R$ of
\[\del_{\mu,S'}+\sC:\xcs'\oplus\R^\ell\to\cf'_\mu\oplus\R^m.\]
(The present situation is slightly different from that in the appendix,
but the construction there carries over.)
Furthermore,
\[\|R\|\le\const(\|D\inv\|+\sum_j\|P_j\|).\]

Let $\ssf:\R\to\R$ be a smooth function such that
$\ssf(t)=0$ for $t\le\frac12$ and $\ssf(t)=1$ for $t\ge1$. Set
\[q_j(t)=\ssf(T_j-\tau+t)\ssf(T_j-\tau-t).\]
Thus, $q_j$ approximates the characteristic function of the interval
$[-T_j+\tau,T_j-\tau]$.
For any $c=(c_1,\dots,c_{r_0})\in\R^{r_0}$
let $\eta(c)$ be the imaginary valued $1$--form
on $\xt$ given by
\[\eta(c)=\be{cases}
0 & \text{outside $\bigcup_{j=1}^{r_0}[-T_j,T_j]\times Y_j$,}\\
-(2T_j)\inv c_jq_ji\,dt & \text{on $[-T_j,T_j]\times Y_j,\quad
j=1,\dots,r_0$.}
\end{cases}\]
For the present purposes it is convenient to rearrange summands
and regard $E$ as mapping into $(\cf'_\mu\oplus\R^m)\oplus\R^{r_0}$.
Set
\[L=R+\eta:(\cf'_\mu\oplus\R^m)\oplus\R^{r_0}\to\xcs'\oplus\R^\ell.\]
Then $EL$ takes the matrix form
\be{equation}\label{eqn:iio}
\left(\be{array}{cc}
I & 0\\
\beta & I
\end{array}\right)+o,
\end{equation}
where for large $\tmin$ one has $\|\beta\|\le\const\|R\|$ and
$\|o\|\le\const\tmin\inv$, the constants being independent of $S,T$.
As in the proof of Lemma~\ref{Xi-bounds1}
we conclude that $EL$ is invertible when $\tmin>\const(\|\beta\|+1)$,
which holds if $\tmin>\const(\|D\inv\|+1)$. Since $E$ has index~$0$,
it is invertible whenever $EL$ is surjective.\square

\be{lemma}\label{lemma:ctsc}
Suppose $\sC,\tsc$ are two maps as in \Ref{eqn:scdef} which define
correctors of $(\del_{\mu,S})_{\ell,m}$,
and let $\ga,\ti\ga$ be two $r_0$--tuples of
paths as in Subsection~\ref{subsec:gluing-statement1}.
Then for sufficiently large
$\tmin$ the following holds: $\sC$ and $\tsc$ define
equivalent correctors of $(\del_{\mu,S})_{\ell,m}$
if and only if $H_\ga+\sC$ and $H_{\ti\ga}+\tsc$ define equivalent
correctors of $(\del_{\mu,S'})_{\ell,r_0+m}$.
\end{lemma}

\proof We will use the same notation as in Lemma~\ref{lemma:CDE} and its proof.
Let $\ti D,\ti E$ be defined as $D,E$, replacing $\sC,\ga$ by
$\tsc,\ti\ga$. Observe that the image of $D-\ti D$ is contained in
$N\oplus\R^m$ for some finite dimensional subspace $N\subset\cf_\mu\Endt\tau$,
and the image of $E-\ti E$ is then contained in $(N\oplus\R^m)\oplus\R^{r_0}$
(again rearranging summands). Moreover,
\[\ti EE\inv=\ti EL(EL)\inv,\]
and $\ti EL$ has the form
\be{equation}\label{eqn:dtq}
\left(\be{array}{cc}
\ti D'R & 0\\
\beta_1 & I
\end{array}\right)+o_1,
\end{equation}
where $\|\beta_1\|$ is bounded and $\|o_1\|\to0$ as $\tmin\to\infty$,
and $\ti D'=\del_{\mu,S'}+\ti\sC$. From the
description \Ref{eqn:iio} of $EL$ we see that $\ti EE\inv$ also has the shape
\Ref{eqn:dtq}.

If $s\in\cf_\mu\Endt\tau$ and $\rho$ denotes restriction
to $X\endt\tau$, then
\be{align*}
\|\rho\ti D'Rs-\rho\ti DD\inv s\|_{\lw pw}
&\le\const\|\ti D\|\cdot\|\rho Rs-\rho D\inv s\|_{\lw pw}\\
&\le\const\|\tmin\inv\|\cdot\|\ti D\|\cdot
\left(\|D\inv\|+\sum\|P_j\|\right)\cdot\|s\|.
\end{align*}
It follows that as $\tmin\to\infty$, the determinant of the endomorphism
of $(N\oplus\R^m)\oplus\R^{r_0}$ induced by $\ti EE\inv$ approaches
the determinant of the endomorphism of $N\oplus\R^m$ induced by
$\ti DD\inv$.\square

Consider again the situation before Lemma~\ref{lemma:CDE}. Choosing a corrector
$\sC$ of $(\del_{\mu,S})_{\ell,m}$ of the kind \Ref{eqn:scdef} we obtain,
by Lemma~\ref{lemma:CDE}, an orientation of $(\del_{\mu,S'})_{\ell,r_0+m}$ for
$\tmin\gg0$, which we can extend by continuity to $\tmin>2\tau$. Now fix
$T$ with $\tmin>2\tau$, and let 
$\lla_{\mu,S}$ (resp.\ $\lla_{\mu,S'}$) denote the set consisting of the
two orientations of $\del_{\mu,S}$ (resp.\ $\del_{\mu,S'}$). From
Lemmas~\ref{lemma:CDE} and \ref{lemma:ctsc} we obtain a natural map
\[\lla_{\mu,S}\to\lla_{\mu,S'}.\]
There are two cases that we are interested in: One is
when $b=0$ (so that $\lla_{\mu,S}=\lla_S$). The other is when $b=|\hx|$
and $(g_0,\dots,g_b)$ defines a positive basis for $\cv_\Phi$.
Now letting $\mu$ refer to the second case, the preceeding discussion
yields the following diagramme of bijections:

\be{equation}\label{eqn:ortn-diagramme}
\be{array}{ccc}
\lla_S & \longrightarrow & \lla_{S'}\\
\downarrow & & \downarrow\\
\lla_{\mu,S} & \longrightarrow & \lla_{\mu,S'},
\end{array}
\end{equation}
which commutes if and only if $br_0$ is even. (The reason for the sign
is that one has to permute the summands $\R^b$ and $\R^{r_0}$ in the
target spaces of the correctors.)

Turning now to the global picture, and taking $b=0$,
let $\cc'_{[\tau]}$ denote
the set of all $L^p_{1,\loc}$ configurations $\ti S$ over $\xt$ such that
$\ti S-S'\rr$ is supported in $X\endt\tau$. Then the gluing operation
$S\mapsto S'$ defines a homeomorphism
$u:\cc\Endt\tau\to\cc'_{[\tau]}$, and Lemmas~\ref{lemma:CDE} and \ref{lemma:ctsc}
establish an isomorphism between the orientation cover
of $\cc\Endt\tau$ and the pull-back by $u$ of the orientation cover
of $\cc'_{[\tau]}$. Combining this with Proposition~\ref{prop:llatriv2} below
we see that any section of 
$\lla\to\cb$ determines a section of the orientation cover
$\lla'\to\cb'$. (Here $\cb,\cb'$ mean the same as in the beginning of
Subsection~\ref{subsec:surj}).

\be{prop}\label{prop:llatriv2}
If $X,\vec\al$ are as in \cite[Subsection~3.4]{Fr10} then
the orientation cover $\lla\to\cb(X;\vec\al)$ is trivial.
\end{prop}

\proof We may assume $X$ is connected. Let $\hx\subset X$ consist of a
single point. Let $\pi:\cb_\hx\to\cb$ be the projection, where
$\cb=\cb(X;\vec\al)$ etc.
Since $\cb$ is the quotient of $\cb_\hx$ by the
natural $\U1$ action, the local slice theorem and
Lemma~\ref{lemma:cov-descend} imply that any section of $\pi^*\lla$ descends
to a section of $\lla$. It therefore suffices to show that $\pi^*\lla$ is
trivial, or equivalently, that for any loop $\ell$ in $\cb$ that lifts
to $\cb_\hx$ the pull-back $\ell^*\lla$ is trivial.
Since $\cc\to\cb_\hx$ is a (principal) fibre bundle,
such a loop is the image of a path $z:[0,1]\to\cc$ such that 
$z(1)=u(z(0))$ for some $u\in\cg_\hx$. After
altering the loop $\ell$ by a homotopy one can arrange that $u=1$ and
$z(t)=S\rr$ (for all $t$) outside a compact subset of $X$.

(Here is one way to construct such a homotopy. For $0\le s\le1$ let
$\xi_s=\zeta_s*\ti\zeta_s$ be the composite of the two paths (both defined
for $0\le t\le1$)
\be{align*}
\zeta_s(t)&=(1-s)z(t)+sS\rr,\\
\ti\zeta_s(t)&=(1-t)\zeta_s(1)+tv_s(\zeta_s(0)),
\end{align*}
where $v_s$ is a path in $\cg$ such that $v_0=u$, and $v_1=1$ outside
a compact subset of $X$. Clearly, $\xi_0=\zeta_0$ is homotopic to $z$
relative to $\{0,1\}$. Moreover, $v_s(\xi_s(0))=\xi_s(1)$, and $\xi_1(t)=S\rr$
where $v_1=1$.)

Now let $-X$ be the Riemannian manifold $X$ with the opposite orientation
and corresponding $\spc$ structure. Starting with $X\cup(-X)$ we form, for
any $T>0$, a compact manifold $W\tu$ by gluing the
$j$'th end of $X$ with the $j$'th end of $-X$ to obtain
a neck $[-T_j,T_j]\times Y_j$. Let $S$ be any configuration over $-X$ which
agrees with $\ual_j$ over the $j$'th end, and $z_1(t)$ the configuration
over $W\tu$ obtained by gluing $S$ and $z(t)$. Then $z_1$ maps to a loop
$\ell_1$ in $\cb(W\tu)$. By Proposition~\ref{prop:llatriv2}
the orientation cover $\lla_1\to\cb(W\tu)$ is trivial. Now
Lemmas~\ref{lemma:CDE} and \ref{lemma:ctsc}
yield an isomorphism of $\z/2$--bundles
$\ell^*\lla\to\ell_1^*\lla_1$ when $\tmin$ is large,
hence $\ell^*\lla$ is trivial.\square

We now consider the situation of Theorem~\ref{gluing-thm1}. Let $b=|\hx|$.
Choose an
orientation of $\lla\to\cb$, and let $\lla'\to\cb'$ have
the associated orientation. Given an orientation of $\hx$,
this orients the regular parts of $M_\hx$ and $\mt_\hx$. In terms of 
Diagramme~\ref{eqn:ortn-diagramme} we are here using the top horizontal and
vertical maps. Because the proof of the next theorem will use
the bottom horizontal map, and the diagramme commutes if and only if
$br_0$ is even,
this causes a sign in the theorem. Note that $b$ is the cardinality of the
set $\hx$ of starting-points of the paths $\ga^\pm_j$, which of course is even 
if these points are distinct.

\be{thm}\label{thm:uglmap-ortn}
In the situation of Theorem~\ref{gluing-thm1}, if $\tmin$ is sufficiently
large then the diffeomorphism
\be{equation*}
\FF:\qq\inv G\to\ur\times G,\quad\om\mapsto(\hol(\om),\qq(\om))
\end{equation*}
is orientation preserving if $br_0$ is even and orientation reversing if
$br_0$ is odd.
\end{thm}

\proof Given $\eps=\pm1$, we will say that a map is {\em $\eps$--preserving}
if it changes orientations by the factor $\eps$. A corrector of an
oriented Fredholm operator of index~$0$ is called an
{\em $\eps$--corrector} if it is positive or negative according to the 
sign of $\eps$.

Now set $\eps=(-1)^{br_0}$.
In view of Proposition~\ref{prop:surj} it suffices to show that,
for any given point $(z,\om)\in\ur\times G$, the inverse $\FF\inv$ is
$\eps$--preserving at $(z,\om)$ when $\tmin$ is sufficiently large.

Consider the set-up in Subsection~\ref{subsec:surj}, with $\varpi:\R^d\to M_\hx$ 
orientation preserving. Let
$\pi:\cc(K)\to\tcb(K)$ be the projection. Then $f:=\varpi\inv\circ q\circ\pi$
maps a small neighbourhood of $S_0|_K$ in $\cc(K)$ to $\R^d$. Let
\[\sC:L^p_1(K;i\La^1\oplus\bs^+)\to\R^d\]
be the derivative of $f$ at $S_0$.
Let $\mu$ be as in \Ref{eqn:muf}, with
$\Phi$ the spinor part of $S_0$.
For $0\le t\le1$ set
\be{align*}
S(t)&=(1-t)S_0+tS_{0,\tau-2},\quad\del(t)=\del_{\mu,S(t)},\\
\hat S(t)&=(1-t)\hat S+tS_{0,\tau-2,T},
\quad\del'(t)=\del_{\mu,\hat S(t)}.
\end{align*}
(Thus, the $\tau$ in the proof of Theorem~\ref{gluing-thm1}
corresponds to the present $\tau-2$.)
In the following, constants will be independent of $\tau,T$.
Because $q(\om'|_K)=\om'$ for all $\om'\in G$, we see that
$\sC$ defines a positive corrector of
\be{equation*}
\del(t)_{0,d}:\xcs\to\cf_\mu\oplus\R^d
\end{equation*}
for $t=0$. Hence, if $\tau>\konst$ (for a suitable constant)
then $\sC$ will define a positive corrector
of $\del(t)_{0,d}$ for $0\le t\le1$. We want to show that if $\tau>\konst$ then
\[E_t:=\del'(t)+H_\ga+\sC:\xcs'\to
\cf'_\mu\oplus\R^{r_0}\oplus\R^d\]
is an isomorphism for $0\le t\le1$ when $\tmin\gg0$. This is a Fredholm 
operator of index~$0$, so it suffices to show that it is surjective.
As in the proof of Lemma~\ref{lemma:CDE} we can, for $\tau>\konst$ and
$\tmin>\tau+\text{const}$,
construct a right inverse $R$ of
\[\del'(1)+\sC:\xcs'\to\cf'_\mu\oplus\R^d\]
such that $\|R\|$ is bounded independently of $\tau,T$.
Set $L=R+\eta$ as in the said proof.
For notational convenience we will here regard $E_tL$ as acting on
\[\left(\cf'_\mu\oplus\R^d\right)\oplus\R^{r_0}.\]
Then there is the matrix representation
\[E_tL=
\left(\be{array}{cc}
(\del'(t))+\sC)R & \del'(t)\eta\\
H_\ga R & H_\ga\eta
\end{array}\right).\]
By Lemma~\ref{hzeta-bounds2} one has, for $\tau>\konst$,
\be{align*}
\|(\del'(t)+\sC)R-I\|&=\|(\del'(t)-\del'(1))R\|\\
&\le\const\|\hat S(t)-\hat S(1)\|_{\llw p\ka1}\le\const e^{(2\si-\lla)\tau}.
\end{align*}
Furthermore, for $\tau>\konst$,
\be{align*}
\|\del'(t)\eta\|&\le\const e^{\si \tau}\tmin\inv,\\
\|H_\ga R\|&\le\konst,\\
\|H_\ga\eta-I\|&\le(\tau-\konst)\cdot\tmin\inv.
\end{align*}
Recalling the assumption $0\le6\si<\lla$, we see that if
$\tau>\konst$ then $E_tL$ (and hence $E_t$)
will be invertible for $0\le t\le1$ when
$\tmin\gg0$.
Since $H_\ga+\sC$ is an $\eps$--corrector of
$\del'(1)_{0,r_0+d}$,
it must also be an $\eps$--corrector of $\del'(0)_{0,r_0+d}$, which
in turn is equivalent to $\FF$ being $\eps$--preserving at
$\FF\inv(z,\om)$.\square

\subsection{Homology orientations and gluing}
\label{subsec:homortngl}

In this subsection we will describe the ``gluing of orientations''
of Subsection~\ref{subsec:orient-gluing}
in terms of homology orientations in the simplest cases.
This result will be needed in \cite{Fr12}.

Let $X$ be as in \cite[Subsection~1.4]{Fr10} with $r=1$, ie.\ only one
pair of ends $\R_+\times(\pm Y)$ is being glued. Suppose $Y$ and each $Y'_j$
are rational homology spheres. We assume the glued manifold
$X^\#$ is connected, so that $X$ has at most two components. 
As in Subsection~\ref{subsec:gluing-statement1} let $\ga$ be a path in $\xt$
running once through the neck $[-T,T]\times Y$,
with starting-point $x_0$ and
end-point $x_1$. If $X$ is connected then we assume $x_0=x_1$.

As before in this section, we will denote by $H^+(X)$ the space of self-dual
closed $L^2$ $2$--forms on $X$. It is useful to observe here that orientations
of $H^+(X)$ can be specified solely in terms of the intersection form on $X$. (We made 
implicit use of this already in the definition of homology orientation in 
\cite[Subsection~1.1]{Fr10}.)
To see this, let $V$ be any real vector space with a non-degenerate symmetric
bilinear form $B:V\times V\to\R$ of signature $(m,n)$, where $m>0$ (the case
$m=0$ being trivial).
Let $\cv^+$ denote the space of all linearly independent
$m$--tuples $(v_1,\dots,v_m)$ of
elements of $V$ such that $B$ is positive definite on the linear 
span of $v_1,\dots,v_m$. Then $\cv^+$ has exactly two path-components,
and two such $m$--tuples $(v_1,\dots,v_m)$ and $(w_1,\dots,w_m)$ lie
in the same component if and only if the matrix $(B(v_j,w_k))_{j,k=1,\dots,m}$
has positive determinant. In the case when $B$ is the intersection form 
of $X$, a choice of a component of $\cv^+$ determines orientations of both
$H^+(X)$ and $H^+(\xt)$ (since the intersection
forms of $X$ and $\xt$ are canonically isomorphic).

Given the ordering of the ends $\R_+\times(\pm Y)$ of $X$ there is a natural
1--1 correspondence between homology orientations of $X$ and of $X^\#$.
In general one can specify a homology orientation by choosing ordered
bases for $H^0$, $H^1$, and $H^+$ (or equivalently, for
the dual groups). If $X$ has two components then the 
correspondence is given by replacing the basis $(x_0)$ for $H_0(X^\#)$
with the ordered basis $(x_0,x_1)$ for $H_0(X)$. If $X$ is connected
then we replace a given ordered basis $(e_1,\dots,e_\ell)$ for $H_1(X)$ 
(where $\ell=b_1(X)$) with the ordered basis $(-[\ga],e_1,\dots,e_\ell)$
for $H_1(X^\#)$. (The sign in front of $[\ga]$ is related to
a sign appearing in the formula for $\hol_j$ in \Ref{hol-defn}.)

Now fix homology orientations of $X,X^\#$ which are compatible in the
above sense.
Let $\cb,\cb'$ be the configuration spaces over $X,\xt$ with reducible
limits. According to Proposition~\ref{prop:llatriv}
the chosen homology orientations
determine an orientation $o$ of $\lla\to\cb$ and an
orientation $o'$ of $\lla'\to\cb'$. On the other hand, $\lla'$ inherits
a glued orientation $\ti o$ from $(\lla,o)$ as specified in
Subsection~\ref{subsec:orient-gluing}.

\be{prop}
\be{description}
\item[(i)]If $X$ is connected then $o'=\ti o$.
\item[(ii)]If $X$ has two components, then $o'=\ti o$ if and only if
$b_1(X)+b^+_2(X)$ is odd.
\end{description}
\end{prop}

\proof Let $S\rr=(A\rr,0)$ be a reference configuration over $X$ as in 
Subsection~\ref{subsec:gluing-statement1}
with reducible limit over each end. To simplify notation
we will now write $S,A$ instead of $S\rr,A\rr$. Let $S'=(A',0)$
be the glued reference configuration over $\xt$. Set $L_A=d^+\oplus D_A$,
so that
\[\del_S=-d^*+L_A:\xcs\to\cf.\]
Set $b_1=b_1(X)$, $b^+=b^+_2(X)$, $m=\dim\,\ker(\del_S)$,
and $\ell=\dim\,\coker(L_A)$.

Choose smooth loops $\ell_1,\dots,\ell_{b^1}$ in $X\endt0$ representing
a positive basis
for $H_1(X;\R)$ and define
\[\sB_1:L^p_1(X\endt0,i\La^1)\to\R^{b_1},\qquad
a\mapsto\left(-\int_{\ell_j}ia\right)_{j=1,\dots,b_1}.\]
Choose a bounded complex-linear map
\[\sB_2:L^p_1(X\endt0,\bs^+)\to\co^{m_2}\]
whose composition with the restriction to $X\endt0$ defines an isomorphism 
$\ker(D_A)\to\co^{m_2}$. Set
\[\sB=\sB_1+\sB_2:\xcs\endt0\to\R^{b_1}\oplus\co^{m_2}=\R^m.\]

Choose smooth imaginary-valued closed $2$--forms 
$\om_1,\dots,\om_{b^+}$ on $X$ which are supported in $X\endt0$ and such that
the cohomology classes $[-i\om_1],\dots,[-i\om_{b^+}]$ form a positive basis
of a positive subspace for the intersection form of $X$.
Then the self-dual parts $\om^+_1,\dots,\om^+_{b^+}$ map to a basis for 
$\coker(d^+)$ on both $X$ and $\xt$, which in both cases is compatible
with the chosen orientation of $H^+=\coker(d^+)^*$. Choose smooth sections 
$\om_{b^++1},\dots,\om_\ell$ of $\bs^-_X$ which are supported in $X\endt0$ and map
to a positive basis for the real vector space $\coker(D_A)$ (with its complex orientation).

The remainder of the proof deals separately with the two cases.

{\bf Case (i):} $X$ is connected.

Let $g:X\to i\R$ be a smooth function
supported in $X\endt0$ and with $\int g=i$. Then for large $T$ the orientations
$o',\ti o$ of $\del_{S'}$ are both represented by the following corrector
of $(\del_{S'})_{\ell+1,m+1}$:
\begin{align*}
\xcs'\oplus\R\oplus\R^\ell&\to\cf'\oplus\R\oplus\R^m,\\
(\xi,t,z)&\mapsto(tg+\sum_{j=1}^\ell z_j\om_j,H_\ga\xi,\sB\xi),
\end{align*}
where $H_\ga\xi$ means $H_\ga$ applied to the $1$--form part of $\xi$.

{\bf Case (ii):} $X$ has two components $X_0,X_1$, where $x_j\in X_j$.

Thus, $\rpy\subset X_0$ and $\R_+\times(-Y)\subset X_1$.
For $j=0,1$ choose a smooth function $g_j:X_j\to i\R$ supported in 
$(X_j)\endt0$ and with $\int g_j=i$. Set
\be{align*}
\sC':\xcs'\oplus\R\oplus\R^\ell\oplus\R&\to\cf'\oplus\R^m\oplus\R,\\
(\xi,t,z,t')&\mapsto(tg_0+\sum_{j=1}^\ell z_j\om_j,\sB\xi,t'),\\
\sC_\ga:\xcs'\oplus\R\oplus\R\oplus\R^\ell&\to\cf'\oplus\R\oplus\R^m,\\
(\xi,t,t',z)&\mapsto(tg_0+t'g_1+\sum_{j=1}^\ell z_j\om_j,H_\ga\xi,\sB\xi).
\end{align*}

When $T$ is large, $\sC'$ and $\sC_\ga$ are both correctors of
$(\del_{S'})_{\ell+2,m+1}$ which represent the orientations $o',\ti o$
of $\del_{S'}$, respectively. Let $\sC$ be the corrector of
$(\del_{S'})_{\ell+2,m+1}$ which has the same domain
and target spaces as $\sC_\ga$,
and which is obtained from $\sC'$ by interchanging summands as follows.
If 
\begin{gather*}
(x,y,z)\in(\xcs'\oplus\R)\oplus\R^\ell\oplus\R,\\
\sC'(x,y,z)=(u,v,w)\in\cf'\oplus\R^m\oplus\R
\end{gather*}
then $\sC(x,z,y)=(u,w,v)$. As explained in Subsection~\ref{subsec:ben-furi},
the correctors $\sC,\sC'$ are equivalent if and only if $\ell+m$ is even. Set
\[E=\del_{S'}+\sC,\quad E_\ga=\del_{S'}+\sC_\ga.\]
We have
\[(\sC_\ga-\sC)(\xi,t,t',z)=(t'g_1,H_\ga\xi-t',0),\]
so the image $N$ of $\sC_\ga-\sC$ has dimension $2$. We need to compute the
determinant of the automorphism of $N$ induced by $E_\ga E\inv$.
Let $s,s'\in\R$ and set
\[(\xi,t,t',z)=E\inv(sg_1,s',0).\]
Write $\xi=(a,\phi)\in\Ga(i\La^1\oplus\bs^+)$. Then
\be{equation}\label{eqn:dst}
-d^*a+tg_0=sg_1,\quad t'=s'.
\end{equation}
Integrating the first equation gives $t=s$. The equation 
\[E_\ga E\inv=(\sC_\ga-\sC)E\inv+I\]
now yields
\[E_\ga E\inv(sg_1,s',0)=((s+s')g_1,H_\ga\xi,0).\]
Note that $\xi$ depends on $s$ alone, so we can write $a=a(s)$.
Set $\eta=H_\ga(a(1))$. Then $E_\ga E\inv|_N$ is represented by the matrix
\[\be{pmatrix}
1 & 1\\
\eta & 0
\end{pmatrix},\]
so $\sC,\sC_\ga$ are equivalent correctors if and only if $\eta<0$.
We will show that $\eta>0$ when $T$ is large. This implies that $\sC',\sC_\ga$
are equivalent correctors if and only if $\ell+m$ is odd. Since
$D_A$ is complex linear, this will prove part~(ii) of the proposition.

Let $\{T_n\}$ be a sequence tending to $\infty$ and,
working over $X^{(T_n)}$, set
\[(\xi_n,s_n,0,z_n)=E\inv(s_n g_1,0,0),\]
or more explicitly,
\[\del_{S'}(\xi_n)+\sum_jz_{n,j}\om_j=s_n(g_1-g_0),\quad\sB\xi_n=0,\]
where $s_n>0$ is chosen such that
\[\|\xi_n\|_{L^2(X\endt2)}=1.\]
(Because the supports of $g_0$ and $g_1$ are disjoint, Equation~\Ref{eqn:dst}
shows that $a\neq0$ over $X\endt0$ when $s\neq0$.)
Write $\xi_n=(a_n,\phi_n)$. Equation~\Ref{eqn:dst} yields
\[s_n\|g_1\|^2_2=-\int\la a_n,dg_1\ra\le\|dg_1\|_2,\]
hence the sequence $s_n$ is bounded. An
analogous argument applied to the equation
\[L_A(\xi_n)+\sum_jz_{n,j}\om_j=0\]
shows that the sequence $z_n$ is bounded as well. Thus,
$\del_{S'}(\xi_n)$ is supported in $X\endt0$, and for each $k\ge0$ the
$C^k$--norm of $\del_{S'}(\xi_n)$ is bounded independently of $n$.
Now recall from \cite[Subsection~3.4]{Fr10} that over the
neck $[-T_n,T_n]\times Y$ the operator $\del_{S'}$ can be expressed in the form
$\frac\prtl{\prtl t}+P$, where
\[P=\be{pmatrix}
0 & -d^* & 0 \\
-d & *d & 0 \\
0 & 0 & -\prtl_B
\end{pmatrix}\]
for some $\spc$ connection $B$ over $Y$. Because of our non-degeneracy
assumption on the critical points, the kernel of $P$ consists of the constant
functions in $i\Om^0(Y)$. There is also a similar description of
$\del_{S'}$ over the ends $\rpy'_j$. In general, if
$(\frac\prtl{\prtl t}+P)\zeta=0$ over a band $[0,\tau]\times Y$ and $\zeta$
involves only eigenvectors of $P$ corresponding to positive eigenvalues then
for any non-negative integer $k$ and $1\le t\le\tau-2$, say, there is an
estimate
\[\|\zeta\|_{C^k([t,t+1]\times Y)}\le
\const e^{-\rho t}\|\zeta\|_{L^2(\{0\}\times Y)},\]
where $\rho$ is the smallest positive eigenvalue of $P$. This result
immediately applies to $\xi_n$ over the ends $\rpy'_j$. Over the neck
$[-T_n,T_n]\times Y$ one can write $\xi_n=\const i\,dt+\xi^+_n+\xi^-_n$
where $\xi^\pm_n$ involves only eigenvectors corresponding to
positive/negative eigenvalues of $P$. One then obtains $C^k$--estimates
on $\xi^\pm_n$ in terms of its $L^2$--norm over $\{\mp T_n\}\times Y$.
It follows that after passing to a subsequence we may assume that $\xi_n$
c-converges over $X$ to some pair $\xi=(a,\phi)$ satisfying
$\|\xi\|_{L^2(X\endt2)}=1$. Of course, we may also 
assume that the sequences $s_n,z_n$ converge, with limits $s,z$, say.
Then
\[\del_S\xi+\sum_jz_j\om_j=s(g_1-g_0),\quad \sB\xi=0.\]
Moreover,
\[\xi=\pm ci\,dt+\zeta_\pm\quad\text{on $\R_+\times(\pm Y)$,}\]
where $\zeta_\pm$ decays exponentially and
\[c=-\lim_{n\to\infty}\frac{H_\ga(a_n)}{2T_n}.\]
On the other hand, Stokes' theorem yields
\[\int_{\{-T_n\}\times Y}*a_n=-\int_{(X_0)\endt0}d^*a_n
=-\int_{(X_0)\endt0}s_ng_0=-s_ni,\]
hence
\[ci\cdot\text{Vol}(Y)=
\int_{\{0\}\times Y}*a=\lim_n\int_{\{-T_n\}\times Y}*a_n=-si.\]
Thus, $c\cdot\text{Vol}(Y)=-s\le0$.
If $c=0$ then $\xi\neq0$ would decay exponentially
on all ends of $X$ and satisfy $\del_S\xi=0$,
contradicting $\sB\xi=0$. Therefore, $c<0$, and
$H_\ga(a_n)>0$ for large $n$.

This shows that $\eta>0$ when $T$ is large.\square

\subsection{Orientation of $\cv_0$}
\label{subsec:orient-cv0}

In Subsection~\ref{subsec:orient-mod} the question arose what it means
for $\{g_j\}$ to be a positive
basis for $\cv_\Phi$ when $\Phi=0$. The following proposition
answers this question when $b=1$.
This result is not needed
elsewhere in this paper, but will be used in \cite{Fr12}.

\be{prop}\label{prop:cvbasis}
Let $X$ be as in \cite[Subsection~1.3]{Fr10}. Suppose $X$ is connected,
$b=1$, and consider the bundle $\cv\to\cw=\llw pw1(X;\bs^+)$
with all weights
$\si_j$ positive. Then $g\in C^\infty_c(X;i\R)$,
represents a positive basis for
$\cv_{0}$ if and only if $\int_Xg/i>0$.
\end{prop}

\proof It suffices to prove that $g=ih$ represents a positive basis when
$\int_Xh>0$. Let $\hx=\{x\}$.
By \cite[Proposition~2.3]{Fr10} we may assume $h\ge0$, $h(x)>0$, and
$\supp(h)\subset X\endt0$. The proposition is then a consequence of
the following lemma.

\be{lemma}\label{lemma:Dvfh}
Let $X,h$ be as above and $v$ a smooth positive function on $X$
whose restriction to each end $\rpy_j$ is the pull-back of a function
$\sv_j$ on $\R_+$. Suppose $f$ is a real function on $X$ satisfying
\[(\Delta+v)f=h,\quad df\in\llw pw1.\]
Then $f\ge0$, and $f>0$ where $h>0$.
\end{lemma}

The proof will make use of the following elementary result, whose
proof is left to the reader.

\be{sublemma}
Suppose $a,u$ are smooth real functions on $[0,\infty)$ such that
$u''=au$, $a>0$, $u(0)>0$, and $u$ is bounded. Then $u>0$ and $u'<0$.
\square
\end{sublemma}

{\em Proof of Lemma~\ref{lemma:Dvfh}:} We first study the behaviour
of $f$ on an end $\rpy_j$. We omit $j$ from notation and write $Y=Y_j$
etc. Set $\ssf=f|_{\rpy}$. By \cite[Proposition~2.1]{Fr10} the assumption
$df\in\llw pw1$ implies that $\ssf(t,\cdot)$ converges uniformly towards
a constant function $c$ as $t\to\infty$. Let $\{e_\nu\}$
be a maximal orthonormal set of eigenvectors of $\Delta_Y$ with
corresponding eigenvalues $\lla^2_\nu$. Write
\[\ssf(t,y)=\sum_\nu u_\nu(t)e_\nu(y).\]
Then
\[u''_\nu=(\lla^2_\nu+\sv)u_\nu.\]
By the sublemma, either $u_\nu=0$ or $u_\nu u'_\nu<0$.
Consequently,
\[\int_Y\ssf^2(t,y)\,dy=\sum_\nu u_\nu^2(t)\]
is a decreasing function of $t$, and
\[\max_{y\in Y}|\ssf(t,y)|\ge c\]
for all $t\ge0$. In particular, if $c<0$ then there exists a $(t,y)\in\rpy$
with $\ssf(t,y)\le c$.
Hence, if $\inf f<0$ then the infimum is attained.

Now, at any local minimum of $f$ one has
\[vf=h-\Delta f\ge h,\]
so $f\ge0$ everywhere. But then every zero of $f$ is an absolute minimum,
so $f>0$ where $h>0$. This proves the lemma and thereby also
Proposition~\ref{prop:cvbasis}.\square

\section{Parametrized moduli spaces}
\label{subsec:par-mod}


Parametrized moduli spaces appear in many different situations in gauge
theory, e.g.\ in the construction of $4$--manifold invariants \cite{DK,das2}
and Floer homology \cite{D5}, and in connection with gluing obstructions
\cite{Fr12}. A natural setting here
would involve certain fibre bundles whose fibres are $4$--manifolds. We feel,
however, that gauge theory for such bundles in general
deserves a separate treatment,
and will therefore limit ourselves, at this time, to the case 
of a product bundle over a vector space. However, we take care to set up
the theory in such a way that it would easily carry over to more
general situations.

The main goal of this section is to extend the gluing theorem and
the discussion of orientations to the parametrized case.

\subsection{Moduli spaces}
\label{subsec:mod-sp}

As in \cite[Subsection~1.3]{Fr10} let $X$ be a $\spc$ $4$--manifold with 
Riemannian metric $\bar g$ and tubular ends $\orp\times Y_j$, $j=1,\dots,r$.
Let $\sW$ be a finite-dimensional Euclidean vector space and
$\Bg=\{g_\sw\}_{\sw\in\sW}$
a smooth family of Riemannian metrics on $X$ all of which agree with $\bar g$
outside $X\endt0$. 
We then have a principal $\SO4$--bundle $\pso(\Bg)\to X\times\sW$ whose fibre
over $(x,\sw)$ consists of all positive $g_\sw$--orthonormal frames in
$T_xX$.

In the notation of \cite[Subsection~3.1]{Fr10} let $\pglc\to\pglp$ be 
the $\spc$ structure on $X$. Denote by $\pspc(\Bg)$ the pull-back of 
$\pso(\Bg)$ under the projection $\pglc\times\sW\to\pglp\times\sW$. Then
$\pspc(\Bg)$ is a principal $\Spinc4$--bundle over $X\times\sW$. 

For $j=1,\dots,r$ let $\al_j\in\cc(Y_j)$ be a non-degenerate smooth
monopole. Let $\cc(g_\sw)$ denote the $\llw pw1$ configuration space over $X$
for the metric $g_\sw$ and limits $\al_j$, where $p,w$ are as in 
\cite[Subsection~3.4]{Fr10}. We will provide the disjoint union
\[\cc(\Bg)=\bigcup_{\sw\in\sW}\cc(g_\sw)\times\{\sw\}\]
with a natural structure of a (trivial) smooth fibre bundle over $\sW$.
Let
\be{equation}\label{eqn:v-iso}
\sv:\pspc(\Bg)\to\pspc(g_0)\times\sW
\end{equation}
be any isomorphism of $\Spinc4$--bundles which covers the identity on
$X\times\sW$ and which outside $X\endt1\times\sW$ is
given by the identification $\pspc(g_\sw)=\pspc(g_0)$.
There is then an induced isomorphism of $\SO4$--bundles
\[\pso(\Bg)\to\pso(g_0)\times\sW,\]
since these are quotients of the 
corresponding $\Spinc4$--bundles by the $\U1$--action.
Such an isomorphism $\sv$
can be constructed by means of the holonomy along rays of the form
$\{x\}\times\orp\sw$ where $(x,\sw)\in X\times\sW$, with respect to any
connection in $\pspc(\Bg)$ which outside $X\endt1\times\sW$ is the pull-back
of a connection in $\pspc(g_0)$. Then $\sv$ induces a 
$\cg=\cg(X;\vec\al)$--equivariant diffeomorphism
\be{equation}\label{viso}
\cc(g_\sw)\to\cc(g_0)
\end{equation}
for each $\sw$, where the map on the spin connections is obtained by
identifying these with connections in the respective determinant line bundles
and applying the isomorphism between these bundles induced by $\sv$. Putting
together the maps \Ref{viso} for all $\sw$ yields a bijection
\[\sv_*:\cc(\Bg)\to\cc(g_0)\times\sW.\]
If $\ti\sv$ is another isomorphism as in \Ref{eqn:v-iso}
then $\sv_*(\ti\sv_*)\inv$
is smooth, hence we have obtained the desired structure on
$\cc(\Bg)$. Furthermore, because of the gauge equivariance
of $\sv_*$ we also get a similar smooth fibre bundle structure on
\be{equation}\label{cbpar-defn}
\cb^*_\hx(\Bg)=
\bigcup_{\sw\in\sW}\cb^*_\hx(g_\sw)\times\{\sw\}
\end{equation}
for any finite subset $\hx\subset X$. The image of
$(S,\sw)\in\cc(\Bg)$ in $\cb(\Bg)$ will be denoted $[S,\sw]$.

We consider the natural smooth action of $\bt$ on $\cb^*_\hx(\Bg)$ where an element of $\bt$
maps each fibre $\cb_\hx(g_\sw)$ into itself in the standard way. (There is another version
of the gluing theorem where $\bt$ acts non-trivially on $\sW$, see below.)

The principal bundle $\pspc(\Bg)$ also
gives rise to Banach vector bundles $\xcs(\Bg),\cf(\Bg),\cf_2(\Bg)$ over $\sW$ whose 
fibres over $\sw\in\sW$ are the spaces $\xcs(g_\sw),\cf(g_\sw),\cf_2(g_\sw)$ resp.\ defined as
in \Ref{eqn:ce-def} using the metric $g_\sw$ on $X$. 

Let $\dth:\cc(\Bg)\to\cf_2(\Bg)$ be the fibre-preserving monopole
map whose effect on the fibre
over $\sw\in\sW$ is the left hand side of \cite[Equation (20)]{Fr10}, 
interpreted in terms of the metric $g_\sw$. If we conjugate $\dth$
by the appropriate diffeomorphisms induced by $\sv$ then we obtain the smooth
$\cg$--equivariant map
\begin{gather*}
\dth_{\sv}:\cc(g_0)\times\sW\to\cf_2(g_0),\\
(A,\Phi,\sw)\mapsto
\left((\sv_\sw(\hatf_A+\mm(A,\Phi)))^+-Q(\Phi)\,,\,
\sum_j\sv_\sw(e_j)\cdot\nabla^{A+\sa_\sw}_{e_j}(\Phi)\right)
\end{gather*}
where the perturbation $\mm$ is smooth, hence $\dth$ is smooth.
Here $\sv_\sw$ denotes the isomorphism that $\sv$ induces from
the Clifford bundle of $(X,g_\sw)$ to the Clifford bundle of $(X,g_0)$,
and $\{e_j\}$ is a local $g_\sw$--orthonormal frame on $X$. Finally,
if we temporarily let $\nabla^{(\sw)}$ denote the
$g_\sw$--Riemannian connection
in the tangent bundle of $X$ then
\[\sa_\sw=\sv_\sw(\nabla^{(\sw)})-\nabla^{(0)}.\]
Note that $\sa_\sw$ is supported in $X\endt1$.

In situations involving parametrized moduli spaces there will often
be an additional perturbation which affects the equations only over some
compact part of $X$. For the gluing theory one can consider quite generally
perturbations given by an isomorphism $\sv$ and a smooth $\cg$--equivariant map
\[\fo:\cc(X\endt\ft,g_0)\times\sW\to(\cf_2)\Endt\ft(g_0)\]
for some $\ft\ge0$, using notation introduced in
Subsection~\ref{subsec:orient-gluing}. We require that the derivative of
$\fo$ at any point be a compact operator. Let
\be{equation}\label{eqn:Theta-cc}
\Theta:\cc(\Bg)\to\cf_2(\Bg)
\end{equation}
be the map corresponding to $\Theta_\sv:=\dth_\sv+\fo$. We define the
parametrized moduli space $M_\hx(\Bg)$ to be the image of $\Theta\inv(0)$
in $\cb_\hx(\Bg)$. By construction, $\sv_*$ induces a homeomorphism
\[M_\hx(\Bg)\oset\approx\to\Theta_{\sv}\inv(0)/\cg_\hx.\]
A point in $M_\hx(\Bg)$ is called {\em regular} if the
corresponding zeros of $\Theta_{\sv}$ are regular
(a regular zero being one where the derivative of $\Theta_{\sv}$ is
surjective). This notion is independent of $\sv$.
By the local slice theorem, the set of regular points in
$M^*_\hx(\Bg)$
is a smooth submanifold of $\cb^*_\hx(\Bg)$.

\subsection{Orientations}

Fix orientations of the vector space $\sW$ and of the set $\hx$.
For any $S\in\cc(g_\sw)$ let
\[\del_{S,\sw}:\xcs(g_\sw)\to\cf(g_\sw)\]
be the Fredholm operator $\del_S$ defined in terms of the metric $g_\sw$,
now using
the perturbed monopole map \Ref{eqn:Theta-cc}.
The orientation cover
of this family descends to a double cover $\lla(\Bg)\to\cb(\Bg)$.
(Note that the perturbation $\fo$ can
be scaled down, so that an orientation of $\lla(\Bg)$ for $\fo=0$ determines
an orientation for any other $\fo$.)
Clearly, any section of $\lla(\Bg)$ over $\cb(g_0)$
extends uniquely to all of $\cb(\Bg)$.
On the other hand, a section of $\lla(\Bg)$
determines an orientation of the regular part of $M^*_\hx(\Bg)$,
as we will now explain. 

Let $T^v\cc(\Bg)\subset T\cc(\Bg)$ be the subbundle of vertical
tangent vectors. We can identify
$T_{(S,\sw)}^v\cc(\Bg)=\xcs(g_\sw)$.
A choice of an isomorphism $\sv_1$ as in \Ref{eqn:v-iso}
determines a 
bundle homomorphism
\[P_1:T\cc(\Bg)\to T^v\cc(\Bg)\]
which is the identity on vertical tangent vectors. This yields a splitting
\[T_{(S,\sw)}\cc(\Bg)=\xcs(g_\sw)\oplus\sW\]
into vertical and horizontal vectors (the latter making up the kernel
of $P_1$ and being identified with $\sW$ through the projection).

In general, a connection in a vector bundle $E\to\sW$ determines for every
element $u$ of a fibre $E_\sw$ a linear map $T_uE\to E_\sw$, namely the projection onto the
vertical part of the tangent space. Moreover, if $u=0$ then this projection is
independent of the connection. Together these projections form a smooth map
$TE\to E$. Let
\[P_2:T\cf_2(\Bg)\to\cf_2(\Bg)\]
be such a map for $E=\cf_2(\Bg)$ determined by some isomorphism $\sv_2$.

Now let
\[\ci^*:T^v\cc(\Bg)\to\cf_1\]
be the map which sends
$s\in T^v_{(S,\sw)}\cc(\Bg)$ to $\ci_S^*(s)$,
where the ${}^*$ refers to the metric $g_\sw$. Set
\[\uline\del:=\ci^*\circ P_1+P_2\circ D\Theta:T\cc(\Bg)\to\cf(\Bg),\]
where $D\Theta$ is the derivative of the map \Ref{eqn:Theta-cc}.
By restriction of $\udel$ we obtain bounded operators
\[\udel_{S,\sw}:T_{(S,\sw)}\cc(\Bg)\to\cf(g_\sw).\]
Since the restriction of $\udel_{S,\sw}$ to the vertical tangent space
$\xcs(g_\sw)$ is equal to the Fredholm operator $\del_{S,\sw}$, we conclude
that $\udel_{S,\sw}$ is also Fredholm, and
\[\ind(\udel_{S,\sw})=\ind(\del_{S,\sw})+\sd,\]
where $\sd=\dim\,\sW$.

Choose non-negative integers $\ell,m$ with $\ind(\del_{S,\sw})=m-\ell$
and an orientation preserving linear isomorphism $h:\sW\to\R^\sd$.
If $\sC$ is any corrector of $(\del_{S,\sw})_{\ell,m}$ then
\[\udel_{S,\sw}+h+\sC:\xcs(g_\sw)\oplus\sW\oplus\R^\ell
\to\cf(g_\sw)\oplus\R^\sd\oplus\R^m\]
is an injective Fredholm operator of index~$0$, hence an isomorphism.
The map $\sC\mapsto h+\sC$ respects the equivalence relation for correctors
and therefore defines a 1-1 correspondence between orientations of 
$\del_{S,\sw}$ and orientations of $\udel_{S,\sw}$.

Fix $S\in\cc^*_\hx(g_\sw)$, choose a map $\mu$ as in \Ref{eqn:muf}, 
and let
\[\udel_{\mu,S,\sw}:=\xcs(g_\sw)\oplus\sW\to\cf_\mu(g_\sw)\]
be the operator obtained from $\udel_{S,\sw}$ by replacing $\ci^*_S$
by $\mu\circ\ci^*_S$ (cf.\ \Ref{eqn:delmusdefn}). Just as in the 
unparametrized case one establishes a 1-1 correspondence between orientations
of $\udel_{S,\sw}$ and orientations of $\udel_{\mu,S,\sw}$.

Now suppose $[S,\sw]\in M^*_\hx(\Bg)$. Working in the trivialization $\sv_1$
and using the local slice theorem for the metric $g_\sw$ one finds that
$[S,\sw]$ is a regular point of $M^*_\hx(\Bg)$
if and only if $\udel_{\mu,S,\sw}$ is surjective,
and in that case the projection $\cc^*_\hx(\Bg)\to\cb^*_\hx(\Bg)$
induces an isomorphism
\[\ker(\udel_{\mu,S,\sw})\oset\approx\to T_{[S,\sw]}M^*_\hx(\Bg).\]
This establishes a 1-1 correspondence between orientations of
$\del_{S,\sw}$ and orientations of $T_{[S,\sw]}M^*_\hx(\Bg)$.
This correspondence
is obviously independent of $P_2$, and it is independent of $P_1$ because
the set of such operators form an affine space. It is also independent of
$\mu$ for reasons explained earlier.

This associates to any orientation
of $\lla(\Bg)$ an orientation of the regular part of $M^*_\hx(\Bg)$.

\subsection{The gluing theorem}
\label{subsec:gluing-statement2}


We continue the discussion of the previous subsection, but we now specialize
to the case when the ends of $X$ are $\R_+\times(\pm Y_j)$, $j=1,\dots,r$
and $\R_+\times Y'_j$, $j=1,\dots,r'$, with non-degenerate limits $\al_j$ over
$\R_+\times(\pm Y_j)$ and $\al'_j$ over $\R_+\times Y'_j$, as in
Subsection~\ref{subsec:gluing-statement1}. 
Let the paths $\ga^\pm_j,\ga_j$ and $\hx\subset X$ be as in that
subsection. The family of metrics $\Bg$ on $X$ defines, in a natural way,
a smooth family of metrics $\{g(T,\sw)\}_{\sw\in\sW}$ on $\xt$ for any $T$.
We retain our previous notation for configuration and moduli spaces over $X$,
whereas those over $\xt$
will be denoted $\cc'(\Bg),\cb'_\hx(\Bg),M\tu_\hx(\Bg)$ etc.
Fix an isomorphism $\sv$ as in \Ref{eqn:v-iso}. 

We first discuss gluing of orientations.
The isomorphism $\sv$ defines a 
corresponding isomorphism over $\xt$ and operators $P_1,P_2$ over
both $X$ and $\xt$. We then get families of Fredholm operators
$\udel,\udel'$ parametrized by
$\cc(\Bg),\cc'(\Bg)$ resp. The procedure in
Subsection~\ref{subsec:orient-gluing}
for gluing orientations carries over to this situation
and yields a 1-1 correspondence between orientations of $\udel$ and
orientations of $\udel'$. Given $S\in\cc\Endt0(g_\sw)$, if $\lla_{S,\sw}$
resp.\ $\uline\lla_{S,\sw}$ denote the set of orientations of
$\del_{S,\sw}$ resp.\ $\udel_{S,\sw}$, and similarly for the glued 
configuration $S'\in\cc'(g_\sw)$, then the diagramme of bijections
\be{equation*}
\be{array}{ccc}
\lla_{S,\sw} & \longrightarrow & \lla_{S',\sw}\\
\downarrow & & \downarrow\\
\uline\lla_{S,\sw} & \longrightarrow & \uline\lla_{S',\sw}
\end{array}
\end{equation*}
commutes if and only if $\sd r_0$ is even.

Now fix an orientation of $\lla(\Bg)\to\cb_\hx(\Bg)$ and let
$\lla'(\Bg)\to\cb'_\hx(\Bg)$
have the glued orientation. These orientations determine orientations 
of the regular parts of the moduli spaces $M_\hx(\Bg)$ and $\mt_\hx(\Bg)$,
respectively, as specified in the previous subsection.

As before, a choice of reference
configuration in $\cc(g_0)$ gives rise to a glued reference
configuration in $\cc'(g_0)$
and a holonomy map
\[\cb'_\hx(g_0)\to\ur.\]
Composing this with the map $M\tu_\hx(\Bg)\to\cb'_\hx(g_0)$ defined by the chosen
isomorphism $\sv$ yields a holonomy map
\[\hol:M\tu_\hx(\Bg)\to\ur.\] 

Fix an open $\bt$--invariant subset $G\subset M_\hx(\Bg)$ whose closure
is compact and contains only regular points.

By a kv-pair we mean as before a pair $(K,V)$, where
$K\subset X$ is a compact codimension~$0$ submanifold which contains
$\hx$ and intersects every component of $X$, and
$V$ is an open $\bt$--invariant neighbourhood
of $R_K(\oline G)$ in
\[\tcb_\hx(K,\Bg)=\bigcup_{\sw\in\sW}\tcb_\hx(K,g_\sw)\times\{\sw\}.\]

Now fix a kv-pair $(K,V)$ satisfying similar additional assumptions as before:
firstly, that $V\subset\tcb^*_\hx(K,\Bg)$; secondly, that
if $X_e$ is any component of $X$ which contains a point from $\hx$
then $X_e\cap K$ is connected.

Suppose
\[q:V\to M_\hx(\Bg)\]
is a smooth $\bt$--equivariant map such that $q(\om|_K)=\om$ for all
$\om\in G$.
(We do not require that $q$ commute with the projections to $\sW$.)
Choose $\lla_j,\lla'_j>0$. Let `admissibility of $\vec\al'$' be defined
in terms of the parametrized moduli spaces $\mt_\hx(\Bg)$ (see
\cite[Definition~7.3]{Fr10}).

\be{thm}\label{gluing-thm2}
Theorem~\ref{gluing-thm1} holds in the present situation if one
replaces $M_\hx$ and $M\tu_\hx$ by $M_\hx(\Bg)$ and $M\tu_\hx(\Bg)$,
respectively. Moreover, the diffeomorphism $\FF$ defined as in 
Theorem~\ref{thm:uglmap-ortn}
preserves or reverses orientations according as to whether $(b+\sd)r_0$
is even or odd.
\end{thm}

\proof The proofs carry over without any substantial changes.\square

There is another version of the theorem (which will be used in \cite{Fr12})
where the family of metrics $\Bg$ is constant (ie
$g_\sw=g_0$ for every $\sw$) and $\bt$ acts smoothly on the manifold $\sW$. One then has
a product action of $\bt$ on
\[M_\hx(\Bg)=M_\hx\times\sW,\]
and the theorem holds in this setting as well. In fact, the action of $\bt$ affects the
proof in only one way, namely the requirement that $\ti K$ be $\bt$--invariant. To obtain
this, let $\text{dist}$ be a $\bt$--invariant metric on the set $\sW$ (arising for instance
from a $\bt$--invariant Riemannian metric) and replace the definition of $d_m$ in
\Ref{eqn:dm-defn} by
\[d_m((S,\sw),(\bar S,\bar{\sw}))=
\int_{X\endt m}|\bar S-S|^p+|\nabla_{\bar A}(\bar S-S)|^p + \text{dist}(\sw,\bar{\sw}).\]
Then $V'_m$ will be $\bt$--invariant.

\subsection{Compactness}

In contrast to gluing theory, compactness requires more specific
knowledge of the perturbation $\fo$, so we will here take $\fo=0$.
We observe that the notion of chain-convergence 
has a natural generalization to the 
parametrized situation, and that the compactness theorem 
\cite[Theorem~1.4]{Fr10} carries over to sequences
\[[A_n,\Phi_n,\sw_n]\in M_\hx(\xtn,g(T(n),\sw_n);\vec\al'_n)\]
provided the sequence $\sw_n$ is
bounded (and similarly for \cite[Theorem~1.3]{Fr10}).
The only new ingredient in the proof is the following simple fact:
Suppose $B$ is a Banach space, $E,F$ vector bundles over a compact manifold,
$L,L':\Ga(E)\to\Ga(F)$ differential operators of order $d$, and 
$K:\Ga(E)\to B$ a linear operator.
If $L$ satisfies an inequality 
\[\|f\|_{L^p_k}\le C\left(\|Lf\|_{L^p_{k-d}}+\|Kf\|_B\right)\]
and $L,L'$ are sufficiently close in the sense that 
\[\|(L-L')f\|_{L^p_{k-d}}\le\eps\|f\|_{L^p_k}\]
for some constant $\eps>0$ with $\eps C<1$, then $L'$ obeys the inequality
\[\|f\|_{L^p_k}\le
(1-\eps C)\inv C\left(\|L'f\|_{L^p_{k-d}}+\|Kf\|_B\right).\]

\appendix

\section{Splicing left or right inverses}
\label{app:splicing}

Let $X$ be a Riemannian manifold with tubular ends as in
\cite[Subsection~1.4]{Fr10} but of arbitrary dimension. Let $E\to X$ be a 
vector bundle which over each end $\R_+\times(\pm Y_j)$
(resp.\ $\rpy'_j$) is isomorphic (by a fixed isomorphism) to the pull-back
of a bundle $E_j\to Y_j$ (resp.\ $E'_j\to Y'_j$). Let $F\to X$ be another
bundle of the same kind. Let $D:\Ga(E)\to\Ga(F)$ be a differential operator
of order $d\ge1$ which is translationary invariant over each end and such that
for each $j$ the restrictions of $D$ to $\rpy_j$ and $\R_+\times(-Y_j)$
agree in the obvious sense. The operator $D$ gives rise to a glued
differential operator $D':\Ga(E')\to\Ga(F')$ over $\xt$, where
$E',F'$ are the bundles over $\xt$ formed from $E,F$ resp.
Let $k,\ell,m$ be non-negative integers and $1\le p<\infty$. Let
$L^p_k(X;F)\endt0$ denote the subspace of $L^p_k(X;F)$ consisting of those
elements that vanish a.e.\ outside $X\endt0$. We can clearly also identify
$L^p_k(X;F)\endt0$ with a subspace of $L^p_k(\xt;F')$.
Let 
\[V:L^p_{k+d}(X\endt0;E)\oplus\R^\ell\to L^p_k(X;F)\endt0\oplus\R^m\]
be a bounded operator and set
\be{align*}
P=D+V:L^p_{k+d}(X;E)\oplus\R^\ell&\to L^p_k(X;F)\oplus\R^m\\
(s,x)&\mapsto(Ds,0)+V(s|_{X\endt0},x).
\end{align*}
Define the operator $P'=D'+V$ over $\xt$ similarly.

\be{prop}
If $P$ has a bounded left (resp.\ right) inverse $Q$ then for $\tmin>C_1\|Q\|$
the operator $P'$ has a bounded left (resp.\ right) inverse $Q'$ with
$\|Q'\|<C_2\|Q\|$. Here the constants $C_1,C_2<\infty$ depend on the
restriction
of $D$ to the ends $\R_+\times(\pm Y_j)$ but are otherwise independent of $P$.
\end{prop}

For left inverses this was proved in a special case in
\cite[Lemma~5.4]{Fr10}, and the general case is not very different.
However, we would like to have the explicit expression for the right inverse
on record, since this is used both in Subsection~\ref{subsec:surj} and in 
Subsection~\ref{subsec:orient-gluing}.


\proof Choose smooth functions $f_1,f_2:\R\to\R$
such that $(f_1(t))^2+(f_2(1-t))^2=1$ for all $t$, and $f_k(t)=1$ for
$t\le\frac13$, $k=1,2$.
Define $\beta:X\to\R$ by
\[\beta=\be{cases}
f_1(t/(2T_j)) & \text{on each end $\rpy_j$,}\\
f_2(t/(2T_j)) & \text{on each end $\R_+\times(-Y_j)$,}\\
1 & \text{elsewhere,}
\end{cases}\]
where $t$ is the first coordinate on $\R_+\times(\pm Y_j)$. If $s'$ is a 
section of $F'$ let the section $\obeta(s')$ of $F$ be the result of pulling
$s'$ back to $X^{\{T\}}$ by means of $\pi^{\{T\}}$, multiplying by $\beta$,
and then extending trivially to all of $X$. (The notation $\pi^{\{T\}}$ was
introduced in \cite[Subsection~1.4]{Fr10}.) If $x\in\R^m$ set
$\obeta(s',x)=(\obeta(s'),x)$. For any section $s$ of $E$ we define
a section $\ubeta(s)$ of $E'$ as follows when $\tmin\ge3/2$. Outside
$[-T_j+1,T_j-1]\times Y_j$ we set $\ubeta(s)=s$. Over 
$[-T_j,T_j]\times Y_j$ let $\ubeta(s)$ be the sum of the restrictions
of the product $\beta s$ to $[0,2T_j]\times Y_j$ and $[0,2T_j]\times(-Y_j)$, 
identifying both these bands with $[-T_j,T_j]\times Y_j$ by means of
the projection $\pi^{\{T\}}:X^{\{T\}}\to\xt$. If $x\in\R^\ell$ set
$\ubeta(s,x)=(\ubeta(s),x)$. Note that
\[\ubeta\obeta=I.\]
Now suppose $Q$ is a left or right inverse
of $P$. Define
\[R'=\ubeta Q\obeta:L^p_k(\xt;F')\oplus\R^m
\to L^p_{k+d}(\xt;E')\oplus\R^\ell.\]
If $QP=I$ then a simple calculation yields
\[\|R'P'-I\|\le C\tmin\inv\|Q\|.\]
Therefore, if $\tmin>C\|Q\|$ then $R'P'$ is invertible and
$Q'=(R'P')\inv R'$ is a left inverse of $P'$. Similarly, if $PQ=I$ then
\[\|P'R'-I\|\le C\tmin\inv\|Q\|,\]
hence $Q'=R'(P'R')\inv$ is a right inverse of $P'$ when $\tmin>C\|Q\|$.
In both cases the constant $C$ depends on the restriction
of $D$ to the ends $\R_+\times(\pm Y_j)$ but is otherwise independent of $P$.
As for the bound on $\|Q'\|$, see the proof of Lemma~\ref{Xi-bounds1}.
\square

From the proposition one easily deduces the following version of
the addition formula for the index, which was proved for
first order operators in \cite{D5}.

\be{cor}\label{cor:add-formula}
If
\[D:L^p_{k+d}(X;E)\to L^p_k(X;F)\]
is Fredholm, then for sufficiently large
$\tmin$,
\[D':L^p_{k+d}(\xt;E')\to L^p_k(\xt;F')\]
is Fredholm with $\ind(D')=\ind(D)$.
\end{cor}

\noindent\textsc{Fakult\"at f\"ur Mathematik, Universit\"at Bielefeld,\\
Postfach 100131, D-33501 Bielefeld, Germany.}\\
\\
E-mail:\ froyshov@math.uni-bielefeld.de


\begin{thebibliography}{10}

\bibitem{Ben-Furi}
P.~Benevieri and M.~Furi.
\newblock A simple notion of orientability for {Fredholm} maps of index zero
  between {Banach} manifolds and degree theory.
\newblock {\em Ann. Sci. Math. Qu\'ebec}, 22:131--148, 1998.

\bibitem{wchen1}
W.~Chen.
\newblock Casson's invariant and {Seiberg--Witten} gauge theory.
\newblock {\em Turkish J. Math.}, 21:61--81, 1997.

\bibitem{D1}
S.~K. Donaldson.
\newblock An application of gauge theory to four dimensional topology.
\newblock {\em J.~Differential Geometry}, 18:279--315, 1983.

\bibitem{D7}
S.~K. Donaldson.
\newblock Connections, cohomology and the intersection forms of $4$--manifolds.
\newblock {\em J.~Differential Geometry}, 24:275--341, 1986.

\bibitem{D5}
S.~K. Donaldson.
\newblock {\em Floer Homology Groups in {Yang--Mills} Theory}.
\newblock Cambridge University Press, 2002.

\bibitem{DK}
S.~K. Donaldson and P.~B. Kronheimer.
\newblock {\em The Geometry of Four-Manifolds}.
\newblock Oxford University Press, 1990.

\bibitem{Feehan-Leness2}
P.~M.~N. Feehan and T.~G. Leness.
\newblock {A general $SO(3)$--monopole cobordism formula relating Donaldson and
  Seiberg--Witten invariants}.
\newblock math.DG/0203047.

\bibitem{F1}
A.~Floer.
\newblock An instanton invariant for 3--manifolds.
\newblock {\em Comm. Math. Phys.}, 118:215--240, 1988.

\bibitem{FU}
D.~S. Freed and K.~K. Uhlenbeck.
\newblock {\em Instantons and Four-Manifolds}.
\newblock MSRI Publications. Springer-Verlag, second edition, 1991.

\bibitem{Fr7}
K.~A. Fr{\o}yshov.
\newblock An inequality for the $h$--invariant in instanton {Floer} theory.
\newblock {\em Topology}, 43:407--432, 2004.

\bibitem{Fr10}
K.~A. Fr{\o}yshov.
\newblock {Monopoles over $4$--manifolds containing long necks, I}.
\newblock {\em Geometry \&\ Topology}, 9:1--93, 2005.

\bibitem{Fr12}
K.~A. Fr{\o}yshov.
\newblock A generalized blow-up formula for {Seiberg--Witten} invariants.
\newblock Preprint, math.DG/0604242.

\bibitem{Fr4}
K.~A. Fr{\o}yshov.
\newblock {Monopole Floer homology for rational homology $3$--spheres}.
\newblock In preparation.

\bibitem{Fukaya1}
K.~Fukaya.
\newblock Floer homology for oriented $3$--manifolds.
\newblock In Y.~Matsumoto, editor, {\em Aspects of low dimensional manifolds},
  volume~20 of {\em Adv. Stud. Pure Math.}, pages 1--92. Kinokuniya Company
  Ltd, 1992.

\bibitem{Kelley}
J.~L. Kelley.
\newblock {\em General Topology}.
\newblock Van Nostrand, 1955.

\bibitem{Lang1}
S.~Lang.
\newblock {\em Differential Manifolds}.
\newblock Springer-Verlag, 1985.

\bibitem{Marcolli-Wang1}
M.~Marcolli and B-L. Wang.
\newblock {Equivariant Seiberg--Witten Floer homology}.
\newblock {\em Comm. Anal. and Geom.}, 9:451--639, 2001.

\bibitem{Morgan1}
J.~W. Morgan.
\newblock {\em The {Seiberg--Witten} equations and applications to the topology
  of smooth four-manifolds}.
\newblock Princeton University Press, 1996.

\bibitem{Morgan-Mrowka1}
J.~W. Morgan and T.~S. Mrowka.
\newblock On the gluing theorem for instantons on manifolds containing long
  cylinders, 1994.
\newblock Preprint.

\bibitem{nico1}
L.~I. Nicolaescu.
\newblock {\em Notes on {Seiberg--Witten} theory}.
\newblock American Mathematical Society, 2000.

\bibitem{Safari1}
P.~Safari.
\newblock Gluing {Seiberg--Witten} monopoles.
\newblock math.DG/0311329. To appear in Comm. Anal. and Geom.

\bibitem{das2}
D.~A. Salamon.
\newblock {Spin Geometry and Seiberg--Witten Invariants}.
\newblock Preprint, 1996.

\bibitem{T4}
C.~H. Taubes.
\newblock Self-dual {Yang--Mills} connections on non-self-dual 4-manifolds.
\newblock {\em J.~Differential Geometry}, 17:139--170, 1982.

\bibitem{T5}
C.~H. Taubes.
\newblock Self-dual connections on $4$--manifolds with indefinite intersection
  matrix.
\newblock {\em J.~Differential Geometry}, 19:517--560, 1984.

\bibitem{ShuguangWang}
S.~Wang.
\newblock On orientability and degree of {Fredholm} maps.
\newblock {\em Michigan Math. J.}, 53:419--428, 2005.

\end{thebibliography}
\end{document}